\documentclass[12pt,oneside,reqno]{amsart}
\usepackage{amssymb}
\usepackage{amsmath}
\usepackage{mathtools}
\usepackage{mathrsfs}
\usepackage{amsthm}
\usepackage{amsfonts}
\usepackage{eucal}
\usepackage{dutchcal}
\usepackage[utf8]{inputenc}
\usepackage{comment}
\usepackage{marginnote}
\usepackage[style=alphabetic,backend = bibtex]{biblatex}
\usepackage[dvipsnames]{xcolor}
\usepackage{geometry}
\usepackage{hyperref}
\usepackage{amsaddr}
\usepackage{fancyhdr}
\newtheorem{thm}{Theorem}[section]
\newtheorem{prop}[thm]{Proposition}
\newtheorem{lem}[thm]{Lemma}
\newtheorem{cor}[thm]{Corollary}

\newtheorem{remark}[thm]{Remark}
\title[Estimate and CLT for Diophantine approximation on spheres]{Effective estimate and Central Limit Theorem for Diophantine approximation on spheres}
\author{Zouhair Ouaggag}
\address{University of Zurich, Institute of Mathematics\\ Winterthurerstrasse 190, CH-8057 Zurich, Switzerland}
\email{zouhair.ouaggag@math.uzh.ch}
\date{}
\newtheorem*{thm: first result}{Theorem \ref{thm: first result}}
\newtheorem*{thm: second result}{Theorem \ref{thm: second result}}
\newtheorem*{thm: third result}{Theorem \ref{thm: third result}}
\newtheorem*{thm: fourth result}{Theorem \ref{thm: fourth result}}
\newtheorem*{thm: fifth result}{Theorem \ref{thm: fifth result}}

\newcommand{\SO}{\operatorname{SO}}

\newcommand{\cP}{\mathcal{P}}
\newcommand{\cQ}{\mathcal{Q}}

\newcommand{\cT}{\mathcal{T}}

\newcommand{\cX}{\mathcal{X}}
\newcommand{\cY}{\mathcal{Y}}

\newcommand{\bR}{\mathbb{R}}

\newcommand{\qen}{\enskip \textrm{and} \enskip}
\newcommand{\qand}{\quad \textrm{and} \quad}

\newcommand\subsetsim{\mathrel{%
\ooalign{\raise0.2ex\hbox{$\subset$}\cr\hidewidth\raise-0.8ex\hbox{\scalebox{0.9}{$\sim$}}\hidewidth\cr}}}

\DeclareMathOperator{\supp}{supp}

\begin{document}
\begin{abstract}
We study the counting function of rational approximations with given bounds on the denominator and satisfying the critical Dirichlet exponent on the sphere $S^d$, $d\geq 3$. We give an effective estimate for this counting function, with an error term of square root order, analogous to the optimal estimate in the Euclidean setting. We also show that the counting function has vanishing third and higher correlations and derive a Central Limit Theorem describing its fluctuations. We prove these results using arguments from homogeneous dynamics on the space of orthogonal lattices, in particular effective multiple equidistribution of all orders, which we establish for spherical averages and which could be useful for other applications.
\\

\end{abstract}

\maketitle
\tableofcontents
\newpage

\section{Introduction and main results}

It is well-known in metric Diophantine approximation that for any $c>0$ and Lebesgue-almost all $\alpha \in \mathbb{R}^{m}$, there exist infinitely many solutions $(\textbf{p},q) \in \mathbb{Z}^{m}\times \mathbb{N}$ to the inequality\footnote{$|| \cdot||$ will denote the Euclidean norm.}
\begin{equation}
\label{Diophantine}
 \left\lVert \alpha - \frac{\textbf{p}}{q}\right\rVert~ < ~  \frac{c}{q^{1+\frac{1}{m}}} \quad .
\end{equation}

A refinement of this problem is to count solutions up to a certain bound for the complexity $q$ of the approximants, which leads to consider counting functions such as
$$
\mathrm{N}_{T,c}(\alpha) \coloneqq  \lvert \lbrace (\textbf{p},q) \in \mathbb{Z}^{m}\times \mathbb{N} : 1 \leq q < e^T \text{ and } (\ref{Diophantine}) \text{ holds } \rbrace  \rvert~.
$$

An accurate estimate of the counting function $\mathrm{N}_{T,c}$ was given in the work of W. Schmidt \cite{Schmidt1960AMT}, who proved for more general approximating functions, that for Lebesgue-almost all $\alpha \in [0,1]^m$, 
\begin{equation}
\label{schmidt effective estimate}
\mathrm{N}_{T,c}(\alpha)=\mathsf{C}_{c,m}\,T + O_{\alpha,\varepsilon}\left(  T^{\frac{1}{2}+\varepsilon}\right)~,
\end{equation}
for all $\varepsilon>0$, with a constant $\mathsf{C}_{c,m}>0$ depending only on $c$ and $m$.\\

In recent years, there has been significant interest in the problems of so-called \emph{intrinsic} Diophantine approximation, where one considers approximation by rational points on algebraic varieties, addressing analogues of classical questions in the geometry of numbers and metric Diophantine approximation. Important progress has been achieved in this setting. In particular, Kleinbock and Merrill in \cite{kleinbock2013rational} developed the theory of Diophantine approximation on spheres, which was subsequently generalized to quadratic surfaces with general signatures by Fishman, Kleinbock, Merrill and Simmons in \cite{FKMS2018intrinsic,FKMS2021quadric}. These works established in particular analogues of the classical Dirichlet's and Khinchin's theorems. Then Alam and Ghosh in \cite{alam2020quantitative} proved an asymptotic formula for the number of rational approximants on spheres \eqref{quantitative approximation}. We also mention the works of Ghosh, Gorodnik and Nevo who developed the metric theory of Diophantine approximation on simple algebraic groups, providing estimates for uniform and almost sure Diophantine exponents in \cite{GGN2013Diophantine}, establishing analogues of Khinchin's and Jarnik's theorems in \cite{GGN2014metric}, and deriving an asymptotic formula with an error term for the number of approximants for a range of uniform Diophantine exponents in \cite{GGN2021counting}.\\
 
Let us consider the following intrinsic Diophantine approximation problem. Given $T, c>0$ and $\alpha \in$ $S^{d}$, $d\geq 3$, we consider the inequality \eqref{intrinsic Diophantine intro} (with the critical Dirichlet exponent for intrinsic Diophantine approximation on $S^{d}$)
\begin{equation}
\label{intrinsic Diophantine intro}
 \left\lVert \alpha - \frac{p}{q}\right\rVert~ < ~  \frac{c}{q} \quad ,
\end{equation}
and the counting function given by \eqref{counting function} for intrinsic rational approximations
\begin{equation}\label{counting function}
    \mathsf{N}_{T,c}(\alpha) \coloneqq  \lvert \lbrace (p,q) \in \mathbb{Z}^{d+1}\times \mathbb{N} : \frac{p}{q} \in S^d,~1 \leq q < \cosh T \text{ and } (\ref{intrinsic Diophantine intro}) \text{ holds } \rbrace  \rvert~.
\end{equation}
For this intrinsic Diophantine approximation problem, we have the following important theorems by Kleinbock and Merrill.

\begin{thm}[\cite{kleinbock2013rational}, Theorem 1.1.]
There exists a constant $C \geq 1$ such that for every $\alpha \in S^n$, there exist infinitely many rationals $\frac{p}{q} \in S^d$ such that $$\left\| \alpha - \frac{p}{q} \right\| < \frac{C}{q}.$$
\end{thm}
\begin{thm}[\cite{kleinbock2013rational}, Theorem 1.3.]\label{thm: khintchine type dichotomy KM}
For any $\psi: \mathbb{N} \to (0,\infty)$ non-increasing, the set of well-approximable vectors $\alpha \in S^d$ has full measure (resp. measure zero) if and only if $\sum_{q\geq 1}q^{-1}\psi(q)$ diverges (resp. converges).\\
\end{thm}

Later, Alam and Ghosh proved in \cite{alam2020quantitative} a first quantitative estimate of $\mathsf{N}_{T,c} $, using Birkhoff pointwise ergodicity on the space of orthogonal lattices.

\begin{thm}[\cite{alam2020quantitative}, Theorem 1.2.] There exists a computable constant $C_{c,d}>0$, depending only on $c$ and $d$, such that, for almost every $\alpha \in \text{S}^{d}$, 
\begin{equation}
\label{quantitative approximation}
\lim_{T \to \infty} \frac{\mathsf{N}_{T,c}(\alpha)}{T} = C_{c,d}~.
\end{equation}
\end{thm} 

With a different approach, relying on a crucial result due to Eskin, Margulis and Mozes in \cite{eskin1998upper} about the $L_p$-integrability for $p<2$ of the Siegel transform with respect to the spherical measure, we gave in \cite{ouaggag2022effective} an effective estimate for $\mathsf{N}_{T,c}$, showing that there exists a power saving in the error term, i.e. a constant $\gamma <1$ depending only on the dimension $d$, such that for almost every $\alpha \in \emph{S}^{d}$, 
\begin{equation}
\label{eq:power saving estimate}
\mathsf{N}_{T,c}(\alpha) = C_{c,d}T + O_{\alpha}(T^{\gamma})~.
\end{equation}
 
In order to improve the estimate of the error term in \eqref{eq:power saving estimate} to the order $T^{\frac{1}{2}+\varepsilon}$ as in  \eqref{schmidt effective estimate} for the Euclidean space, we were missing an analog of Roger's formula for the space of orthogonal lattices. A crucial result in this direction was established recently by Kelmer and Yu in \cite{secondmomentKY23}, using spectral theory of spherical Eisenstein series to give an estimate of the second moment of the Siegel transform (see Theorem \ref{thm: second moment siegel tranform}).
 They also derived an effective estimate for $\mathsf{N}_{T,\psi(q)}$, for more general $d$-dimensional quadratic surfaces $\emph{S}$ and approximation function $\psi$. 
 
 \begin{thm}[\cite{secondmomentKY23}, Theorem 1.9.]
     Let $\psi:\mathbb N \rightarrow (0,+\infty)$ decreasing and satisfying $\sum_{q\geq 1}q^{-1}\psi(q)^d=\infty$. For $d\not\equiv 1 (\text{mod } 8)$, for almost every $\alpha \in \emph{S}$, we have
 \begin{equation}
\label{eq:power saving estimate2}
\mathsf{N}_{T,\psi(q)}(\alpha) = C_{d}J_{\psi}(T) + O_{\alpha,\psi}\left(J_{\psi}(T)^{\frac{d+3}{d+4}}\log (J_\psi(T))+I_{\psi}(T)\right)~,
\end{equation}
with $\displaystyle{J_{\psi}(T):=\sum_{1\leq q < T}q^{-1}\psi(q)^d}$, $\displaystyle{I_{\psi}(T)=\sum_{1\leq q < T}q^{-3}\psi(q)^{d+2}}$ and $C_d>0$.
 \end{thm}
 
In order to bound the second moment of $N_T$, we need moreover effective double equidistribution of $\nu$ for spherical averages, for smooth compactly supported test functions. Since the Siegel transform is typically not bounded, we also need to control the contribution in $L_2(\mathcal{Y})$ of smoothing and truncating $\widehat\chi_{_{E}}\circ a_t$ as $t\rightarrow \infty$, using in particular non-escape of mass.

\subsection{First main result: Effective estimate for the counting function $\mathsf{N}_{T,c}$}

Developing our approach from \cite{ouaggag2022effective} and using methods derived from \cite{secondmomentKY23} to analyze the second moment of the Siegel transform, we give in Section \ref{sec:counting} an effective estimate of $\mathsf{N}_{T,c}$ with an error term of square root order, similar to \eqref{schmidt effective estimate}, for all dimensions $d\geq 3$.
 \begin{thm: first result}[first main result]
Let $d\geq 3$. For almost every $\alpha \in \emph{S}^{d}$, for all $\varepsilon>0$, we have 
\begin{equation}
\mathsf{N}_{T,c}(\alpha) = C_{c,d}T + O_{\alpha,\varepsilon}(T^{\frac{1}{2}+\varepsilon})~.
\end{equation}
\end{thm: first result}	

\begin{remark}
    \emph{Some remarks related to Theorem \ref{thm: first result}:
    \begin{enumerate}
        \item The constant $C_{c,d}>0$ in \eqref{quantitative approximation}, \eqref{eq:power saving estimate} and also in Theorems \ref{thm: first result} and \ref{thm: second result}, is equal to the volume of a domain $F_{1,c} \subset \mathbb{R}^{d+2}$ defined explicitly in \eqref{eq:defintion F_1,c}.
        \item  Although the estimate of the error term $T^{\frac{1}{2}+\varepsilon}$ is optimal in the case of the Euclidean space, we cannot conclude about the optimality of this estimate for the sphere. Nevertheless, our analysis of the limit distribution of $\mathsf{N}_{T,c}$ (see Theorem \ref{thm: second result}) suggests that $T^{\frac{1}{2}}$ is the correct normalization, and the error term would be optimal if one could show that the variance $\sigma^2$ is positive.
        \item Our method fails for dimension $d=2$ because of the escape of mass in the space of orthogonal lattices for this dimension (see the estimate \eqref{L2(X) bound} in Proposition \ref{bounds truncated siegel transform}). 
    \end{enumerate}}
\end{remark}

\paragraph{\textbf{Outline of Section \ref{sec:counting}}}
Using the classical Dani correspondence, we first interpret the counting function $\mathsf{N}_{T,c}$ in terms of ergodic averages of a function over a subset of the space of unimodular lattices in $\mathbb{R}^{d+2}$, developing the approach in \cite{kleinbock2013rational}, \cite{alam2020quantitative} and \cite{ouaggag2022effective}. To do so, we embed the sphere $S^d$ in the positive light cone $\mathcal{V}\coloneqq \{ x \in \mathbb{R}^{d+1}\times \mathbb{R}_+ :Q(x)=0\}$ of a quadratic form $Q$ of inertia $(d+1,1)$, and identify good approximations $\frac{p}{q} \in$ $S^d$ for $\alpha \in$ $S^d$ with integer points $(p,q)$ in $\Lambda_0:=\mathbb{Z}^{d+2}\cap \mathcal{V}$ whose images under certain rotations $k_\alpha \in K=\text{SO}(d+1)$ lie in a specific domain $E_{T,c} \subset \mathcal{V}$ (we recall more details about this correspondence in Section \ref{correspondance}). The number of solutions $\mathsf{N}_{T,c}$ is then related to the number of lattice points in the domain $E_{T,c}$, which can be approximated by a more convenient domain $F_{T,c}$ and tessellated by the action of a hyperbolic subgroup $\{a_t, t\in \mathbb{R} \} \subset $ SO$(Q)$. \\
For a bounded and compactly supported function $f$ on $\mathcal{V}$, we introduce the \textit{light-cone Siegel transform} $\widehat f$, defined for any lattice $\Lambda \subset \mathcal{V}$ by
$$\widehat f(\Lambda):= \sum_{x\in \Lambda\setminus \{0\}}f(x) .$$ 
The counting function $\mathsf{N}_{T,c}$ can then be related to the light-cone Siegel transform of the characteristic function $\chi$ of an elementary domain $F_{1,c}\subset \mathcal{V}$, using averages of the form
\begin{equation}\label{approximation N_T,c}
    \mathsf{N}_{T,c}(\alpha) \approx \sum_{t=0}^{T-1} \widehat\chi \circ a_t (k_\alpha \Lambda_0),
\end{equation}
i.e. ergodic averages of the light-cone Siegel transform $\widehat{\chi}$ along $K$-orbits pushed by $\{a_t\}$. The analysis of these averages can then be carried out using dynamics on the space of orthogonal lattices $\mathcal{X}\cong \SO(Q)/\SO_\mathbb{Z}(Q)$.\vspace{7pt}\\
\indent In \cite{alam2020quantitative}, the authors use ergodicity of the $\{a_t\}$-action and Birkoff's Ergodic Theorem to obtain an asymptotic estimate for these ergodic averages. In order to obtain an estimate with an error term, we use effective pointwise equidistribution along $\{a_t\}$-orbits, which we derive from effective double equidistribution of $K$-orbits pushed by $\{a_t\}$ (see the estimate on correlations of all orders in Proposition \ref{double equidistribution K orbits}), and then derive an almost-everywhere bound from an L$^p$-bound on ergodic averages with $p>1$ (Proposition \ref{pointwise equidistribution}). However, using effective equidistribution requires considering smooth and compactly supported test functions, whereas our test functions have typically none of these regularities. We address this issue in two steps.\vspace{7pt}\\ 
\indent We first introduce and study the integrability of the Siegel transform on the space of orthogonal lattices, then show (Proposition \ref{bounds truncated siegel transform}) that the Siegel transform $\widehat{f}$ of a compactly supported function $f$ can be approximated by a truncated Siegel transform $\hat{f}^{(L)}$ in a way to control the approximation on translated $K$-orbits, i.e. control $| \widehat{f}\circ a_t - \widehat{f}^{(L)}\circ a_t|$ with respect to the probability measure on the orbits. In a second step (Proposition \ref{smooth approximation}), we show that the characteristic function $\chi$ of the elementary domain $F_{1,c}$ can be approximated by a family of smooth compactly supported functions $f_{\varepsilon}$, again in a way to keep control of the approximated Siegel transform $\widehat{f}_{\varepsilon}$ over translated $K$-orbits. In this process, we also need non-divergence results for the Siegel transform with respect to the probability measure on $K$-orbits (Propositions \ref{alpha_Lp} and \ref{non-espace of mass}).\\ 

Thus, we can relate the counting function $ \widehat{\chi}$ to a smooth and compactly supported function on the space of orthogonal lattices and then establish an effective almost-everywhere estimate for its ergodic averages (Theorem \ref{pointwise equidistribution counting function}). A crucial ingredient needed to obtain an "almost-everywhere" error term of square root order is an estimate of the variance of $\mathsf{N}_{T,c}$, which we obtain using recent results from Kelmer and Yu in \cite{secondmomentKY23} on the second moment of the Siegel transform on the light cone.\\

\subsection{Second main result: Limit distribution of the counting function $\mathsf{N}_{T,c}$}
Another problem related to the intrinsic Diophantine approximation \eqref{intrinsic Diophantine} is to understand the limit distribution of $\mathsf{N}_{T,c} $ as a random variable on the sphere. If we consider counting functions of approximations with denominators $q$ in ranges $[\cosh (t) , \cosh (t+1) ]$, i.e. functions
$$
\mathcal{N}_{t,c} := \mathsf{N}_{t+1,c}- \mathsf{N}_{t,c},
$$
and if the random variables $\mathcal{N}_{t_1,c}$ and $\mathcal{N}_{t_2,c}$ "decorrelate" for large $t_1$, $t_2$ and $|t_2-t_1|$, then Theorem \ref{thm: first result} would follow by a Law of large numbers for the random variable $\mathsf{N}_{T,c}=\sum_{t=0}^{T-1}\mathcal{N}_{t,c}$. This heuristic was developed by Bjöklund and Gorodnik in \cite{bjoerklund2018central} to show that the counting function $\mathrm{N}_{T,c}$ for the Euclidean space follows a Central Limit Theorem, using higher-order mixing in the space of unimodular lattices as the dynamical translation of quasi-independent random variables. Using a similar approach as in \cite{bjoerklund2018central} and the recent results by Kelmer and Yu in \cite{secondmomentKY23} about the second moment of the Siegel transform, we show in Section \ref{sec:clt counting} that the counting function $\mathsf{N}_{T,c}$ on the sphere also follows a Central limit theorem.
	
\begin{thm: second result}[second main result] 
Let $d\geq 3$. Then, for every $\xi \in \mathbb{R}$, 
\begin{equation}
\mu_{_{S^d}}\left( \left\lbrace \alpha \in S^d:\frac{\mathsf{N}_{T,c}(\alpha)-C_{c,d}\cdot T}{T^{1/2}}<\xi \right\rbrace\right) \rightarrow \emph{Norm}_\sigma(\xi)~,\quad \text{as }T\rightarrow \infty,
\end{equation}
where $\emph{Norm}_\sigma$ denotes the normal distribution with variance $\sigma^2\geq 0$.
\end{thm: second result}

\begin{remark}
\emph{The variance $\sigma^2$ in Theorem \ref{thm: second result} is given explicitly in \eqref{eq: explicite variance}.}
\end{remark}

\paragraph{\textbf{Outline of Section \ref{sec:clt counting}}}

\indent In order to show that $\mathsf{N}_{T,c} $ follows a Central Limit Theorem, we will use the method of cumulants (Section \ref{sec:cumulants}), which is equivalent to the more widely known method of moments. The normal distribution is characterized by vanishing cumulants of orders $r\geq 3$, which can be expressed in the dynamical language as higher-order correlations of an averaging function of the form
    $$
        \mathsf{F}_N \approx \frac{1}{N^{1/2}}\sum_{t=0}^{N-1} \widehat\chi \circ a_t\;.
    $$
The "quasi-independence" of sampling observables $\widehat\chi \circ a_t$ along $K$-orbits pushed by $\{a_t,t\geq 0\}$ corresponds to multiple equidistribution of these orbits in the space of orthogonal lattices, which we establish in Section \ref{sec:higher order correlations} (Proposition \ref{double equidistribution K orbits}). As explained in the outline of Section \ref{sec:counting} earlier, using effective equidistribution requires considering smooth and compactly supported test functions. We address this issue by considering a family of smooth approximations $f_\varepsilon$ of the test function $\chi$ and truncated Siegel transforms $\widehat f_\varepsilon^{(L)}$, with parameters $\varepsilon, L>0$, and keep track of the approximation error for averages along translated $K$-orbits (Propositions \ref{bounds truncated siegel transform} and \ref{smooth approximation}). We also need to show that these approximations still give the same limit distribution as the original averaging function $\mathsf{F}_N$ (Section \ref{approximations of F_N}). \vspace{7pt}\\
\indent In Section \ref{sec:cumulants}, we use exponential multiple equidistribution established previously to show that the cumulants of $\mathsf{F}_N$ of orders $r\geq 3$ vanish as $N\rightarrow \infty$. To do so, we follow the argument developed in  \cite{bjoerklund2020} and \cite{bjoerklund2018central}, analyzing the joint cumulants through a decomposition into sub-sums of correlations corresponding to ``separated" or ``clustered" tuples $t_1,\dots,t_r$ and controlling their size in terms of the parameters $\varepsilon,L>0$ related to the smoothing and truncation of the test function.\\

The limit variance obtained in Section \ref{sec:counting}, together with vanishing cumulants of orders $r\geq 3$ from Section \ref{sec:cumulants} complete the characterization of the normal distribution for $\mathsf{F}_N$. In Section \ref{sec:proof of the first main theorem} we relate the distribution of $\mathsf{F}_N$ to the distribution of our counting function $\mathsf{N}_{T,c}$ and conclude the proof of the CLT-Theorem.\\

\paragraph{\textbf{Acknowledgment}} I am grateful to Dubi Kelmer for sharing his ideas and for the useful discussions. Many thanks also to René Pfitscher for our discussions and his useful observations. I am very grateful to Alex Gorodnik for his guidance throughout this research work.

\newpage

\section{Effective estimate for counting Diophantine approximations on spheres}\label{sec:counting}

In this section, we study the counting function of rational approximations on the sphere $S^d$ defined, for $T>0$, $c>0$, $\alpha \in S^{d}$ and $d\geq 3$, by
\begin{equation*}
    \mathsf{N}_{T,c}(\alpha) \coloneqq  \lvert \lbrace (p,q) \in \mathbb{Z}^{d+1}\times \mathbb{N} : \frac{p}{q} \in S^d,~1 \leq q < \cosh T \text{ and } (\ref{intrinsic Diophantine}) \text{ holds } \rbrace  \rvert~,
\end{equation*}
for the intrinsic Diophantine approximation problem (with the critical Dirichlet exponent for intrinsic Diophantine approximation on $S^{d}$) given by the inequality
\begin{equation}
\label{intrinsic Diophantine}
 \left\lVert \alpha - \frac{p}{q}\right\rVert~ < ~  \frac{c}{q} \;.
\end{equation}
We prove an effective estimate of the counting function $N_{T,c}$ with an error term of square root order.

\begin{thm}[]
    \label{thm: first result}
Let $d\geq 3$. For almost every $\alpha \in \emph{S}^{d}$, for all $\varepsilon>0$, we have 
\begin{equation}
\mathsf{N}_{T,c}(\alpha) = C_{c,d}T + O_{\alpha,\varepsilon}(T^{\frac{1}{2}+\varepsilon})~.
\end{equation}
\end{thm}

Our approach is to use an analog of the classical Dani correspondence for the space of orthogonal lattices, and show that $\mathsf{N}_{T,c}$ can be related to the Siegel transform of the characteristic function $\chi$ of an elementary domain on the light cone, using averages of the form
\begin{equation*}
    \mathsf{N}_{T,c}(\alpha) \approx \sum_{t=0}^{T-1} \widehat\chi \circ a_t (k_\alpha \Lambda_0),
\end{equation*}
i.e. ergodic averages of the Siegel transform $\widehat{\chi}$ along compact orbits $k_\alpha \Lambda_\alpha$ pushed by the hyperbolic subgroup $\{a_t\}$. The analysis of these averages can then be carried out using dynamics on the space of orthogonal lattices, as detailed in the following section.
 
\subsection{Diophantine approximation and dynamics on the space of lattices}
\label{correspondance}

We present in this section the particular formulation of Dani's correspondence relating intrinsic Diophantine approximation on the sphere $S^{d}$ and the dynamics of orthogonal lattices in $\mathbb{R}^{d+2} $(see also \cite{kleinbock2013rational}, \cite{alam2020quantitative} and \cite{ouaggag2022effective}).\\

We consider the quadratic form $Q: \mathbb{R}^{d+2} \rightarrow \mathbb{R}$ defined by
\begin{equation*}
Q(x) := \sum_{i=1}^{d+1} x_i^{2} - x_{d+2}^{2}, \quad \text{ for } x= (x_1, \dots , x_{d+2})~,
\end{equation*}
and the embedding of $S^{d}$ in the positive light cone
$$
\mathcal{V} := \{  x \in \mathbb{R}^{d+2} : Q(x)=0,~ x_{d+2} > 0 \}~,
$$
via $\alpha \mapsto (\alpha,1)$, which yields a one-to-one correspondence between primitive integer points on the positive light cone, $(p,q) \in  \mathcal{V} \cap \mathbb{Z}^{d+2}_{\text{prim}}$, and rational points on the sphere, $ \frac{p}{q} \in S^{d}$.\\

We denote by $G$=SO$(Q)^\circ \cong$ SO$(d+1,1)^\circ$ the connected component of the group  of orientation-preserving linear transformations which preserve $Q$. We denote by $\Lambda_0 := \mathcal{V} \cap \mathbb{Z}^{d+2}$ the set of integer points on the positive light cone. By a \emph{lattice $\Lambda$ in $\mathcal{V}$} we mean a set of the form $g\Lambda_0$ for some $g\in G$. 
If we denote by $\Gamma$ the stabilizer of $\Lambda_0$ in $G$, then $\Gamma$ is a lattice in $G$ containing the subgroup SO$(Q)_{\mathbb{Z}}^\circ$ of integer points in $G$, as a finite index subgroup.
The space of lattices in $\mathcal{V}$ can be identified with the homogeneous space $\mathcal{X} := G/\Gamma$, endowed with the $G$-invariant probability measure $\mu_{\mathcal{X}}$.\\

We introduce, for a bounded and compactly supported function $f$ on $\mathcal{V}$, the \textit{light-cone Siegel transform} of $f$, denoted by $\widehat f$ and defined for any lattice $\Lambda \subset \mathcal{V}$ by
\begin{equation}\label{def Siegel Transform}
    \widehat f(\Lambda):= \sum_{x\in \Lambda\setminus \{0\}}f(x) .
\end{equation}

Let $K$ denote the subgroup of $G$ that preserves the last coordinate in $\mathbb{R}^{d+2} $, i.e.
$$
K= \begin{pmatrix}
\text{SO}(d+1) &  \\
 & 1
\end{pmatrix} \cong \text{SO}(d+1)~,
$$
equipped with the Haar probability measure $\mu_K$.\\

The sphere $S^d$ can be realized as a quotient of $K$, endowed with a unique left $K$-invariant probability measure, giving a natural correspondence between full-measure sets in $K$ and those in $S^d$. \\

For $k \in K$ we define $\alpha_k \in  $ $S^d$ by \begin{equation}\label{eq:definition alpha_k}
    k(\alpha_k,1) = (0,\dots,0,1,1) \in \mathcal{V}.
\end{equation}
For $(p,q) \in \Lambda_0$, we write $$k(p,q)= (\xi_1,\xi_2, \dots , \xi_{d+1},\xi_{d+2}) \in \mathcal{V},\quad \text{with }\xi_{d+2}=q $$ and observe the following correspondence\footnote{$\lVert \cdot \rVert$ denotes the Euclidean norm.} (\cite{alam2020quantitative}, Lemma 2.2.):
\begin{align*}
 \left\lVert \alpha_k - \frac{p}{q}\right\rVert < \frac{c}{q} \quad &\Leftrightarrow \quad \lVert q( \alpha_k,1) -(p,q) \rVert < c \\
&\Leftrightarrow \quad \lVert qk(\alpha_k,1) - k(p,q)\rVert < c, \\
&\Leftrightarrow \quad \lVert(\xi_1,\xi_2, \dots , \xi_d, \xi_{d+1}-\xi_{d+2},0)\rVert < c, \\
&\Leftrightarrow \quad \xi_1^2+\dots +\xi_{d}^2+\left( \xi_{d+1}-\xi_{d+2}\right)^2 < c^2, \\
&\Leftrightarrow \quad 2\xi_{d+2}(\xi_{d+2}-\xi_{d+1}) < c^{2} \quad (\text{since }\xi_1^2+\dots +\xi_{d+1}^2=\xi_{d+2}^2). 
\end{align*}
Hence, if we write 
$$E_{T,c}\coloneqq  \{ x \in \mathcal{V} : 2x_{d+2}(x_{d+2}-x_{d+1}) <c^{2}, 1\leq x_{d+2} < \cosh T \},
$$
then we have
\begin{equation}
\label{N_T,c as volume of E_T,c}
\mathsf{N}_{T,c}(\alpha_k) =\lvert E_{T,c} \cap k\Lambda_0 \lvert.
\end{equation}

We write $
\mathcal{Y} \coloneqq K \Lambda_0$, equipped with the Haar probability measure $\mu_{\mathcal{Y}}$.\\

We also consider elements 
$$
a_t =  \begin{pmatrix}
I_d & & \\
 &  \cosh t & -\sinh t \\
 & -\sinh t & \cosh t \end{pmatrix} \quad \in G
$$
and the corresponding one-parameter subgroup
$$
A = \left\lbrace  a_t : t \in \mathbb{R}\right\rbrace
$$
endowed with the natural measure $dt$.\\

In order to use the dynamics of translates of $\mathcal{Y}$ for the Diophantine approximation problem (\ref{N_T,c as volume of E_T,c}), we first approximate $E_{T,c}$ by a domain offering a convenient tessellation under the action of the subgroup $A$. We start with a similar approach as in \cite{alam2020quantitative} and improve the approximation by $F_{T,c}$ in order to satisfy the accuracy obtained later for the counting function.

\subsubsection{Tesselation of $E_{T,c}$}
For $T>0$, we consider the following domain on the light-cone
\begin{equation}\label{eq:defintion F_1,c}
    F_{T,c} := \big\{ x \in \mathcal{V} : x_{d+2}^{2}- x_{d+1}^{2} < c^{2}~,~ c\leq x_{d+2}+ x_{d+1} < ce^{T} \big\},
\end{equation}
and the sequence of domains $\left(F_{T,c,t}\right)_{t\geq 0}$ defined by 
\begin{equation}\label{eq:defintion F_1,c_l}
    F_{T,c,t} := \big\{ x \in \mathcal{V} : x_{d+2}^{2}- x_{d+1}^{2} < c_t^{2}~,~ c\leq x_{d+2}+ x_{d+1} < ce^T \big\},
\end{equation}
with $c_t\coloneqq c(t)= c\cdot \left( \frac{t+1}{t+2}\right)^{1/2}$.\\

We first observe that the transformation $a_t\in A$ preserves the difference $ x_{d+2}^2-x_{d+2}^2$ and scales the sum $x_{d+2}+x_{d+1}$ by a factor $e_t$, which gives the following \emph{tesselation}, for all integers $N \geq 1$,
\begin{align}
    &F_{N,c}= \bigsqcup_{t=0}^{N-1} a_{-t}F_{1,c}.\label{tessaltion of F_T,c}
    \end{align}
    
We show that $E_{N,c}$ can be approximated by sets of the form $F_{N,c}$ and $F_{N,c,N-1}$ up to some bounded sets in $\mathcal{V}$. 
\begin{lem} There exist positive integers $N_0$ and $r_0$, depending only on $c$, and bounded sets $B_0,B_0'\subset \mathcal{V}$ such that, for all $N\geq N_0$, we have
\begin{equation}
\label{sandwiching}
F_{N-r_0,c, N-r_0-1} \setminus B_{0}' \quad \subseteq \quad
E_{N,c} \setminus B_0 \quad \subseteq \quad F_{N+r_0,c}~.
\end{equation} 
\end{lem}

\begin{proof}
Consider the sets
\begin{align*}
    &B_0 \coloneqq \Big\{ x \in \mathcal{V}: x_{d+2} \leq \frac{c^2+c}{2}+1\Big\},\\
    \text{and }\quad &C_t:=a_{-t}F_{1,c,t}\cap \left(\Big\{ x \in \mathcal{V}: \frac{|x_{d+1}|}{x_{d+2}}\leq \frac{t+1}{t+2} \Big\}\cup B_0\right), \quad \text{for all } t\geq 0.
\end{align*}

For $x\in E_{N,c} \setminus B_0 ~$, we have
$$2x_{d+2}(x_{d+2}-x_{d+1})<c^2 ~\text{ and, since }~ 1\leq x_{d+2}
 ~\text{, this implies }~ x_{d+2}-x_{d+1}< c^2~,$$
thus 
\begin{equation}
\label{eq:condition 1}
    x_{d+2}+x_{d+1} ~>~ 2x_{d+2}-c^2 ~>~2\left(\frac{c^2+c}{2}+1\right)-c^2 ~>~c~.
\end{equation}
Further, 
\begin{equation}
\label{eq:condition 2}
    x_{d+1}\leq x_{d+2} ~\text{ implies }~ x_{d+2}+x_{d+1}\leq 2x_{d+2}<2\cosh N < c e^{N+r_0},
\end{equation}
for all $N\geq N_0$, for some integers $N_0\geq1$ and $r_0\geq 1$ both depending only on $c$.\\ \\
We also have
\begin{equation}
\label{eq:condition 3}
    x_{d+2}^2-x_{d+1}^2= (x_{d+2}+x_{d+1})(x_{d+2}-x_{d+1}) \leq 2x_{d+2}(x_{d+2}-x_{d+1}) < c^2.
\end{equation}
The inequalities \eqref{eq:condition 1}, \eqref{eq:condition 2} and \eqref{eq:condition 3} prove the second inclusion in \eqref{sandwiching}.\\ \\
To show the first inclusion in \eqref{sandwiching}, we first show that, for all $N\geq N_0$, we have
\begin{equation}\label{first inclusion for all T}
    \bigsqcup_{t=0}^{N-r_0-1} \left(F_{1,c,t}\setminus C_{t}\right)\; \subseteq \;E_{N,c}\setminus B_0.
\end{equation}
For $\displaystyle x\in  \bigsqcup_{t=0}^{N-r_0-1} \left(F_{1,c,t}\setminus C_{t}\right) $, we have, for some $t\in \{0,\dots, N-r_0-1\}$, 
$$x_{d+2}^2-x_{d+1}^2<c_{t}^2<c^2 ~\text{ and, since }~ x_{d+2}+x_{d+1} \geq c
 ~\text{, this implies }~ x_{d+2}-x_{d+1}< c~,$$
hence 
\begin{equation}
\label{eq:condition 4}
    2x_{d+2}~<~x_{d+2}+x_{d+1}+c ~<~ ce^{N-r_0}+c~\leq ~2\cosh N, 
\end{equation}
for all $N\geq N_0$, for some integers $N_0\geq 1$ and $r_0\geq 1$ depending on $c$ and large enough.\\ 
We also have, since $x\notin B_0$, 
\begin{align}
    &x_{d+2}\geq 1 \label{lower bound x_d+2}\\
    \qand &x_{d+1} > x_{d+2}-c > \frac{c^2+c}{2}+1 -c > 0~
,\nonumber
\end{align}
which implies,
\begin{align}
2x_{d+2}(x_{d+2}-x_{d+1}) &= (x_{d+2}^2-x_{d+1}^2)\frac{2x_{d+2}}{x_{d+2}+x_{d+1}}\nonumber \\
&< c_t^2\frac{2x_{d+2}}{x_{d+2}+x_{d+1}} \nonumber\\
&= c^2\left(\frac{t+1}{t+2}\right)\frac{2x_{d+2}}{x_{d+2}+x_{d+1}} \nonumber\\
&< c^2\frac{x_{d+2}+|x_{d+1}|}{x_{d+2}+x_{d+1}}\nonumber\\
&= c^2~.
\label{eq:condition 5}
\end{align}
The inequalities \eqref{eq:condition 4},\eqref{lower bound x_d+2} and \eqref{eq:condition 5} prove inclusion \eqref{first inclusion for all T}.\\

Moreover, for $t\geq N_0$ large enough and $x\in a_{-t}F_{1,c,t}$, we have
\begin{equation}
    ce^t\leq x_{d+2}+x_{d+1} \leq 2 x_{d+2}, \label{2. bound x_{d+2}},
\end{equation} 
hence $a_{-t}F_{1,c,t}\cap B_0 =\emptyset$, which implies that $C_t = a_{-t}F_{1,c,t}\cap \left\{\frac{|x_{d+1}|}{x_{d+2}}\leq \frac{t+1}{t+2}\right\}$. We have then, for all $x\in C_t$, 
\begin{align}
&x_{d+2}^2~ <~ c_t^2+x_{d+1}^2 ~\leq~ c_t^2+\left( \frac{t+1}{t+2}\right)^2x_{d+2}^2 ~,\nonumber\\
\text{thus}\qquad &x_{d+2} < \frac{t+2}{(2t+3)^{1/2}}c_t~\ll_c t^{1/2}. \label{1. bound on x_{d+2}}
\end{align}
From \eqref{2. bound x_{d+2}} and \eqref{1. bound on x_{d+2}}, we have that $C_t=\emptyset$ for all $t\geq N_0$ large enough, and therefore
$$ \bigsqcup_{t=0}^{N-r_0-1} \left(F_{1,c,t}\setminus C_{t}\right)\; = \bigsqcup_{t=0}^{N-r_0-1} F_{1,c,t}\setminus B_{0}',$$
with $\displaystyle B_0':= \bigsqcup_{t=0}^{N_0-1}C_t$.
\end{proof}

It follows from (\ref{N_T,c as volume of E_T,c}) and (\ref{sandwiching}) that, for $T\geq T_0>0$ large enough,
\begin{equation}
\label{sandwiching N_T,c}
\Big|\bigsqcup_{t=0}^{\left\lfloor T\right\rfloor-r_0-1} F_{1,c,t}\cap k\Lambda_0\Big|+O(1) ~\leq ~\mathsf{N}_{T,c}(\alpha_k) + O(1) ~\leq~ \Big|F_{\left\lceil T \right\rceil+r_0,c} \cap k\Lambda_0\Big|. \\
\end{equation}

We will need later an estimate of the counting error in \eqref{sandwiching N_T,c} related to the approximation \eqref{sandwiching}, in terms of the parameter $t>0 $. This is given in the following lemma. 

\begin{lem}\label{lem: estimate errors sandwishing}
We have, for all $t\geq 1$,
\begin{align}
\label{volum F_1,c approximation}
&\emph{vol}(F_{1,c}) =\emph{vol}( F_{1,c,t})+ O_{c,d}\left(t^{-1} \right).
\end{align}
\end{lem}
\begin{proof} For all $t\geq 0$, we have $c_t<c$, hence $F_{1,c,t}\subset F_{1,c}$ and
\begin{align*}
     F_{1,c}\setminus F_{1,c,t}&= \{ x \in \mathcal{V} : c_t^2\leq x_{1}^{2}+\dots+ x_{d}^{2}< c^2~,~ c\leq x_{d+2}+ x_{d+1} < ce \}
\end{align*}
The invariant volume in coordinates $x_1,\ldots,x_{d+1}$ is given by $x_{d+2}^{-1}dx_1\ldots dx_{d+1}$.
We use the new coordinates 
$$
y_1=x_1,\;\ldots,\; y_d=x_d,\; y_{d+1}=x_{d+2}+x_{d+1},
$$
We have $(y_{d+1}-x_{d+1})^2 =\sum_{i=1}^d y_i^2 + x_{d+1}^2$, so that 
$$
x_{d+1} = -\frac{1}{2}y_{d+1}^{-1}\left(\sum_{i=1}^d y_i^2-y^2_{d+1}\right).
$$
Also,
$$
x_{d+2} = y_{d+1}-x_{d+1} = \frac{1}{2}y_{d+1}^{-1}\left(\sum_{i=1}^d y_i^2+y^2_{d+1}\right).
$$
Hence, the invariant volume in the new coordinates is given by
$$
2y_{d+1}\left(\sum_{i=1}^d y_i^2+y^2_{d+1}\right)^{-1} \cdot \left|\frac{\partial x_{d+1}}{\partial y_{d+1}}\right|  dy_1\ldots dy_{d+1}=y_{d+1}^{-1} dy_1\ldots dy_{d+1}.
$$
The domain of integration is given by the inequalities 
$$
c_t^2\leq \sum_{i=1}^d y_i^2<c^2\qand c\leq y_{d+1}<ce.
$$
We denote by $\omega_d$ the volume of the $d$-dimensional unit ball, and compute 
\begin{align*}
\mu_\mathcal{V}(F_{1,c}\setminus F_{1,c,t}) &= \int_c^{ce} y_{d+1}^{-1} \left( \int_{c_t^2\leq \sum_{i=1}^d y_i^2 < c^2} dy_1 \ldots dy_d \right) dy_{d+1}\\
&=
\int_c^{ce} y_{d+1}^{-1} \omega_d (c_t^d-c^d) dy_{d+1}\\
&=\omega_d (c_t^d-c^d) = O_d\left( c-c_t \right).
\end{align*}
Further, we have for $t\geq 1$,
\begin{align*}
&c-c_{t}=c\left(1-\frac{t^{1/2}}{(t+1)^{1/2}}\right)=O_c(t^{-1})~,
\end{align*}
hence the claim.
\end{proof}

\subsubsection{The counting function $\mathsf{N}_{T,c}$ as ergodic averages on the light-cone}\label{subsec:The counting function as ergodic averages on the light-cone}

We denote by $\chi_{1,c}$ the characteristic function of $F_{1,c}$. By definition \eqref{def Siegel Transform} the light-cone Siegel transform of $\chi_{1,c}$ is given by
$$\widehat{\chi}_{1,c}(\Lambda) \coloneqq \sum_{z \in \Lambda \setminus \{0\}} \chi_{1,c} (z),~ \text{ for all } \Lambda \in \mathcal{X}.$$

The tessellation \eqref{tessaltion of F_T,c} implies, for all $N\geq 1$ and all $\Lambda \in \mathcal{X}$,
\begin{align}    
 &|F_{N,c} \cap \Lambda| ~=~ \sum_{t=0}^{N-1} \widehat{\chi}_{1,c}(a_t\Lambda). \label{sandwiching F_T,c}
\end{align}

It follows from (\ref{sandwiching N_T,c}) and (\ref{sandwiching F_T,c}), for $T>T_0>0$ large enough and $k \in K$ as in (\ref{N_T,c as volume of E_T,c}),
\begin{align}
\label{resandwiching N_T,c}
\sum_{t=0}^{\left\lfloor T\right\rfloor-r_0 -1} \widehat{\chi}_{1,c,t}(a_t k\Lambda_0) +O\left( 1\right) ~\leq ~\mathsf{N}_{T,c}(\alpha_k) +O(1) ~\leq~ \sum_{t=0}^{\left\lceil T\right\rceil+r_0 -1} \widehat{\chi}_{1,c}(a_t k\Lambda_0).
\end{align}

The estimate \eqref{resandwiching N_T,c} shows that estimating $\mathsf{N}_{T,c}(\alpha)$ amounts to analyzing ergodic sums of the form \mbox{$\sum_{t=0}^{N} \widehat{\chi}_{1,c}\circ a_t$} on $\mathcal{Y}= K\Lambda_0$. We will use for this purpose effective higher order equidistribution results for unimodular lattices, specialized to $\mathcal{Y}$, which we discuss in the following section.
	
\subsection{Estimates on higher order correlations}
\label{sec:higher order correlations}

We prove in this section an effective equidistribution of $K$-orbits by relating it to effective equidistribution of unstable horospherical orbits established in a more general setting in \cite{bjoerklund2021equidistribution}. We recall the notations
\begin{align*}
G &= \text{SO}(Q)^\circ \cong \text{SO}(d+1,1)^\circ,\\
K &= \begin{pmatrix}
\text{SO}(d+1) &  \\
 & 1
\end{pmatrix},\\
a_t &=  \begin{pmatrix}
I_d & & \\
 &  \cosh t & -\sinh t \\
 & -\sinh t & \cosh t \end{pmatrix} ~ \in G, \quad \text{and }~ A = \left\lbrace  a_t : t \in \mathbb{R}\right\rbrace.
\end{align*}

We also consider the corresponding horospherical subgroups
\begin{align*}
&U = \left\lbrace  g\in G ~ : ~ a_{-t}ga_{t}\rightarrow e \text{ as }t\rightarrow \infty \right\rbrace,\\
&U^- = \left\lbrace  g\in G ~ : ~ a_{t}ga_{-t}\rightarrow e \text{ as }t\rightarrow \infty \right\rbrace,\\
&H = \left\lbrace  g\in G ~ : ~ a_{t}g=ga_{t}\right\rbrace,
\end{align*}
and the Haar measures $d\mu_K$, $dt$ and $d\mu_U$ on $K$, $A$ and $U$ respectively.

It will be important in our argument later that the error term in the effective equidistribution is explicit in terms of the $C^l$-norm, for some $l\geq 1$, of the test functions on $\mathcal{X}$. We introduce below the required notations.\\

Every $Y \in$ Lie$(G)$ defines a first-order differential operator $D_Y$ on $C_c^{\infty}(\mathcal{X})$ by
$$
D_Y(\phi)(x) \coloneqq \frac{d}{dt}\phi(\exp (tY)x)|_{t=0}.
$$

If $\{ Y_1,\dots, Y_r\}$ is a basis of Lie$(G)$, then every monomial $Z=Y_1^{l_1}\dots Y_f^{l_r}$ defines a differential operator by 
\begin{equation}
\label{diff operator}
D_Z\coloneqq D_{Y_1}^{l_1}\dots D_{Y_r}^{l_r},
\end{equation}
of degree deg$(Z)=l_1+\dots+l_r$. For integers $l\geq 0$ and $\phi \in C_c^{\infty}(\mathcal{X})$, we write
\begin{equation}\label{def C_l norm}
   ||\phi||_l \coloneqq ||\phi||_{C^l}= \sum_{\text{deg}(Z)\leq l} ||D_Z(\phi)||_{\infty}.
\end{equation}

A crucial ingredient for our analysis is the following effective equidistribution result for higher order correlations on translated $U$-orbits (Theorem \ref{double equidistribution U orbits}) and the analogous result we derive for translated $K$-orbits (Theorem \ref{double equidistribution K orbits}).

\begin{thm}[specializes Theorem 1.2. in \cite{bjoerklund2021equidistribution}]
\label{double equidistribution U orbits}
For every $r\geq 1$ there exist $\gamma_r>0$ and $l_r\geq 1$ such that, for every $f \in C_{c}^{\infty}(U)$ and $\varphi_1,\dots \varphi_r \in C_{c}^{\infty}(\mathcal{X})$ and every compact subset $L \subset \mathcal{X}$, there exists $C>0$ such that for every $\Lambda \in L$ and $t_1,\dots t_r > 0$, we have
\begin{align*}
    &\left|  \int_{U} f(u)\left(\prod_{i=1}^{r}\varphi_i (a_{t_i} u\Lambda) \right) d\mu_U(u)- \left(\int_{U}fd\mu_U\right) \left(\prod_{i=1}^{r}\int_{\mathcal{X}} \varphi_i d\mu_\mathcal{X}\right) \right|\\
    &\qquad\qquad\qquad\qquad\qquad\qquad\qquad\qquad\qquad\leq C e^{-\gamma_r D(t_1,\dots,t_r)}||f||_{l_r}\prod_{i=1}^{r}||\varphi_i||_{l_r}~,
\end{align*}
 
where $D(t_1,\dots,t_r)\coloneqq \min\{t_i,|t_i-t_j|:1\leq i\neq j\leq r\}.$
\end{thm}

\begin{thm}
\label{double equidistribution K orbits}
For every $r\geq 1$ there exist $\delta_r>0$ and $l_r\geq 1$ such that, for every $f \in C^{\infty}(K)$ and $\varphi_1,\dots \varphi_r \in C_{c}^{\infty}(\mathcal{X})$ and every compact subset $L \subset \mathcal{X}$, there exists $C>0$ such that for every $\Lambda \in L$ and $t_1,\dots t_r > 0$, we have
\begin{align*}
    &\left|  I_{\Lambda,f,\varphi_1,\dots,\varphi_r}(t_1,\dots,t_r) - \left(\int_{K}fd\mu_K\right) \left(\prod_{i=1}^{r}\int_{\mathcal{X}} \varphi_i d\mu_\mathcal{X}\right) \right|\\
    &\qquad\qquad\qquad\qquad\qquad\qquad\qquad\qquad\qquad
 \leq C e^{-\delta_r D(t_1,\dots,t_r)}||f||_{l_r}\prod_{i=1}^{r}||\varphi_i||_{l_r}~,\end{align*}
where $~I_{\Lambda,f,\varphi_1,\dots,\varphi_r}(t_1,\dots,t_r) \coloneqq \int_{K} f(k)\left(\displaystyle\prod_{i=1}^r\varphi_i (a_{t_i} k\Lambda)\right) d\mu_K(k)$.
\end{thm}

\begin{proof}
We consider the centralizer of $A$ in $K$,  
\begin{equation*}
 M\coloneqq \text{cent}_K(A)=K\cap H = \begin{pmatrix}
\text{SO}(d) &  \\
 & I_2
\end{pmatrix} \cong \text{SO}(d),
\end{equation*}
and the submanifold $S\subset K$ defined via the exponential map by
$$
\text{Lie}(S)= \left\lbrace \begin{pmatrix}
0_d &\textbf{s}& &  \\
-\textbf{s}^T &0 &\\
& & 0
\end{pmatrix} : \textbf{s}\in \mathbb{R}^d\right\rbrace.$$
We have $\text{Lie}(K)=\text{Lie}(M)\oplus\text{Lie}(S)$ and the map $M\times S\rightarrow K$ is a diffeomorphism in a neighborhood of the identity, giving a unique decomposition $k= m(k)s(k)$ and also a decomposition of the measure $d\mu_K$, in the sense that $\int_K f d\mu_K=\int_{M\times S}  f d\mu_S d\mu_M$ for any $f$ bounded and compactly supported in this neighborhood, where we denote by $d\mu_M$ the Haar measure on $M$ and by $d\mu_S$ a smooth measure defined on a neighborhood of the identity in $S$.\\
Further, we consider the decomposition  of $G$ as the product $U^-HU$ in a neighborhood of the identity, giving a unique decomposition $s=u^-(s)h(s)u(s)$. We verify that the coordinate map $S\rightarrow U$, $s\mapsto u(s)$ is a diffeomorphism in a neighborhood of the identity. We first observe that
\begin{equation*}
\text{dim}(S)=\text{dim}(K)-\text{dim}(M)=\frac{(d+1)d}{2}-\frac{(d-1)d}{2}=d=\text{dim}(U).
\end{equation*}
Moreover, for the product map $p:U^-\times H \times U \rightarrow G,~ (u^-,h,u)\mapsto u^-hu$, the derivative at the identity is given by $D(p)_e(x,y,z)=x+y+z$, for all $(x,y,z)\in \text{Lie}(U^-) \times\text{Lie}(H) \times\text{Lie}(U)$. Hence, for all $w \in \text{Lie}(G)$, the $U$-component of $D(p)_e^{-1}(w)$ is zero if and only if $w \in \text{Lie}(U^-)+\text{Lie}(H)$. Since $\text{Lie}(S) \cap \left(\text{Lie}(U^-)+\text{Lie}(H)\right)= 0$, the derivative of $s\mapsto u(s)$ is injective. Since dim$(S)$=dim$(U)$, this is a local diffeomorphism.\\
We denote by $B_\mathcal{r}^K$ (by $B_\mathcal{r}^G$ resp.) the ball of radius $\mathcal{r}>0$ centered at the identity in $K$ (in $G$ resp.) and localize the problem to a neighborhood of the identity by considering the partition of unity $1= \sum_{j=1}^{N} \phi_j(kk_j^{-1})$ for all $k \in \text{supp}(f)$ and some $k_j \in$ supp$(f)$, with non-negative functions $\phi_j \in C^{\infty}(K)$ such that supp$(\phi_j) \subseteq B_\mathcal{r}^K$, $|| \phi_j ||_l \ll \mathcal{r}^{-\nu}$ and $N \ll \mathcal{r}^{-\lambda}$, for some $\nu$, $\lambda >0$, and for $\mathcal{r}>0$ small enough to be fixed later.\\
We write for simplicity $k=m_ks_k=m_ku_{s_k}^-h_{s_k}u_{s_k}$, the unique decompositions of $k$ and $s$ in a neighborhood of the identity in $K$ and $S$. We also write $f_j(k)\coloneqq f(kk_j)$ and $\Lambda_j \coloneqq k_j\Lambda$. We compute
\begin{align*}
&I_{\Lambda,f,\varphi_1,\dots,\varphi_r}(t_1,\dots,t_r) = \sum_{j=1}^N \int_K \phi_j(k) f(kk_j)\left(\prod_{i=1}^r\varphi_i(a_{t_i}kk_j\Lambda)\right)d\mu_K(k)\\
&=\sum_{j=1}^N \int_K \phi_j(k) f_j(k)\left(\prod_{i=1}^r\varphi_i(m_ka_{t_i}u^-_{s_k}a_{-t_i}h_{s_k}a_{t_i}u_{s_k}\Lambda_j)\right)d\mu_K(k).
\end{align*}
By Lipschitz continuity of the coordinate maps $m_k$, $u^{-}_{s_k}$ and $h_{s_k}$ on $B_\mathcal{r}^K$ with $\mathcal{r}$ small enough, there exists a constant $C_1>0$ such that for all $k \in B_\mathcal{r}^K$, we have 
$$ a_tu^-_{s_k}a_{-t} \in B_{C_1\mathcal{r}e^{-2t}}^G~ ~~~ \text{and}~~~ m_k,h_{s_k} \in B_{C_1\mathcal{r}}^G.$$
By Lipschitz continuity of $\varphi_1,\dots,\varphi_r$, it follows
\begin{align*}
   & \left|~I_{\Lambda,f,\varphi_1,\dots,\varphi_r}(t_1,\dots,t_r) - \sum_{j=1}^N \int_K \phi_j(k) f_j(k)\left(\prod_{i=1}^r\varphi_i(a_{t_i}u_{s_k}\Lambda_j)\right)d\mu_K(k) ~\right|~\\
    &=~\left|\sum_{j=1}^N \int_K \phi_j(k) f_j(k)\left(\prod_{i=1}^r\varphi_i(a_{t_i}u_{s_k}\Lambda_j)-\prod_{i=1}^r\varphi_i(m_ka_{t_i}u^-_{s_k}a_{-t_i}h_{s_k}a_{t_i}u_{s_k}\Lambda_j\right)d\mu_K(k) ~\right|~\\
    &\ll_l ~  \mathcal{r}\left\lVert f \right\rVert_l\left\lVert\prod_{i=1}^r\varphi_i\right\rVert_l \int_K \left|\sum_{j=1}^{N}\phi_j(k)\right|d\mu_K(k)\\
    &= ~ \mathcal{r} \left\lVert f \right\rVert_l\left\lVert\prod_{i=1}^r\varphi_i\right\rVert_l \qquad \text{(by $K$-invariance and since $\phi_j$ is a partition of unity)}\\
    &\ll_r ~ \mathcal{r} \left\lVert f \right\rVert_l\prod_{i=1}^r\left\lVert\varphi_i\right\rVert_l.
\end{align*}

We now use the decomposition of $\mu_K$ and apply the change of variable $u \mapsto s(u)=s_{u}$, with a density $\rho$ defined in a neighborhood of the identity in $U$ by
\begin{equation*}
\int_S \Phi(s)d\mu_S(s)=\int_U\Phi(s(u))\rho(u)d\mu_U(u) \quad \text{for all } \Phi \in C_c(S) \text{ with supp}(\Phi)\subset B_\mathcal{r}^S. 
\end{equation*}
We have
\begin{align}
&\sum_{j=1}^N \int_K \phi_j(k) f_j(k)\left(\prod_{i=1}^r\varphi_i(a_{t_i}u_{s_k}\Lambda_j)\right)d\mu_K(k) \nonumber \\
&=~\sum_{j=1}^N \int_{M\times S} \phi_j(ms) f_j(ms)\left(\prod_{i=1}^r\varphi_i(a_{t_i}u_{s}\Lambda_j)\right)d\mu_S(s)d\mu_M(m) \nonumber\\
&= ~\int_M\left(\sum_{j=1}^N \int_U \phi_j(ms_u) f_j(ms_u)\left(\prod_{i=1}^r\varphi_i(a_{t_i}u\Lambda_j)\right)\rho(u) d\mu_U(u)\right)d\mu_M(m). \label{integral}
\end{align}
Using Theorem \ref{double equidistribution U orbits} with the function $f_{m,j}(u) \coloneqq \phi_j(ms_u)\rho(u)f_j(mu)$ and observing that $||f_{m,j} ||_l \ll ||\phi_j||_l||\rho||_l ||f_j||_l$ and that $||\rho||_l\ll 1 $, it follows that the integral (\ref{integral}) is equal to
\begin{align}
&\int_M\sum_{j=1}^N \Biggl(\int_U \phi_j(ms_u) f_j(ms_u)\rho(u)d\mu_U(u)\left(\prod_{i=1}^r \int_{\mathcal{X}}\varphi_i d\mu_\mathcal{X}\right) \Biggl.\nonumber\\
& \quad\quad\quad\quad\quad\quad\quad\quad\quad\quad\quad\quad+ \Biggr. O\left( e^{-\gamma D(t_1,\dots,t_r)}||\phi_j||_l||f||_l\prod_{i=1}^r||\varphi_i||_l\right)\Biggr)d\mu_M(m)\nonumber\\
&= \int_M\left(\sum_{j=1}^N \int_S \phi_j(ms)f_j(ms)d\mu_S(s)\right)d\mu_M(m) \left(\prod_{i=1}^r \int_{\mathcal{X}}\varphi_id\mu_\mathcal{X}\right)  \nonumber\\
& \quad\quad\quad\quad\quad\quad\quad\quad\quad\quad\quad\quad\quad\quad\quad+  O\left( Ne^{-\gamma D(t_1,\dots,t_r)}||\phi_j||_l||f||_l\prod_{i=1}^r||\varphi_i||_l\right).\label{estimate}
\end{align}
Using again the decomposition of $\mu_K$, $K$-invariance and the partition of unity, we have
\begin{align*}
    \int_M\left(\sum_{j=1}^N \int_S \phi_j(ms)f_j(ms)d\mu_S(s)\right)d\mu_M(m) &= \int_K\left(\sum_{j=1}^N  \phi_j(kk_j^{-1})\right)f(k)d\mu_K(k)\\
    &= \int_K f d\mu_K
\end{align*} 
which simplifies the estimate \eqref{estimate} to
$$\left( \int_K fd\mu_K\right) \left(\prod_{i=1}^r \int_{\mathcal{X}}\varphi_id\mu_\mathcal{X}\right)  + O\left( \mathcal{r}^{-\lambda} e^{-\gamma D(t_1,\dots,t_r)}\mathcal{r}^{-\nu}||f||_l\prod_{i=1}^r||\varphi_i||_l\right).
$$
Altogether, we obtain
\begin{align*}
I_{\Lambda,f,\varphi_1,\dots,\varphi_r}(t_1,\dots,t_r) =&\left(\int_K fd\mu_K\right) \left(\prod_{i=1}^r \int_{\mathcal{X}}\varphi_id\mu_\mathcal{X}\right)   \\
&+ O\left( \left( \mathcal{r}^{-\lambda-\nu} e^{-\gamma D(t_1,\dots,t_r)}+\mathcal{r} \right)||f||_l\prod_{i=1}^r||\varphi_i||_l\right).
\end{align*} 
We take $\mathcal{r}=e^{-\delta D(t_1,\dots,t_r)}$ with $\delta = \frac{\gamma}{1+\lambda+\nu}$, which yields the claim. 
\end{proof}

We will use the following simplified version of Theorem \ref{double equidistribution K orbits}.

\begin{cor}
\label{cor:multiple equidistribution K orbits}
For every $r\geq 1$, there exist $\delta_r>0$ and $l_r\geq 1$ such that for every $\varphi_0,\dots \varphi_r \in C_{c}^{\infty}(\mathcal{X})$ and $t_1,\dots t_r > 0$, we have
\begin{align*}
    \int_{\mathcal{Y}} \varphi_0(y)\left(\displaystyle\prod_{i=1}^r\varphi_i (a_{t_i} y)\right) d\mu_\mathcal{Y}(y) =& \int_{\mathcal{Y}}\varphi_0 d\mu_\mathcal{Y} \left(\prod_{i=1}^{r}\int_{\mathcal{X}} \varphi_i d\mu_\mathcal{X}\right)\\
    &+ O\left(e^{-\delta_r D(t_1,\dots,t_r)}\prod_{i=0}^{r}||\varphi_i||_{l_r}\right)~.
\end{align*}
\end{cor}

We recall in the following section some properties of the Siegel transform that we use later to analyze the ergodic averages $\sum_{t=0}^N \widehat \chi \circ a_t $.

\subsection{Siegel transform and approximation of the counting function}
\label{sec:approx counting function}

\subsubsection{Siegel transform}
\label{siegel transform}

Given a bounded measurable function $f:\mathbb{R}^{d+2} \rightarrow \mathbb{R}$ with compact support, its (standard) Siegel transform on the space $t$ of unimodular lattices in $\mathbb{R}^{d+2}$ is defined by
\begin{equation}
\label{def standard siegel transform}
{\widehat f}^{\text{st.}}(\Lambda) \coloneqq \sum_{z \in \Lambda \setminus \{0\}} f(z), \quad \text{for } \Lambda \in t.
\end{equation}
Its restriction to $\mathcal{X}$ is called the \emph{light-cone Siegel transform}, defined for a bounded and compactly supported function $f$ on $\mathcal{V}$ by
\begin{equation}
\label{def siegel transform}
\widehat {f}(\Lambda) \coloneqq \sum_{z \in \Lambda \setminus \{0\}} f(z), \quad \text{for } \Lambda \in \mathcal{X}.
\end{equation}
The Siegel transform of a bounded function is typically unbounded, but its growth rate is controlled by an explicit function $\alpha$ defined as follows.\\

Given a lattice $\Lambda \in t$, we say that a subspace $V$ of $\mathbb{R}^{d+2}$ is $\Lambda$-\emph{rational} if the intersection $ V\cap\Lambda$ is a lattice in $V$. If $V$ is $\Lambda$-rational, we denote $d_\Lambda(V)$ the covolume of $V\cap \Lambda$ in $V$. We define then
$$
\alpha(\Lambda) \coloneqq \sup \left\lbrace d_{\Lambda}(V)^{-1}: V \text{ is a } \Lambda\text{-rational subspace of } \mathbb{R}^{d+2}  \right\rbrace.
$$ 

It follows from Mahler's Compactness Criterion that $\alpha$  is a proper map $t \rightarrow [1, +\infty)$. We recall below some important properties.

\begin{prop}[\cite{Schmidt1968AsymptoticFF}]
\label{alpha_growth}
If $f:\mathbb{R}^{d+2} \rightarrow \mathbb{R}$ is a bounded function with compact support, then
$$|\widehat{f}^{\text{st.}}(\Lambda)| \ll_{supp(f)} ||f||_{\infty} \alpha(\Lambda), \quad \text{ for all } \Lambda \in t.
$$
\end{prop}

We restrict this function to the space $\mathcal{X}$ of lattices on the positive light cone and denote it also by $\alpha$. An important property of $\alpha$ is its L$^{p}$-integrability in $t$ (see \cite{eskin1998upper}) and also in $\mathcal{X}$ with an explicit non-escape of mass.

\begin{prop}[\cite{ouaggag2022effective}]
\label{alpha_Lp}
The function $\alpha$ is in \emph{L}$^{p}(\mathcal{X})$ for $1\leq p < d$. In particular,
$$ \mu_{\mathcal{X}} (\{\alpha \geq L \}) \ll_p L^{-p},
\quad \text{ for all } p<d.$$
\end{prop}

We also have an analog of the Siegel Mean Value Theorem (\cite{Siegel1945AMV}) for the space $\mathcal{X}$ .

\begin{prop}[\cite{ouaggag2022effective}]
\label{siegel mean value thrm}
If $f:\mathcal{V} \rightarrow \mathbb{R}$ is a bounded Riemann integrable function with compact support, then
$$ \int_{\mathcal{X}} \widehat{f}(\Lambda) d\mu_{\mathcal{X}}(\Lambda) = \int_{\mathcal{V}} f(z) dz  $$
for some $G$-invariant measure $dz$ on $\mathcal{V}$.
\end{prop}

\subsubsection{Non-divergence estimates}
\label{subsection non-divergence}
We recall here important estimates for the Siegel transform $\widehat{f}$ on translated $K$-orbits by analyzing the escape of mass on submanifolds $a_t\mathcal{Y} \subset \mathcal{X}$.

Following the same argument as in \cite{bjoerklund2018central} and using effective equidistribution of translated $K$-orbits and L$^p$-integrability of the function $\alpha$, we verified in \cite{ouaggag2022effective} the following analogous non-escape of mass for $a_t\mathcal{Y}$.

\begin{prop}[\cite{ouaggag2022effective}]
\label{non-espace of mass}
There exists $\kappa >0$ such that for every $L\geq 1$ and $t\geq \kappa \log L$, 
$$\mu_{\mathcal{Y}} (\{ y \in \mathcal{Y} : \alpha(a_t y) \geq L  \}) \ll_p L^{-p}, \quad \text{ for all } p<d.
$$
\end{prop}

A crucial ingredient in our argument later is the integrability of the Siegel transform $\widehat{f}$ on $a_t\mathcal{Y}$ uniformly in $t$. This is an important result of Eskin, Margulis and Mozes in \cite{eskin1998upper} establishing the following estimate for the function $\alpha$.

\begin{prop}[\cite{eskin1998upper}]
\label{bounds siegel transform}
If $d \geq 2$ and $0<p<2$, then for any lattice $\Lambda$ in $\mathbb{R}^{d+2}$,
$$ \sup_{t> 0} \int_{K} \alpha ( a_t k \Lambda)^p d\mu_{K}(k) < \infty.$$
\end{prop}

\subsubsection{Truncated Siegel transform}
	The Siegel transform of a smooth compactly supported function is typically not bounded. To be able to apply equidistribution results, we truncate the Siegel transform using a smooth cut-off function $\eta_L$ built on the function $\alpha$. We use the same construction as in \autocite[Lemma 4.9]{bjoerklund2018central} which yields the following lemma.

\begin{lem}
For every $\xi>1$, there exists a family $(\eta_L)$ in C$_c^{\infty}(\mathcal{X})$ satisfying:
$$ 0 \leq \eta_L \leq 1 , \quad \eta_L = 1 \text{ on } \{ \alpha \leq \xi^{-1}L \} , \quad \eta_L = 0 \text{ on } \{ \alpha > \xi L \} , \quad || \eta_L||_{C^{l}} \ll 1.\\
$$
\end{lem}

For a bounded function $f:\mathcal{V} \rightarrow \mathbb{R}$ with compact support, we define the \emph{truncated Siegel transform} of $f$ by
$$
\widehat{f}^{(L)} \coloneqq \widehat{f}\cdot \eta_L.
$$

We recall in the following proposition some properties of the truncated Siegel transform $\widehat{f}^{(L)}$ which we use later in our arguments.  

\begin{prop}[\cite{ouaggag2022effective}, except estimate \eqref{L2(X) bound}]
\label{bounds truncated siegel transform}
For a bounded measurable function $f:\mathcal{V} \rightarrow \mathbb{R}$ with compact support, the truncated Siegel transform $\widehat{f}^{(L)}$ satisfies the following bounds:
\begin{align}
& ||\widehat{f}^{(L)} ||_{\infty} \ll_{\emph{supp}(f)} L||f ||_{\infty},\label{second estimate}\\
& ||\widehat{f}^{(L)} ||_{L^p_\mathcal{X}} \leq ||\widehat{f} ||_{L^p_\mathcal{X}}\ll_{\emph{supp}(f),p} ||f ||_{\infty}~, ~~ \text{for all }  p<d,\label{first estimate}\\
& \sup_{t\geq 0}||\widehat{f}^{(L)}\circ a_t ||_{L^{p}_\mathcal{Y}} \leq \sup_{t\geq 0}||\widehat{f}\circ a_t ||_{L^{p}_\mathcal{Y}} < \infty ~, ~~ \text{for all } 1\leq p<2,\label{third estimate}\\
& ||\widehat{f} - \widehat{f}^{(L)} ||_{L^{1}_\mathcal{X}} \ll_{\emph{supp}(f),\tau} L^{-(\tau-1)}||f ||_{\infty} \quad \text{, for all } \tau < d,\\
&\lVert \widehat{f}-\widehat{f}^{(L)}\rVert_{L^{2}_\mathcal{X}} \ll_{\emph{supp}(f),\tau}  L^{-\frac{\tau-2}{2}}||f ||_{\infty}\quad \text{, for all } \tau < d,\label{L2(X) bound}\\
&||\widehat{f}\circ a_t - \widehat{f}^{(L)}\circ a_t ||_{L^{p}_\mathcal{Y}} \ll_{\emph{supp}(f),\tau} L^{-\frac{\tau(2-p)}{2p}}||f ||_{\infty} ~, ~~ \text{for all } 1\leq p<2, ~\tau<d \text{ and } t\geq \kappa \log L. \label{L1 bound}
\end{align}
Moreover, if $f \in C_c^{\infty}(\mathcal{V})$ then $\widehat{f}^{(L)} \in C_c^{\infty}(\mathcal{X})$ and satisfies
\begin{align}
\label{bound for Cl norm truncated}
||\widehat{f}^{(L)} ||_{C^l} \ll_{\emph{supp}(f)} L||f ||_{C^l} \quad \text{, for all } l\geq 1.
\end{align}
\end{prop}

\begin{proof}
All but estimate \eqref{L2(X) bound} were proven in \cite{ouaggag2022effective}.\\
To show \eqref{L2(X) bound}, we apply Hölder's Inequality with $1 \leq p < n$ and $q=(1/2-1/p)^{-1}$ and deduce
$$
\lVert \widehat{f}-\widehat{f}^{(L)} \rVert_{L^2_\mathcal{X}} \ll_{\emph{supp}(f)}  ||\alpha ||_{L_{\mathcal{X}}^p} ~ \mu_{\mathcal{X}}(\{\alpha \geq \xi^{-1}L \})^{1/q}~||f ||_{\infty}.
$$
Then Proposition \ref{alpha_Lp} implies
$$
\lVert \widehat{f}-\widehat{f}^{(L)}\rVert_{L^{2}_\mathcal{X}} \ll_{\emph{supp}(f),p}  L^{-\frac{p-2}{2}}||f ||_{\infty}.
$$
\end{proof}

\subsubsection{Smooth approximation}
For simplicity we write $\chi\coloneqq \chi_{F_{1,c}}$ and $\chi_t\coloneqq \chi_{F_{1,c,t}}$ for the characteristic functions of the sets $F_{1,c}$ and $F_{1,c, t}$ respectively. We approximate $\chi$ and $\chi_t$ by a family of non-negative functions $f_{\varepsilon} , f_{t,\varepsilon}\in C_{c}^{\infty}(\mathcal{V})$ with support in an $\varepsilon$-neighborhood of $F_{1,c}$ and $F_{1,c, t}$ respectively, such that
\begin{align}
\label{epsilon approximation of chi}
&\chi \leq f_{\varepsilon} \leq 1, \quad ||f_{\varepsilon} - \chi||_{\text{L}^{1}_\mathcal{V}} \ll \varepsilon , \quad ||f_{\varepsilon} - \chi||_{\text{L}^{2}_\mathcal{V}} \ll \varepsilon^{1/2}, \quad ||f_{\varepsilon}||_{C^{l}} \ll \varepsilon^{-l},\\
\text{and} \qquad & \chi_t \leq f_{t,\varepsilon} \leq 1, \quad ||f_{t,\varepsilon} - \chi_t||_{\text{L}^{1}_\mathcal{V}} \ll \varepsilon , \quad ||f_{t,\varepsilon} - \chi_t||_{\text{L}^{2}_\mathcal{V}} \ll \varepsilon^{1/2}, \quad ||f_{t,\varepsilon}||_{C^{l}} \ll \varepsilon^{-l}.\label{epsilon approximation of chi_t}
\end{align}

We reformulate in the following proposition a previous result in \cite{ouaggag2022effective}, in order to take into account the parameter $t\geq 1$, and show that the smooth approximation of $\chi_t$ ($\chi$ respectively) yields a good approximation of its Siegel transform $\widehat{\chi_t}$ ($\widehat{\chi}$ respectively) on translated $K$-orbits, uniformly in the parameter $t\geq 0$.
\begin{prop}
\label{smooth approximation}
There exists $\theta>0$ such that for every $t\geq 0$ and every $\varepsilon>0$,
$$\int_{\mathcal{Y}} \left|\widehat{\chi_t}\circ a_t - \widehat{f_{t,\varepsilon}} \circ a_t \right| d\mu_{\mathcal{Y}} \ll_{c,d} \varepsilon +e^{-\theta t}.
$$
\end{prop}

\begin{proof}
Let $t\geq 1$. We first recall the definition of the set $F_{1,c,t}$,
\begin{align*}
    &F_{1,c,t} = \{ x \in \mathcal{V} : x_{d+2}^{2}- x_{d+1}^{2} < c_t^{2}~,~ c\leq x_{d+2}+ x_{d+1} < ce \},\\
    \text{with}\quad & c_t=c\cdot \left(\frac{t}{t+1} \right)^2,
\end{align*}
and observe that there exists $c_{t,\varepsilon}>c_t$ such that $c_{t,\varepsilon}=c_t+ O(\varepsilon)$ and $f_{t,\varepsilon}\leq \chi_{t,\varepsilon}$, where $\chi_{t,\varepsilon}$ denotes the characteristic function of the set
$$
\left\lbrace x \in \mathcal{V} ~:~  c-\varepsilon \leq x_{d+2}+x_{d+1}\leq ce+\varepsilon ,~   x_{d+2}^2-x_{d+1}^2 < c_{t,\varepsilon}^2    \right\rbrace.
$$
The difference $\chi_{t,\varepsilon}-\chi_t$ is bounded by the sum $\chi^{(1)}_{t,\varepsilon}+\chi^{(2)}_{t,\varepsilon}+\chi^{(3)}_{t,\varepsilon}$ of the characteristic functions of the sets
\begin{align*}
&\left\lbrace x \in \mathcal{V} ~:~ c-\varepsilon \leq x_{d+2}+x_{d+1}\leq c ,~   x_{d+2}^2-x_{d+1}^2 < c_{t,\varepsilon}^2    \right\rbrace,\\
&\left\lbrace x \in \mathcal{V} ~:~ ce \leq x_{d+2}+x_{d+1}\leq ce+\varepsilon,~   x_{d+2}^2-x_{d+1}^2 < c_{t,\varepsilon}^2    \right\rbrace,\\
&\left\lbrace x \in \mathcal{V} ~:~ c \leq x_{d+2}+x_{d+1}\leq ce ,~   c_t^2<x_{d+2}^2-x_{d+1}^2 < c_{t,\varepsilon}^2    \right\rbrace.
\end{align*} 
Since $0\leq \chi_t \leq f_{t,\varepsilon} \leq \chi_{t,\varepsilon}$, it follows in particular that 
$$
\widehat{f_{t,\varepsilon}} (a_t\Lambda) - \widehat{\chi_t}(a_t\Lambda) \leq \widehat{\chi_{t,\varepsilon}^{(1)}}(a_t\Lambda)+\widehat{\chi_{t,\varepsilon}^{(2)}}(a_t\Lambda)+\widehat{\chi_{t,\varepsilon}^{(3)}}(a_t\Lambda).
$$
We first consider $\chi^{(1)}_{t,\varepsilon}$. For $x$ in the corresponding set, we also have
\begin{equation*}
0\leq x_{d+2}-x_{d+1} < c_{t,\varepsilon}^2/(c-\varepsilon) \quad \text{and}\quad x_1^2+ \dots +x_{d}^2 < c_{t,\varepsilon}^2.\\
\end{equation*}
We write $I_{0,t,\varepsilon}\coloneqq [0,c_{t,\varepsilon}] $, $I_{1,t,\varepsilon}\coloneqq [-c_{t,\varepsilon}^2/(c-\varepsilon),0]$, $I_{2,\varepsilon}\coloneqq [c-\varepsilon,c] $, $ k=(k_1, \dots , k_{d+2})^T \in K$, and compute
\begin{align*}
&\int_{\mathcal{Y}} |\widehat{\chi_{t,\varepsilon}^{(1)}}\circ a_t| ~d\mu_{\mathcal{Y}} = \int_{K} \widehat{\chi_{t,\varepsilon}^{(1)}}(a_tk\Lambda_0) ~d\mu_K(k) = \int_{K}\sum_{z \in \Lambda_0} \chi_{t,\varepsilon}^{(1)}(a_tkz) ~d\mu_K(k) \\
&=   \sum_{z \in \Lambda_0 }\int_{K} \chi^{(1)}_{t,\varepsilon}\left(\begin{matrix}
 \langle k_{1},z \rangle\\ \dots\\  \langle k_{d},z \rangle\\
 \langle k_{d+1},z \rangle\cosh t - z_{d+2}\sinh t\\
 \langle k_{d+1},z \rangle(-\sinh t) + z_{d+2}\cosh t 
\end{matrix}\right)~d\mu_K(k) \\
&\leq\sum_{z \in \Lambda_0} 
\int_K 
\chi_{I_{0,t,\varepsilon}}\left( ||\langle k_1,z\rangle ,\dots, \langle k_d,z \rangle || \right) 
\chi_{I_{1,t,\varepsilon}}\left(e^{t}\left(\langle k_{d+1},z \rangle -z_{d+2}\right) \right)\cdot\\
&\qquad\qquad\qquad\qquad\qquad\qquad\qquad\qquad\quad
\cdot \chi_{I_{2,\varepsilon}}\left(e^{-t}\left( \langle k_{d+1},z \rangle +z_{d+2}\right) \right)
d\mu_K(k).
\end{align*}
We observe that the intersection $(e^{-t}I_{1,t,\varepsilon}+z_{d+2}) \cap (e^{t}I_{2,\varepsilon}-z_{d+2})$ is non-empty only if $(c-\varepsilon)e^t\leq 2z_{d+2} \leq ce^t+\frac{c_{t,\varepsilon}^2}{c-\varepsilon}e^{-t}$, i.e. $z_{d+2}=ce^t/2 +O_c(\varepsilon e^t+e^{-t})$ where the implicit constant is independent from $t\geq 1$. Moreover, writing each $z \in \Lambda_0$ as $z = z_{d+2}k_zv_0$ with some $k_z \in K$ and $v_0=(0,\dots,0,1,1)\in \mathcal{V}$, and using invariance under $k_z$, we have
\begin{align}
&\int_{\mathcal{Y}} |\widehat{\chi_{t,\varepsilon}^{(1)}}\circ a_t| ~d\mu_{\mathcal{Y}}\nonumber\\
&\leq \sum_{\substack{z\in \Lambda_0\\ z_{d+2}=\frac{ce^t}{2} +O(\varepsilon e^t+e^{-t})}} \int_K 
 \chi_{I_{0,t,\varepsilon}}\left( z_{d+2}||\langle k_1,v_0\rangle ,\dots, \langle k_d,v_0 \rangle || \right)
\chi_{e^{-t}I_{1,t,\varepsilon}}\left(z_{d+2}(\langle k_{d+1},v_0 \rangle -1) \right)\cdot \nonumber\\
&  \qquad \qquad \qquad \qquad \qquad \qquad \qquad \qquad \qquad \qquad \cdot \chi_{e^{t}I_{2,\varepsilon}}\left(z_{d+2}(\langle k_{d+1},v_0 \rangle +1) \right)
d\mu_K(k) \nonumber \\
&\leq \sum_{\substack{z\in \Lambda_0\\ z_{d+2}=\frac{ce^t}{2} +O(\varepsilon e^t+e^{-t})}} \int_K 
\chi_{e^{-t}\frac{2}{c-\varepsilon}I_{0,t,\varepsilon}}\left(||\langle k_1,v_0\rangle ,\dots, \langle k_d,v_0 \rangle || \right)
\chi_{e^{-2t}\frac{2}{c-\varepsilon} I_{1,t,\varepsilon}}\left(\langle k_{d+1},v_0\rangle - 1 \right) \cdot \nonumber\\
& \qquad \qquad \qquad \qquad \qquad \qquad \qquad \qquad \qquad \qquad \cdot \chi_{\frac{2}{c-\varepsilon}I_{2,\varepsilon}}\left(\langle k_{d+1},v_0 \rangle +1 \right)
d\mu_K(k) \nonumber \\
&\leq \sum_{\substack{z\in \Lambda_0\\ z_{d+2}=\frac{ce^t}{2} +O(\varepsilon e^t+e^{-t})}} \mu_K \left( \left\lbrace k \in K ~ : \begin{matrix}  |k_{i,d+1}| \ll_c e^{-t}, ~ i=1,\dots,d  , \\
|k_{d+1,d+1}-1| \ll_c e^{-2t}.
\end{matrix} \right\rbrace \right)\nonumber \\
&\leq \sum_{\substack{z\in \Lambda_0\\ z_{d+2}=\frac{ce^t}{2} +O(\varepsilon e^t+e^{-t})}} \mu_K \left( \left\lbrace k \in K ~ : ||kv_0-v_0 ||\ll_c e^{-t} \right\rbrace \right)\nonumber \\
&\ll_d \sum_{\substack{z\in \Lambda_0\\ z_{d+2}=\frac{ce^t}{2} +O(\varepsilon e^t+e^{-t})}} \mu_{S^d} \left( \left\lbrace v \in S^d ~ :  ||v-v_0 ||\ll e^{-t} \right\rbrace \right)\nonumber \\
&\ll_{c,d} \sum_{\substack{z\in \Lambda_0\\ z_{d+2}=\frac{ce^t}{2} +O(\varepsilon e^t+e^{-t})}} e^{-dt}.\label{sum light-cone}
\end{align}
To estimate the number of summands in \eqref{sum light-cone}, we use a classical upper bound for the number $r_{d+1}(N)$ of representations of a square $N^2\geq 1$ as a sum of $d+1$ squares (Hardy-Littlewood for $d\geq 4$ in \cite{hardy1920partitioI}, \cite{HardyLittlewood1922PNIV}, \cite{HardyLittlewood1925PNVI}, and Kloosterman for $d= 3$ in \cite{kloosterman1927representation}). There exists a constant $C>0$, depending only on the dimension, such that
$$
r_{d+1}(N)= C N^{d-1}+O_{\varepsilon}\left(N^{\frac{d+1}{2}+\varepsilon}\right).
$$
Summing over $N=0,\dots,T-1$, we obtain that there exists $\theta>0$ such that, for all $d\geq 3$, we have
\begin{equation*}
\Big|\big\{z \in \mathcal{V}\cap \mathbb{Z}^{d+2}~: 0\leq z_{d+2}<T \big\}\Big| = CT^{d}+O\left(T^{d-\theta}\right), 
\end{equation*}
hence 
\begin{equation*}
\Big|\big\{z \in \mathcal{V}\cap \mathbb{Z}^{d+2}~: (c-\varepsilon)e^t \leq 2z_{d+2}<ce^t+\frac{c_{t,\varepsilon}^2}{c-\varepsilon}e^{-t} \big\}\Big| \leq \varepsilon e^{dt}+O_c\left(e^{(d-\theta)t}\right).
\end{equation*}
It follows
\begin{equation}
\int_{\mathcal{Y}} |\widehat{\chi_{t,\varepsilon}^{(1)}}\circ a_t| ~d\mu_{\mathcal{Y}} \ll_{c,d} \varepsilon + e^{-\theta t}.
\end{equation}
\\
We proceed similarly for $\chi^{(3)}_{t,\varepsilon}$. For $x$ in the corresponding set, we also have
\begin{equation*}
\frac{c_t^2}{ce}\leq x_{d+2}-x_{d+1} < c_{t,\varepsilon}^2 \quad \text{and}\quad c_t^2<x_1^2+ \dots +x_{d}^2 < c_{t,\varepsilon}^2.\\
\end{equation*}
We write $I'_{0,t,\varepsilon}\coloneqq [c_t,c_{t,\varepsilon}] $, $I'_{1,t,\varepsilon}\coloneqq [-c_{t,\varepsilon}^2,-c_t^2/ce]$, $I'_{2}\coloneqq [c,ce] $ and compute similarly,
\begin{align*}
&\int_{\mathcal{Y}} |\widehat{\chi_{t,\varepsilon}^{(3)}}\circ a_t| ~d\mu_{\mathcal{Y}}  = \int_{K}\sum_{z \in \Lambda_0} \chi_{t,\varepsilon}^{(3)}(a_tkz) ~d\mu_K(k) \\
&=\sum_{z \in \Lambda_0} 
\int_K 
\chi_{I'_{0,t,\varepsilon}}\left( ||\langle k_1,z\rangle ,\dots, \langle k_d,z \rangle || \right) 
\chi_{I'_{1,t,\varepsilon}}\left(e^{t}\left(\langle k_{d+1},z \rangle -z_{d+2}\right) \right)\cdot
\\
&\qquad\qquad\qquad\qquad\qquad\qquad\qquad\qquad\qquad\cdot \chi_{I'_{2}}\left(e^{-t}\left( \langle k_{d+1},z \rangle +z_{d+2}\right) \right)
d\mu_K(k).
\end{align*}
We observe again that the intersection $(e^{-t}I'_{1,t,\varepsilon}+z_{d+2}) \cap (e^{t}I'_{2}-z_{d+2})$ is non-empty only if $C_1e^t \leq z_{d+2} \leq C_2e^t$ for some positive constants $C_1$ and $C_2$ depending only on $c>0$. Moreover, writing each $z \in \Lambda_0$ as $z = z_{d+2}k_zv_0$ with some $k_z \in K$ and $v_0=(0,\dots,0,1,1)\in \mathcal{V}$, and using invariance under $k_z$, we have
\begin{align}
&\int_{\mathcal{Y}} |\widehat{\chi_{t,\varepsilon}^{(3)}}\circ a_t| ~d\mu_{\mathcal{Y}}\nonumber\\
&\leq \sum_{\substack{z\in \Lambda_0\\ z_{d+2} \asymp_c e^t}} \int_K 
 \chi_{I'_{0,t,\varepsilon}}\left( z_{d+2}||\langle k_1,v_0\rangle ,\dots, \langle k_d,v_0 \rangle || \right)
\chi_{e^{-t}I'_{1,t,\varepsilon}}\left(z_{d+2}(\langle k_{d+1},v_0 \rangle -1) \right)\cdot \nonumber\\
& \qquad \qquad \qquad \qquad \qquad \qquad \qquad \qquad \qquad \cdot \chi_{e^{t}I'_{2}}\left(z_{d+2}(\langle k_{d+1},v_0 \rangle +1) \right)
d\mu_K(k) \nonumber \\
&\leq \sum_{\substack{z\in \Lambda_0\\ z_{d+2} \asymp_c e^t}} \int_K 
\chi_{e^{-t}\frac{1}{C_1}I'_{0,t,\varepsilon}}\left(||\langle k_1,v_0\rangle ,\dots, \langle k_d,v_0 \rangle || \right)
\chi_{e^{-2t}\frac{1}{C_1}I'_{1,t,\varepsilon}}\left(\langle k_{d+1},v_0\rangle - 1 \right) \cdot \nonumber\\
& \qquad\qquad \qquad \qquad \qquad \qquad \qquad \qquad \qquad \qquad \cdot \chi_{\frac{1}{C_1}I'_{2}}\left(\langle k_{d+1},v_0 \rangle +1 \right)
d\mu_K(k) \nonumber \\
&\leq \sum_{\substack{z\in \Lambda_0\\ z_{d+2} \asymp_c e^t}} \mu_K \left( \left\lbrace k \in K ~ : ||kv_0-v_0 ||\ll_c \varepsilon
e^{-t} \right\rbrace \right)\nonumber \\
&\ll_d \sum_{\substack{z\in \Lambda_0\\ z_{d+2} \asymp_c e^t}} \mu_{S^d} \left( \left\lbrace v \in S^d ~ :  ||v-v_0 ||\ll_c \varepsilon e^{-t} \right\rbrace \right)\nonumber \\
&\ll_{c,d} \sum_{\substack{z\in \Lambda_0\\ z_{d+2} \asymp_c e^t}} \varepsilon^d e^{-dt}. \nonumber 
\end{align}
We use again the estimate
\begin{equation*}
|\{z \in \mathcal{V}\cap \mathbb{Z}^{d+2}~: z_{d+2}\asymp e^{t} \}|= O(e^{dt}), 
\end{equation*}
hence 
\begin{equation*}
\int_{\mathcal{Y}} |\widehat{\chi_{t,\varepsilon}^{(3)}}\circ a_t| ~d\mu_{\mathcal{Y}} ~\ll_{c,d}~ \varepsilon.
\end{equation*}
The bound for $||\widehat{\chi_{t,\varepsilon}^{(2)}}\circ a_t ||_{L^1_\mathcal{Y}}$ is obtained similarly to that for $\chi^{(1)}_{t,\varepsilon}$.\\
Altogether, we obtain
$$||\hat{f}_{t,\varepsilon} \circ a_t -\hat{\chi_t}\circ a_t ||_{L^1_\mathcal{Y}} \ll_{c,d} \varepsilon+ e^{-\theta t}.
$$
\end{proof}

The same argument as in the proof of Proposition \ref{smooth approximation}, with $c>0$ instead of $(c_t)_{t\geq 0}$ gives the same estimate for the test function $\chi$ instead of $\chi_t$.

\begin{prop}
\label{smooth approximation for c}
There exists $\theta>0$ such that for every $t\geq 0$ and every $\varepsilon>0$,
$$\int_{\mathcal{Y}} \left|\widehat{\chi}\circ a_t - \widehat{f_{\varepsilon}} \circ a_t \right| d\mu_{\mathcal{Y}} \ll_{c,d} \varepsilon +e^{-\theta t}.
$$
\end{prop}

\subsection{Estimating the variance}
\label{subsec:estimating variance}
In this section, we shall estimate and show the convergence of the variance of the averaging function $\mathsf{F}_{N,M}^{(\varepsilon, L)}$, given by
\begin{align*}
	\left\|\mathsf{F}_{N,M}^{(\varepsilon, L)}\right\|_{L^2(\cY)}^2 =&\frac{1}{N-M}\sum_{t_1=M}^{N-1}\sum_{t_2=M}^{N-1} \int_{\cY}
	\psi_{t_1}^{(\varepsilon, L)}\psi_{t_2}^{(\varepsilon, L)}\, d\mu_\cY,
\end{align*}
with
$$
\psi_{t}^{(\varepsilon, L)} :=  \widehat{f}_{\varepsilon}^{(L)}\circ a_t -\mu_\mathcal{Y}(\widehat{f}_{\varepsilon}^{(L)}\circ a_t).
$$
We first observe that this expression is symmetric with respect to $t_1$ and $t_2$,
writing $t_1=s+t$ and $t_2=s$ with $0\le t\le N-M-1$ and $M\le s \le N-t-1$,
we obtain that
\begin{align}\label{eq:FF_N_L}
\left\|\mathsf{F}_{N,M}^{(\varepsilon, L)}\right\|_{L^2(\cY)}^2 
=& \Theta_{N,M}^{(\varepsilon,L)}(0)
+2\sum_{t=1}^{N-M-1}\Theta_{N,M}^{(\varepsilon,L)}(t),
\end{align}
where
$$
\Theta_{N,M}^{(\varepsilon,L)}(t):=\frac{1}{N-M}\sum_{s=M}^{N-1-t} \int_{\cY}\psi_{s+t}^{(\varepsilon, L)}\psi_s^{(\varepsilon, L)}\, d\mu_\cY,
$$
with
\begin{align*}\int_{\cY}\psi_{s+t}^{(\varepsilon, L)}\psi_s^{(\varepsilon, L)}\, d\mu_\cY
=
\int_{\cY} (\widehat f_{\varepsilon}^{(L)}\circ a_{s+t})(\widehat f_{\varepsilon}^{(L)}\circ a_s)\, d\mu_\cY -
\mu_\cY(\widehat f_{\varepsilon}^{(L)}\circ a_{s+t})\mu_\cY(\widehat f_{\varepsilon}^{(L)}\circ a_s).
\end{align*}
We shall first show that with a suitable choice of parameters $\varepsilon$ and $L$ we have:
\begin{equation}
\label{eq:final estimate variance F_N,L}
\left\| \mathsf{F}_{N,M}^{(\varepsilon, L)}\right\|_{L^2(\cY)}^2 \\
=\Theta_\infty^{(\varepsilon, L)} (0) + 2\sum_{t=1}^{N-1}\Theta_\infty^{(\varepsilon, L)} (t)+o(1),
\end{equation}
where
$$
\Theta_\infty^{(\varepsilon, L)}(t):=\int_{\cX} (\widehat f_{\varepsilon}^{(L)}\circ a_t)\widehat f_{\varepsilon}^{(L)}\, d\mu_\cX -\mu_\cX(\widehat f_{\varepsilon}^{(L)})^2.
$$
To estimate $\Theta_{N,M}(t)$, we introduce an additional parameter $K=K(N)\to \infty$ (to be specified later) with $K\le M$ and consider separately the cases when $t< K$ and when $t\ge  K$. \\

First, we consider the case when $t\ge K$. 
By Corollary \ref{cor:multiple equidistribution K orbits}, we have 
\begin{align}\label{eq:second0}
\int_{\cY} (\widehat f_{\varepsilon}^{(L)}\circ a_{s+t})(\widehat f_{\varepsilon}^{(L)}\circ a_s)\, d\mu_\cY
&=\mu_{\cX}(\widehat f_{\varepsilon}^{(L)})^2
 +O\left(e^{-\delta \min(s,t)} \left\|\widehat f_{\varepsilon}^{(L)}\right\|_{C^l}^2 \right).
\end{align}
and also
\begin{align}\label{eq:second}
	\int_{\cY}  (\widehat f_{\varepsilon}^{(L)}\circ a_s)\, d\mu_\cY
	&=\mu_{\cX} (\widehat f_{\varepsilon}^{(L)}) +O\left( e^{-\delta s} \left\|\widehat f_{\varepsilon}^{(L)}\right\|_{C^l} \right).
\end{align}
Hence, combining \eqref{eq:second0} and \eqref{eq:second}, we deduce that
\begin{align*}
	\int_{\cY}\psi_{s+t}^{(\varepsilon, L)}\psi_s^{(\varepsilon, L)}\, d\mu_\cY
	=O \left(e^{-\delta \min(s,t)} \left\|\widehat f_{\varepsilon}^{(L)}\right\|_{C^l}^2\right).
\end{align*}
Since
\begin{align*}
\sum_{t=K}^{N-M-1}\left(\sum_{s=M}^{N-1-t}e^{-\delta \min(s,t)}\right)
&\le \sum_{t=K}^{N-1}\sum_{s=M}^{N-1} (e^{-\delta t}+e^{-\delta s})\ll N e^{-\delta K},
\end{align*}
we conclude that
\begin{equation}
\label{eq:Theta1}
\sum_{t=K}^{N-M-1}\Theta_{N,M}(t) \ll \, e^{-\delta K}\, \left\|\widehat f_{\varepsilon}^{(L)}\right\|_{C^l}^2
\ll \, e^{-\delta K} L^{2}\, \varepsilon^{-2l},
\end{equation}
where we used Lemma \ref{bounds truncated siegel transform} and (\ref{epsilon approximation of chi}).
The implied constants here and below in the proof depend only on $\supp(f_\varepsilon)$, which is uniformly bounded, hence the dependence is only on the constant $c$ from the Diophantine approximation \eqref{intrinsic Diophantine}. \\

Let us now consider the case $t<K$. We observe that Corollary \ref{cor:multiple equidistribution K orbits} (for $r = 1$)
applied to the function $\phi_{t}:= (\widehat f_{\varepsilon}^{(L)}\circ a_t)\widehat f_{\varepsilon}^{(L)}$ yields, 
\begin{align*}
\int_{\cY} (\widehat f_{\varepsilon}^{(L)}\circ a_{s+t})(\widehat f_{\varepsilon}^{(L)}\circ a_s)\, d\mu_\cY
&=\int_{\cY} (\phi_{t}\circ a_s) \, d\mu_\cY\\
&=\int_{\cX} \phi_{t}\, d\mu_\cX+O\left(e^{-\delta s}\, \left\|\phi_{t}\right\|_{C^l} \right).
\end{align*}
Furthermore, for some $\xi=\xi(n,l)>0$, we have
$$
\left\|\phi_{t}\right\|_{C^l}\ll \left\|\widehat f_{\varepsilon}^{(L)}\circ a_t\right\|_{C^l}\, \left\|\widehat f_{\varepsilon}^{(L)}\right\|_{C^l}\ll e^{\xi t}\,
\left\|\widehat f_{\varepsilon}^{(L)}\right\|_{C^l}^2.
$$
Therefore, we deduce that
$$
\int_{\cY} (\widehat f_{\varepsilon}^{(L)}\circ a_{s+t})(\widehat f_{\varepsilon}^{(L)}\circ a_s)\, d\mu_\cY
=\int_{\cX} (\widehat f_{\varepsilon}^{(L)}\circ a_t)\widehat f_{\varepsilon}^{(L)} \,d\mu_\cX +O \left(e^{-\delta s} e^{\xi t}\, \left\|\widehat f_{\varepsilon}^{(L)}\right\|_{C^l}^2 \right).
$$
Combining this estimate with \eqref{eq:second}, we obtain that
\begin{align*}
\int_{\cY}\psi_{s+t}^{(\varepsilon,L)}\psi_s^{(\varepsilon,L)}\, d\mu_\cY
=\Theta_\infty^{(\varepsilon, L)}(t)
+ O \left(e^{-\delta s} e^{\xi t}\, \left\|\widehat f_{\varepsilon}^{(L)}\right\|_{C^l}^2  \right).
\end{align*}
Using further the estimates from Lemma \ref{bounds truncated siegel transform} and (\ref{epsilon approximation of chi}), it follows, for the case $t<K$,
\begin{align*}
\Theta_{N,M}^{(\varepsilon,L)}(t)=&\, \frac{N-M-t}{N-M} 
\Theta_\infty^{(\varepsilon,L)}(t)+O \left((N-M)^{-1} e^{-\delta M} e^{\xi t}\, \left\|\widehat f_{\varepsilon}^{(L)}\right\|_{C^l}^2  \right)\\
=& \, \Theta_\infty^{(\varepsilon,L)}(t)
+O \left((N-M)^{-1}t\left\|\widehat f_{\varepsilon}^{(L)}\right\|_{L^2(\cX)}^2+ (N-M)^{-1} e^{-\delta M} e^{\xi t}\, \left\|\widehat f_{\varepsilon}^{(L)}\right\|_{C^l}^2   \right)\\
=& \, \Theta_\infty^{(\varepsilon,L)}(t)
+O \left((N-M)^{-1}t+ (N-M)^{-1} e^{-\delta M} e^{\xi t}\, \varepsilon^{-2l}L^2  \right).
\end{align*}
It follows
\begin{align}
\label{eq:Theta2}
\Theta_{N,M}^{(\varepsilon,L)}(0)+2\sum_{t=1}^{K-1}\Theta_{N,M}^{(\varepsilon,L)}(t) =& \,
\Theta_\infty^{(\varepsilon,L)}(0)+2\sum_{t=1}^{K-1}\Theta_\infty^{(\varepsilon,L)} (t) \\
&+O \left ((N-M)^{-1}K^2+ (N-M)^{-1} e^{-\delta M}  e^{\xi K}\varepsilon^{-2l} L^{2}\right).\nonumber
\end{align}
Combining \eqref{eq:Theta1} and \eqref{eq:Theta2}, it follows from \eqref{eq:FF_N_L} that 
\begin{align}\label{eq:starr}
\left\| \mathsf{F}_{N,M}^{(\varepsilon,L)}\right\|_{L^2(\cY)}^2 
&= \,
\Theta_\infty^{(\varepsilon,L)}(0)+2\sum_{t=1}^{K-1}\Theta_\infty^{(\varepsilon,L)} (t) \nonumber\\
& + O \left((N-M)^{-1}K^2 \,+ ((N-M)^{-1} e^{-\delta M}  e^{\xi K}+e^{-\delta K}) L^{2}\,\varepsilon^{-2l} \right).
\end{align}
We will choose later in \eqref{eq:choice M} to \eqref{eq:choice K} the parameters $K(N)$, $M(N)$, $\varepsilon(N)$ and $L(N)$ so that
\begin{align}
e^{-\delta K} L^{2}\varepsilon^{-2l} \to 0, \label{eq:cc7}\\
(N-M)^{-1} e^{-\delta M}  e^{\xi K} L^{2}\varepsilon^{-2l} \to 0, \label{eq:cc8}\\
(N-M)^{-1}K^2 \to 0, \label{eq:cc9}
\end{align}
as $N\to \infty$, which gives 
\begin{equation*}
\left\| \mathsf{F}_{N,M}^{(\varepsilon,L)}\right\|_{L^2(\cY)}^2 =\Theta_\infty^{(\varepsilon,L)} (0) + 2\sum_{t=1}^{K-1}\Theta_\infty^{(\varepsilon,L)} (t)+o(1).
\end{equation*}
We shall show next that with a suitable choice of parameters we have
\begin{equation}
\label{eq:final estimate variance F_epsilon^N}
\left\| \mathsf{F}_{N,M}^{(\varepsilon,L)}\right\|_{L^2(\cY)}^2 =\Theta_\infty^{(\varepsilon)} (0) + 2\sum_{t=1}^{K-1}\Theta_\infty^{(\varepsilon)} (t)+o(1),
\end{equation}
where 
$$
\Theta_\infty^{(\varepsilon)}(t):=\int_{\cX} (\widehat f_{\varepsilon}\circ a_t)\widehat f_{\varepsilon}\, d\mu_\cX -\mu_\cX(\widehat f_{\varepsilon})^2.
$$
Using again the estimates from Lemma \ref{bounds truncated siegel transform} 
$$
||\widehat f_\varepsilon - \widehat f_\varepsilon^{(L)} ||_{L^{1}_\mathcal{X}} \ll_{\text{supp}(f_\varepsilon),\tau} L^{-(\tau-1)} \quad \hbox{and }\quad \lVert \widehat f_\varepsilon-\widehat f_\varepsilon^{(L)}\rVert_{L^{2}_\mathcal{X}} \ll_{\text{supp}(f_\varepsilon),\tau}  L^{-\frac{\tau-2}{2}}
$$
we have (since the supports of the functions $f_\varepsilon
$ are uniformly bounded) 
\begin{align*}
    \mu_\cX(\widehat f_\varepsilon^{(L)})&=\mu_\cX(\widehat f_\varepsilon)+O_{\tau}\left( L^{-\tau-1}\right),\\
    \int_{\cX} (\widehat f_{\varepsilon}^{(L)}\circ a_t)\widehat f_{\varepsilon}^{(L)}\, d\mu_\cX &= \int_{\cX} (\widehat f_{\varepsilon}\circ a_t)\widehat f_{\varepsilon}\, d\mu_\cX + O_{\tau}\left( L^{-\frac{\tau-2}{2}}\right),
\end{align*}
which yields
$$
\Theta_\infty^{(\varepsilon,L)}(t)= \Theta_\infty^{(\varepsilon)}(t) + O_{\tau}\left( L^{-\frac{\tau-2}{2}}\right),
$$
and \eqref{eq:starr} then gives
\begin{align}\label{eq:starrr}
\left\| \mathsf{F}_{N,M}^{(\varepsilon,L)}\right\|_{L^2(\cY)}^2 
&= \,
\Theta_\infty^{(\varepsilon)}(0)+2\sum_{t=1}^{K-1}\Theta_\infty^{(\varepsilon)} (t) \nonumber\\
& + O \left((N-M)^{-1}K^2 \,+ ((N-M)^{-1} e^{-\delta M}  e^{\xi K}+e^{-\delta K}) L^{2}\,\varepsilon^{-2l} +KL^{-\frac{\tau-2}{2}}\right).
\end{align}

We will choose the parameters $K(N)$ and $L(N)$ such that, besides the conditions on $L,\varepsilon,M,K$ assumed previously, we also have 
\begin{equation}
    KL^{-\frac{\tau-2}{2}}\rightarrow 0\quad \text{as }~N\rightarrow\infty, \quad \hbox{for some}\quad \tau<d, \label{eq:condition KL}
\end{equation}
which gives \eqref{eq:final estimate variance F_epsilon^N}. 
\\ \\
In order to analyze further the correlations $ \int_{\cX} (\widehat f_{\varepsilon}\circ a_t)\widehat f_{\varepsilon}\, d\mu_\cX$ and show the convergence of the series $\sum_{t=1}^{K-1}\Theta_\infty^{(\varepsilon)} (t)$, we shall use results from a recent work by Kelmer and Yu in \cite{secondmomentKY23},  where incomplete Eisenstein series are used to analyze the second moment of the light-cone Siegel transform. We briefly recall in the following section some preliminaries to this approach.\\

\subsubsection{Preliminaries about Eisenstein's series}

Before recalling the approach and results of Kelmer and Yu, we reproduce below some preliminaries from \cite{secondmomentKY23} about Eisenstein's series and adapt the notations to our coordinate system from Section \ref{correspondance}.\\
To simplify notations, we will denote in this section the elements in the subgroup $A$ by $$a_y:= \begin{pmatrix}
 I_d  & & \\
  & \frac{y+y^{-1}}{2}  & -\frac{y-y^{-1}}{2} \\
 &-\frac{y-y^{-1}}{2}  & \frac{y+y^{-1}}{2} \end{pmatrix} \quad \text{, for } y>0, $$\\
the $\mathbb{R}$-split torus with $a_y$ acting on $e_0=(0,\dots,1,1)\in \mathbb{R}^{d+2}$ as $a_ye_0=y^{-1}e_0$, and denote by
$$
U=\left\{u_{x}=\left(\begin{matrix} I_d & -x & x\\
 x^t&1-\frac{1}{2}\lVert x \rVert^2  &\frac{1}{2}\lVert x \rVert^2 \\ 
 x^t& -\frac{1}{2}\lVert x \rVert^2& 1+\frac{1}{2}\lVert x \rVert^2\end{matrix}\right): x\in \mathbb{R}^d\right\}
$$
its unstable horospherical subgroup. We denote by
\begin{align}\label{equ:macomptgpk}
K:=\left\{k=\begin{pmatrix} \tilde{k} &  \\  &1\end{pmatrix}: \tilde{k}\in \SO_{d+1}(\mathbb{R})\right\}
\end{align}
a maximal compact subgroup. Let $L$ be the stabilizer of $e_0$ in $G$ and let $P$ be the parabolic subgroup fixing the line spanned by $e_0$. More precisely, $P=MAU$ and $L=MU$ with
$$
M=\left\{m=\left(\begin{smallmatrix}  \tilde{m}&  &  \\  &1& \\  & &1\end{smallmatrix}\right): \tilde{m}\in \SO_d(\mathbb{R})\right\}
$$ 
the centralizer of $A$ in $K$, and $U$ the unipotent subgroup given by the Iwasawa decomposition $G=KAU$.\\ 

Any $g\in G$ can be written as $g=ka_yu_{x}$ with $k\in K$, $u_x\in U$, and in these coordinates, the Haar measure of $G$ is given (up to scaling) by
\begin{align}
\label{equ:Haar1}
d\mu_G(g)=y^{-(d+1)}\, d x d y d\mu_K(k),
\end{align}
where $dx$ is the usual Lebesgue measure on $\mathbb{R}^d$ and $\mu_K$ is the probability Haar measure of $K$. 

The subgroup $L$ is unimodular, with its Haar measure given by 
\begin{align*}
d\mu_L(u_{x}m)=dxd\mu_M(m),
\end{align*}
where $\mu_M$ is the probability Haar measure of $M\cong \SO_d(\mathbb R)$. 
Since $L$ is the stabilizer of $e_0$ and $G$ acts transitively on $\mathcal{V}$, we can identify $\mathcal{V}$ with the homogeneous space $G/L$, which gives a natural right $G$-invariant measure on $\mathcal{V}$. Explicitly, further identifying $G/L$ with $K/M \times \mathbb{R}_+$ gives natural polar coordinates on $\mathcal{V}$: Every $x\in \mathcal{V}$ can be written uniquely as $x=ka_ye_0$ for some $y>0$ and $k\in K/M$. In these coordinates, the measure 
\begin{align}\label{equ:lebmea}
d \mu_{\mathcal{V}}(ka_ye_0):=y^{-(d+1)}\,d yd\mu_{K/M}(k)
\end{align}
is such an invariant measure. Here $\mu_{K/M}$ is the unique right $K$-invariant probability measure on the homogeneous space $K/M$ which is homeomorphic to the unit sphere $S^d$.\\

The cusps of $\Gamma$ are the $\Gamma$-conjugacy classes of rational parabolic subgroups of $G$. Let $m$ be the number of these cusps and $P_1,\dots,P_m$ a set of representatives of these classes, each of which having a Langlands decomposition $P_i=M_iA_iU_i$, $i=1,\dots,m$. We denote by $\Gamma_{P_i}:=\Gamma\cap P_i$ and by $\Gamma_{U_i}:=\Gamma\cap U_i$, where $\Gamma_{U_i}$ is by definition a finite index subgroup of $\Gamma_{P_i}$ (see \cite{secondmomentKY23}).\\
 For each $P_i$ we fix the scaling matrix $\tau_i=a_{y_i}k_i$, where $k_i\in K$ is such that $P_i=k_i^{-1}Pk_i$ and where $y_i>0$ is the unique number such that $\mu_L( L /\tau_i\Gamma_{P_i}\tau_i^{-1})=1$. \\
We define the \emph{(spherical) Eisenstein series} corresponding to the $i$-th cusp for $\textrm{Re}(s)>d $ and $g \in G$ by the convergent series
 $$ 
 E_i(s,g) := \sum_{\gamma \in \Gamma_{P_i}\backslash \Gamma } y( g\gamma\tau_i^{-1})^s,
 $$
 where $y(g)$ is given by the Iwasawa decomposition $g=ka_{y(g)}u_{x}\in KAU$.\\
 For each $1\leq  j\leq m$ the constant term of $E_i(s,g)$ with respect to the $j$-th cusp is defined by
\begin{equation}\label{e:ConstSph} 
	c_{ij}(s,g):=\frac{1}{\text{vol}(U/\tau_j\Gamma_{U_j}\tau_j^{-1})}\int_{U/\tau_j\Gamma_{U_j}\tau_j^{-1}}E_i(s, \tau_j u_{x} g)\,dx,
\end{equation}
which is known to be of the form
\begin{equation}\label{e:ConstSph1}
c_{ij}(s,g)=\delta_{ij}y(g)^s+\varphi_{ij}(s)y(g)^{d-s}
\end{equation}
for some holomorphic functions $\varphi_{ij}$, defined for $\textrm{Re}(s)>d$ and building the \emph{scattering matrix} $(\varphi_{ij})_{1\leq i,j\leq m}$.\\

The series $E_i(s, g)$ (and hence also $\varphi_{ij}$) has a meromorphic continuation to the whole $s$-plane, which, on the half plane $\textrm{Re}(s)\geq \frac{d}{2}$, is holomorphic except for a simple pole at $s=d$ (called the \textit{trivial pole}) and possibly finitely many simple poles on the interval $(\frac{d}{2}, d)$ (called \textit{exceptional poles}). We denote by $C_{\Gamma}\subseteq (\frac{d}{2}, d)$ the finite set of exceptional poles of all Eisenstein series of $\Gamma$. \\

All poles of the Eisenstein series come from the poles of the scattering matrix $\varphi_{ij}(s)$, and for any $1\leq i,j\leq m$,
\begin{align}\label{equ:constomeg}
\text{Res}_{s=n}\varphi_{ij}(s)=\text{Res}_{s=n}E_i(s,g)=\omega_{\Gamma}.
\end{align}
Moreover, for any of the exceptional poles $s_l \in C_\Gamma$, it follows from the Maass-Selberg relations (given in \cite[Equation (7.44)]{selberg1989harmonic} and
\cite[Chapter 6, (1.62)]{sarnak1980}) that 
\begin{equation}\label{eq:scatteringres}
\text{Res}_{s=s_l}\varphi_{ij}(s)=\left\langle \text{Res}_{s=s_l}E_i(s,\cdot),\text{Res}_{s=s_l}E_j(s,\cdot)\right\rangle,
\end{equation}
where the inner product is with respect to $\mu_G$.\\

The residue of $E_i(s,g)$ at $s=d$ is a constant which is the same for Eisenstein series at all cusps, given by the reciprocal of the measure of the homogeneous space $G/\Gamma$, that is, for each $1\leq i\leq m$ and $g\in G$,
\begin{align}\label{equ:rescovol}
\omega_{\Gamma}:=\text{Res}_{s=d}E_i(s,g)=\mu_G(G/\Gamma)^{-1}.
\end{align}
 For any bounded and compactly supported function $f:\mathcal{V}\rightarrow \mathbb{C}$, the \emph{incomplete Eisenstein series attached to $f$ at $P_i$} is defined for any $g \in G$ by
 $$
 E_i(g,f) := \sum_{\gamma \in \Gamma_{P_i}\backslash \Gamma } f(g\gamma \tau_i^{-1}e_0).
 $$
 Since $E_i(\cdot,f)$ is left $\Gamma$-invariant, it can be viewed as a function on the homogeneous space $\mathcal{X}$. The light-cone Siegel transform of $f$ can then be expressed in terms of incomplete Eisenstein series as follows.
 
\begin{lem}[\cite{secondmomentKY23}, Lemma 3.1]
\label{siegel transform as incomplete eisenstein series}
    There exist constants $\lambda_1,\dots,\lambda_m >0$ such that for any bounded and compactly supported function $f:\mathcal{V} \rightarrow \mathbb{C}$, 
    $$
    \widehat f = \sum_{i=1}^{m} E_i(\cdot,f_{\lambda_i}),
  $$
  where $f_\lambda(x):=f(\lambda^{-1}x) $ for any $\lambda>0$.
\end{lem}
We also introduce the \emph{light-cone Eisenstein series} of the quadratic form $Q$ (see \cite{kelmer2023fourier}), defined for Re$(s)>n$ by the series
\begin{align}\label{equ:ligconeeise}
E_Q(s,g)&:=\|e_0\|^s\sum_{v\in \Lambda_0^{\text{prim}}}\|gv\|^{-s}.
\end{align}

Similar to the light-cone Siegel transform, $E_Q(s,g)$ can be written as a sum of Eisenstein series of $G$ at all cusps; see \cite[Corollary 3.2]{kelmer2023fourier}. Explicitly, we have for any Re$(s)>n$,
\begin{align}\label{equy:relaeisen}
E_Q(s,g)&=\sum_{i=1}^{m}\lambda_i^sE_i(s,g),
\end{align}
where $\lambda_1,\ldots, \lambda_{m}$ are as in \ref{siegel transform as incomplete eisenstein series}.
In particular, we see from this relation that $E_Q(s,g)$ has a meromorphic continuation to the whole $s$-plane and is holomorphic on the half plane Re$(s)\geq \frac{d}{2}$ except for a simple pole at $s=d$ with constant residue, 
$$\omega_Q:=\text{Res}_{s=d}E_Q(s,g),$$ 
and at most finitely many simple poles in $(\tfrac{d}{2}, d)$ contained in the set $C_{\Gamma}$. Moreover, it follows from \cite[Theorem 1.8]{kelmer2023fourier} that, in fact, there is at most one such pole located at $s_d=\left \lfloor{\tfrac{d+2}{2}}\right \rfloor$ if $d\geq 3$, and no poles for $d<3$.
The residue $\omega_Q$ then satisfies the relation
\begin{align}\label{equ:omegaq}
\omega_Q=\omega_{\Gamma}\sum_{i=1}^m\lambda_i^n.
\end{align}
We observe, moreover, that $\omega_Q$ coincides with the scaling factor relating the measure $\mu_\mathcal{V}$ in \eqref{equ:lebmea} to the $G$-invariant measure $dz$ introduced in Proposition \ref{siegel mean value thrm}, and we have
\begin{equation*}\label{relation dz and dmu_V}
    dz=\omega_Q d\mu_\mathcal{V}.
\end{equation*}
We also observe (see \cite{secondmomentKY23}, proof of Theorem 1.1), using the expansion \eqref{equy:relaeisen} and taking residues, that 
$$\text{Res}_{s=s_l}E_Q(s,g)=\sum_{j=1}^{m} \lambda_j^{s_l} \text{Res}_{s=s_l}E_j(s,g).$$
Taking an inner product and using \eqref{eq:scatteringres} gives
\begin{align*}
	\|\text{Res}_{s=s_l}E_Q(s,g)\|^2&=\sum_{i,j=1}^{m}\lambda_i^{s_l}\lambda_j^{s_l}\langle \text{Res}_{s=s_l}E_i(s,\cdot),\text{Res}_{s=s_l}E_j(s,\cdot)\rangle\\
	&=\sum_{i,j=1}^{m}\lambda_i^{s_l}\lambda_j^{s_l}\text{Res}_{s=s_l}\varphi_{ij}(s).
\end{align*}
For the exceptional pole $s_d$, we write 
\begin{equation}
    \label{def c_Q}
    c_Q:=\omega_\Gamma\|\text{Res}_{s=s_d}E_Q(s,g) \|^2=\sum_{i,j=1}^{m}\lambda_i^{s_l}\lambda_j^{s_l}\text{Res}_{s=s_d}\varphi_{ij}(s) .
\end{equation}

\subsubsection{Preliminaries about spherical analysis} 

By the classical spherical harmonic analysis, the function space $L^2(S^d)$ decomposes into irreducible $\SO_{d+1}(\mathbb R)$-representations as follows:
\begin{align*}
L^2(S^d)=\bigoplus_{n\geq 0}L^2(S^d,n),
\end{align*}
where $L^2(S^d, n)$ is the space of degree $n$ harmonic polynomials in $d+1$ variables restricted to $S^d$. This, in turn, induces the following decomposition of $L^2(K/M)$ into irreducible $K$-representations
\begin{align*}
	L^2(K/M)=\bigoplus_{n\geq 0} L^2(K/M, n),
\end{align*}
where $L^2(K/M, n)$ is the pre-image of $L^2(S^d,n)$ under the isomorphism between $L^2(K/M)$ and $L^2(S^d)$.
For each $n\geq 0$, we fix an orthonormal basis 
$
\{\psi_{n,l}: 0\leq l\leq \dim_{\mathbb C}L^2(K/M, n)-1\}
$
for $L^2(K/M, n)$. 
For any $f: \mathcal{V}\to \mathbb C$ bounded and compactly supported, let 
\begin{align}\label{equ:spexpan}
	f_{n,l}(y):=\int_{K/M}f(ka_ye_0)\overline{\psi_{n,l}(k)}d\mu_{K/M}(k).
\end{align}
so that $f$ has a \textit{spherical expansion} 
\begin{align}\label{equ:sphexpan}
f(ka_ye_0)=\sum_{n, l \geq 0}f_{n,l}(y)\psi_{n,l}(k)
\end{align}
in $L^2$ and also pointwise if $f$ is smooth.\\
For any function $f$ on $\mathbb R^{+}$, we denote by $$\widetilde{f}(s):=\int_0^{\infty}f(y)y^{-(s+1)}\,d y, \quad \text{, for }s \in \mathbb C$$ its \textit{Mellin transform},  whenever this defining integral is absolutely convergent.\\
Using the spherical expansion we define the following bilinear form for any $f,f': \mathcal{V}\rightarrow \mathbb C$ bounded and compactly supported and any $s\in  (\frac{d}{2}, d)$,
\begin{align}\label{def:mff}
	M_{f,f'}(s):=\sum_{n,l\geq 0}P_n(s)\widetilde{f_{n,l}}(s)\overline{\widetilde{f'_{n,l}}(s)},
\end{align}
with $P_0(s):=1$ and $P_n(s):=\prod_{i=0}^{n-1}\frac{d-s+i}{s+i}$ if $n\geq 1$.\\

The following lemmas give estimates related to the operator $ M_{f,f'}$ which will be useful later for the analysis of $\Theta_\infty(t)$. We write for simplicity $$ M_{f}:=M_{f,f}.$$
\begin{lem}
\label{lem: epsilon estimate for M_(chi,f_epsilon)}
    Let $f$ be a bounded function on the light-cone with bounded support. For every $s\in \left(\frac{d}{2},d\right)$, we have
    $$\left| M_{f}(s)\right| \ll_{s,\text{supp}(f)} \lVert f \rVert_2^2 ~.$$
\end{lem}
\begin{proof}
By definition of the Mellin transform and using the spherical expansion  $ f(ka_ye_0)=\sum_{n,l}f_{n,l}(y)\psi_{n,l}(k)$, we have for any $s\in (\frac{d}{2},d)$ 
\begin{align*}
    &M_{f}(s) = \sum_{n,l}P_n(s)\left|\widetilde{{f}_{n,l}}(s)\right|^2\\
    &= \sum_{n,l}P_n(s) \left|\int_0^{+\infty} {f}_{n,l}(y)y^{-(s+1)}\,dy\right|^2\\
    &= \sum_{n,l}P_n(s) \left|\int_0^{+\infty} \left(\int_{K/M}f(ka_ye_0)\overline{\psi_{n,l}(k)}\, dk\right)y^{-(s+1)}\,dy\right|^2.
\end{align*}
Using the decomposition $\mathcal{V} \cong  K/M \times \mathbb{R}_+$ given by $x=ka_ye_0$ with the spherical coordinates $y\in \mathbb{R}_+$ and $k\in K/M$, we can write $f(ka_ye_0)=\phi_{y}(k)\rho(y)$, where $(\phi_{y})_{y>0}$ is a family of bounded functions on $K/M$ and $\rho$ is the characteristic function of an interval away from $y=0$ (since by the parametrization of $\mathcal{V}$ we have $a_ye_0=y^{-1}e_0$). \\
We also introduce the projection operator $pr_n:L^2(K/M)\rightarrow L^2(K/M,n)$ on the space of degree $n$ harmonic polynomials in $d+1$ variables restricted to $K/M$ and write $f_n:=pr_n(f)$ for $f\in L^2(K/M)$. Using that $(\psi_{n,l})_{l\geq 0} $ is an orthonormal basis of $L^2(K/M,n)$ for every $n\geq 0$, it follows
\begin{align}
    &M_{f}(s) = \sum_{n,l\geq 0}P_n(s) \left|\int_0^{+\infty} \left(\int_{K/M}\phi_{y}(k)\rho(y)\overline{\psi_{n,l}(k)}\, dk\right)y^{-(s+1)}\,dy\right|^2 \nonumber\\
    &=\sum_{n\geq 0}P_n(s)\int_0^{+\infty}\int_0^{+\infty}\langle {\phi_{y_1}}_n, {\phi_{y_2}}_n\rangle_{_{K/M}}\; \rho(y_1)\rho(y_2)y_1^{-(s+1)}y_2^{-(s+1)}\, dy_1dy_2.\label{computation of M_f_epsilon}
\end{align}
Using that $ \left|P_n(s) \right|\ll 1 $, Cauchy-Schwarz inequality, the decomposition in spherical harmonics given by $L^2(K/M)= \bigoplus_{n\geq 0} L^2(K/M,n)$, that $\rho^2=\rho $ and that the support of $\rho$ is an interval away from $y=0$,  we obtain
\begin{align*}
    \left|M_{f}(s)\right| &\ll \int_0^{+\infty}\int_0^{+\infty}\sum_{n\geq 0}\lVert {\phi_{y_1}}_n\rVert_2\lVert {\phi_{y_2}}_n\rVert_2\; \rho(y_1)\rho(y_2)y_1^{-(s+1)}y_2^{-(s+1)}\, dy_1dy_2\\
    &\leq \int_0^{+\infty}\int_0^{+\infty}\lVert \phi_{y_1}\rVert_2 \lVert \phi_{y_2}\rVert_2\; \rho(y_1)\rho(y_2)y_1^{-(s+1)}y_2^{-(s+1)}\, dy_1dy_2\\
    &= \left(\int_0^{+\infty}\lVert \phi_{y}\rVert_2  \rho(y)y^{-(s+1)}\, dy\right)^2\\
    &\leq \left(\int_0^{+\infty}\lVert \phi_{y}\rVert_2^2  \rho(y)^2y^{-(d+1)}\, dy\right)\left(\int_0^{+\infty}  \rho(y)^2y^{d-2s-1}\, dy\right)\\
    &\ll_{s,\text{supp}(\rho)} \left(\int_0^{+\infty}\int_{K/M} \left|\phi_{y}(k)\rho(y)\right|^2y^{-(d+1)}\,dk dy\right)\\
    &= \lVert f \rVert_2^2.
\end{align*}
\end{proof}

\begin{lem}[\cite{secondmomentKY23}, sections 2.3 and 2.4]\label{estimates for M_lambda,f and mu(lambda.f)} For any $f,f':\mathcal{V}\to \mathbb{C}$ bounded and compactly supported functions, we have, with the same notation as in Lemma \ref{siegel transform as incomplete eisenstein series},
    \begin{align*}
    M_{f_{\lambda_1},f'_{\lambda_2}}(s)&=\lambda_1^s\lambda_2^sM_{f,f'}(s),\quad \text{for any } \lambda_1,\lambda_2>0, \text{ any }s\in (\frac{d}{2},d),\\
    \qand \qquad \mu_\mathcal{V} (f_\lambda)&= \lambda^d \mu_\mathcal{V} (f), \; \text{for any }\; \lambda>0.
\end{align*}
\end{lem}
\begin{lem}[\cite{LogLaws2017Yu}, Lemma 4.4]
For any $s\in (\frac{d}{2},d)$, we have $P_n(s)\asymp_s (n+1)^{d-2s}$ for all $n\geq 1$    \label{estimate P_n(s) real}
\end{lem}
It will be useful for our argument later to have an estimate of $P_n(s)$ also for $s\in \mathbb C$ with real part in $(\frac{d}{2},d)$.
\begin{lem}
For any $s=r+it \in \mathbb{C}$, with $r\in(\frac{d}{2},d)$, we have $\left|P_n(s)\right|\ll_{\delta} |t|^{1/2}(n+1)^{d-2r+\delta}$ for any $\delta>0$.   
\label{estimate P_n(s) complexe}
\end{lem}
\begin{proof}
    We have
    \begin{align*}
        \left|P_n(s) \right|^2
        &=\prod_{k=0}^{n-1}  \frac{|d-s+k|^2}{|s+k|^2}\\
        &=\prod_{k=0}^{n-1} \frac{(d-r+k)^2}{(r+k)^2}\cdot \frac{1+\frac{t^2}{(d-r+k)^2}}{1+\frac{t^2}{(r+k)^2}}\\
        &\ll (n+1)^{2(d-2r)}\cdot \prod_{k=0}^{n-1}\frac{1+\frac{t^2}{(d-r+k)^2}}{1+\frac{t^2}{(r+k)^2}} \qquad (\text{by Lemma } \ref{estimate P_n(s) real})
    \end{align*}
    Further, we have
    \begin{align*}
        \log \left( \prod_{k=0}^{n-1}\frac{1+\frac{t^2}{(d-r+k)^2}}{1+\frac{t^2}{(r+k)^2}}\right) &= \sum_{k=1}^{n-1}\left( \log \left(1+\frac{t^2}{(d-r+k)^2} \right)- \log \left(1+\frac{t^2}{(r+k)^2} \right)\right)\\
        &\qquad\qquad\qquad+\log \left(\frac{1+\frac{t^2}{(d-r)^2}}{1+\frac{t^2}{r^2}} \right)\\
        &\ll \sum_{k=1}^{n-1}\frac{t^2}{1+\frac{t^2}{(r+k)^2}}\left( \frac{1}{(d-r+k)^2}-\frac{1}{(r+k)^2} \right)+O(1)\\
        &\ll \sum_{k=1}^{n-1}\frac{t^2k^2}{k^2+t^2}\cdot \frac{1}{k^3}+O(1)\\
        &\ll \sum_{k\leq \alpha(t)}\frac{1}{k}+\sum_{k= \alpha(t)+1}^{n-1}\frac{t^2}{(k^2+t^2)k} +O(1)\\
        &\ll \log \alpha(t)+\frac{t^2}{\alpha(t)^2+t^2}\log (n+1)+O(1).
    \end{align*}
    Choosing $\alpha(t)=\left \lceil{\beta |t|}\right \rceil $ with $\beta >0$ large enough, we obtain
    \begin{align*}
        \prod_{k=0}^{n-1}\frac{1+\frac{t^2}{(d-r+k)^2}}{1+\frac{t^2}{(r+k)^2}} \ll \beta |t| \cdot (n+1)^{\frac{1}{\beta^2+1}},
    \end{align*}
    hence $|P_n(s)|\ll_\delta |t|^{1/2}\cdot(n+1)^{d-2r+\delta} $ for any $\delta>0$.
\end{proof}

Let us introduce the projection operator  on the space of degree $n$ harmonic polynomials in $d+1$ variables restricted to $K/M$ and write 
$$pr_n:L^2(K/M)\rightarrow L^2(K/M,n) \quad\text{and}\quad f_n:=pr_n(f)\text{, for any } f\in L^2(K/M).$$
From Lemma \ref{estimate P_n(s) complexe} we deduce the following estimate.

\begin{lem}\label{estimate P(s)-weighted correlations in K/M}
    Let $\phi_1,\phi_2:K/M\to \mathbb{C}$ be two characteristic functions. We have, for any $s=r+it\in \mathbb{C}$ and any $0<\delta<2r-d$,
   $$ \left|\sum_{n\geq 0} P_n(s)\langle {(\phi_{y_1})}_{n}, {(\phi_{t,y_2})}_{n}\rangle_{_{K/M}}\right| \ll_\delta |t|^{1/2}||\phi_{y_1}||_2||\phi_{t,y_2}||_1^{\frac{1}{2}+\frac{2r-d-\delta}{d+1}}.$$
\end{lem}

\begin{proof}
We generalize a similar argument from \cite[Proof of Theorem 1.2.]{LogLaws2017Yu}. We split the summation into two parts with a parameter $D\geq 1$ to be fixed later, and estimate the two parts using the estimate of $P_n(s)$ from Lemma \ref{estimate P_n(s) complexe} and the following inequality from spherical harmonic analysis \cite[inequality (4.4)]{sogge1986oscillatory}, 
$$|| \phi_n||_2 \ll (n+1)^{\frac{d-1}{2}}||\phi||_1, \quad \text{for any } \phi \in L^2(K/M). $$
For the first part of the summation, we use the orthogonality of the projections $pr_n$ and Cauchy-Schwarz inequality to obtain
\begin{align}
   \left|\sum_{n\leq D} P_n(s)\langle {(\phi_{y_1})}_{n}, {(\phi_{t,y_2})}_{n}\rangle_{_{K/M}}\right| &= \left|\sum_{n\leq D} P_n(s)\langle \phi_{y_1}, {(\phi_{t,y_2})}_{n}\rangle_{_{K/M}}\right| \\
   &\ll_\delta |t|^{1/2}\sum_{n\leq D}  (n+1)^{d-2r+\delta}\lVert\phi_{y_1}\rVert_2\lVert{(\phi_{t,y_2})}_n\rVert_2 \nonumber\\
   &\ll |t|^{1/2}\sum_{n\leq D}  (n+1)^{d-2r+\delta+\frac{d-1}{2}}\lVert\phi_{y_1}\rVert_2\lVert\phi_{t,y_2}\rVert_1 \nonumber\\
   &\ll |t|^{1/2}D^{d-2r+\delta+\frac{d-1}{2}+1}\lVert\phi_{y_1}\rVert_2\lVert\phi_{t,y_2}\rVert_1. \label{eq:estimate D1}
\end{align}
For the second part of the summation, we use Cauchy-Schwarz inequality for the sum and the convergence given by $L^2(K/M)= \bigoplus_{n\geq 0} L^2(K/M,n)$. We choose moreover $\delta>0$ small enough, such that $d-2r+\delta<0$. This gives
\begin{align}
   \left|\sum_{n> D} P_n(s)\langle {(\phi_{y_1})}_{n}, {(\phi_{t,y_2})}_{n}\rangle_{_{K/M}}\right| &\ll_\delta |t|^{1/2}\sum_{n> D}  (n+1)^{d-2r+\delta}||{(\phi_{y_1})}_{n}||_2||{(\phi_{t,y_2})}_{n}||_2\nonumber\\
   &\leq |t|^{1/2}D^{d-2r+\delta}\sum_{n\geq 0}||{\phi_{y_1}}_{n}||_2||{\phi_{t,y_2}}_{n}||_2\nonumber\\
   &\leq |t|^{1/2}D^{d-2r+\delta}||\phi_{y_1}||_2||\phi_{t,y_2}||_2.\label{eq:estimate D2}
\end{align}
Using that $||\phi||_2=||\phi||_1^{1/2}$ for any characteristic function $\phi$, we optimize both estimates \eqref{eq:estimate D1} and \eqref{eq:estimate D2} by taking $ D=\max (1,||\phi_{t,y_2}||_1^{-\frac{1}{d+1}} )$ and obtain
\begin{align}
\label{eq:estimate M_f,f' in terms of norms}
   \left|\sum_{n\geq 0} P_n(s)\langle {(\phi_{y_1})}_{n}, {(\phi_{t,y_2})}_{n}\rangle_{_{K/M}}\right| &\ll_\delta |t|^{1/2}||\phi_{y_1}||_2||\phi_{t,y_2}||_1^{\frac{1}{2}+\frac{2r-d-\delta}{d+1}}.
\end{align}
\end{proof}

\subsubsection{Moment formulas of incomplete Eisenstein series}

With the preliminaries introduced above, the correlations of incomplete Eisenstein series can then be estimated as follows.

\begin{thm}[\cite{secondmomentKY23}, Theorem 2.3]
\label{inner product Eisenstein series}
    For any $1\leq i,j\leq m$, there exists a bounded linear operator $\cT_{ij}: L^2(\mathcal{V})\rightarrow L^2(\mathcal{V})$ with operator norm $|| \cT_{ij}||_{\text{op}}\leq 1$ such that for any  $f,f' \in C_c^\infty(\mathcal{V})$,
    \begin{align*}
        \langle E_i(\cdot,f),E_j(\cdot,f') \rangle= \omega_\Gamma^2 \mu_\mathcal{V}(f)\mu_\mathcal{V}(f') 
        &+\omega_\Gamma \langle \delta_{ij}f+\cT_{ij}(f),f'\rangle\\&+\omega_\Gamma \sum_{s_l \in C_\Gamma}M_{f,f'}(s_l)\text{Res}_{s=s_l}\varphi_{i,j}(s)~,
    \end{align*}
    
    where the two inner products are with respect to $\mu_\cX$ and $\mu_\mathcal{V}$ respectively.
    \end{thm}

From Theorem \ref{inner product Eisenstein series} Kelmer and Yu derived the following mean value theorem and effective estimate of the second moment of the light-cone Siegel transform.

\begin{thm}[\cite{secondmomentKY23}, Theorem 1.1]
    \label{thm: second moment siegel tranform}
    Let $f:\mathcal{V} \rightarrow \mathbb C$ be a measurable, bounded and compactly supported function. Then we have
    \begin{equation}
        \int_\cX \widehat f d \mu_\cX = \omega_Q \mu_\mathcal{V}(f).\label{mean value identity}
    \end{equation}
    If we assume moreover that $f$ is smooth, then
    \begin{equation}
        \int_\cX \left|\widehat f\right|^2 d \mu_\cX = \left|\omega_Q \mu_\mathcal{V}(f)\right|^2+c_Q M_{f,f}(s_d) +O\left( \mu_\mathcal{V} \left( |f|^2\right)\right),
    \end{equation}
    where $s_d:= \lfloor \frac{d+2}{2} \rfloor$, the term $M_{f,f}(s)$ is a quadratic form on $f$ given by \eqref{def:mff} and $c_Q$ given by \eqref{def c_Q}.
\end{thm}

We generalize the second moment formula in Theorem \ref{thm: second moment siegel tranform} to measurable, bounded and compactly supported functions in the following Proposition, using our estimate in Lemma \ref{lem: epsilon estimate for M_(chi,f_epsilon)} and a similar argument as in \cite[Proof of corollary 1.2]{secondmomentKY23}.

\begin{prop}
    \label{generalized second moment formula}
    Let $f$ be measurable, bounded and compactly supported functions on $\mathcal{V}$. Then we have
    \begin{equation}
        \int_\cX \left|\widehat f\right|^2 d \mu_\cX = \left|\omega_Q \mu_\mathcal{V}(f)\right|^2+c_Q M_{f,f}(s_d) +O\left( \mu_\mathcal{V} \left( |f|^2\right)\right).
    \end{equation}
\end{prop}

\begin{proof}
    There exists a sequence $(f_i)_{i\in \mathbb N}$ in $C_c^{\infty}(\mathcal{V})$ converging to $f$ in $L^1(\mathcal{V})$ and in $L^2(\mathcal{V})$. By the mean value identity \eqref{mean value identity} it follows that $(\widehat f_i)_{i\in \mathbb{N}}$ converges to $\widehat f$ in $L^1(\mathcal{X})$, hence also pointwise almost everywhere for some subsequence. To show that this convergence is also in $L^2$, we use Theorem \ref{thm: second moment siegel tranform} with the estimate from Lemma \ref{lem: epsilon estimate for M_(chi,f_epsilon)}, and write, for any $i,j \in \mathbb N$,
    \begin{align}
        \left\| \widehat f_i - \widehat f_j \right\|_{L_\mathcal{X}^2} &\ll  \left\| f_i - f_j \right\|_{L^1_\mathcal{V}}^2 + \left\| f_i - f_j \right\|_{L^2_\mathcal{V}}^2 + M_{f_i-f_j}(s_d)\nonumber\\
        &\ll \left\| f_i - f_j \right\|_{L^1_\mathcal{V}}^2 + \left\| f_i - f_j \right\|_{L^2_\mathcal{V}}^2 \quad (\text{by Lemma \ref{lem: epsilon estimate for M_(chi,f_epsilon)}}).\label{bound cauchy sequence}
    \end{align}
    Since $(f_i)_{i\in \mathbb N}$ is a Cauchy sequence in $L^1(\mathcal{V})\cap L^2(\mathcal{V})$, it follows from \eqref{bound cauchy sequence}, that $(\widehat f_i)_{i\in \mathbb N}$ is a Cauchy sequence in $L^2(\mathcal{X})$ and therefore converges to $\widehat f$ in $L^2(\mathcal{X})$. Hence, the second moment formula for $\widehat f$ follows from the second moment formula for $(\widehat f_i)_{i\in \mathbb N}$.  
\end{proof}

\subsubsection{Final estimate of the variance of the counting function}\label{subsec:final estimate variance}

We shall show next that
\begin{equation}
\label{eq:final estimate variance F_N}
\left\| \mathsf{F}_{N,M}^{(\varepsilon,L)}\right\|_{L^2(\cY)}^2 =\Theta_\infty (0) + 2\sum_{t=1}^{K-1}\Theta_\infty (t)+o(1),
\end{equation}
where 
$$
\Theta_\infty(t):=\int_{\cX} (\widehat \chi\circ a_t)\widehat \chi\, d\mu_\cX -\mu_\cX(\widehat \chi)^2.
$$
Using Proposition \ref{generalized second moment formula}, Lemma \ref{lem: epsilon estimate for M_(chi,f_epsilon)} and the estimates in \eqref{epsilon approximation of chi}, we have for any $t\geq 0$
\begin{align*}
    \left|\Theta_\infty(t) - \Theta_\infty^{(\varepsilon)}(t) \right| &\leq  \int_\cX \left|\left(\widehat{f_\varepsilon}\circ a_t \right)\widehat{f_\varepsilon}- \left(\widehat{\chi_t}\circ a_t \right)\widehat{\chi_t}\right| d\mu_\cX + \left| \mu_\cX \left( \widehat{f_\varepsilon}\right)^2- \mu_\cX \left( \widehat{\chi_t}\right)^2\right|\\
    &\ll \left\| \widehat{f_\varepsilon}-\widehat{\chi_t}\right\|_{L_\cX^2}\left(\left\| \widehat{f_\varepsilon}\right\|_{L_\cX^2} +\left\| \widehat{\chi_{ }}\right\|_{L_\cX^2}\right)+ \left\| \widehat{f_\varepsilon}-\widehat{\chi_t}\right\|_{L_\cX^1}\\
    &\ll \left\| \widehat{f_\varepsilon-\chi}\right\|_{L_\cX^2}+ \left\| f_\varepsilon-\chi\right\|_{L_\cX^1}\\
    &\ll \mu_\mathcal{V}\left( \left| f_\varepsilon - \chi\right|\right)^2+ \left|M_{f_\varepsilon-\chi}(s_d)\right|+O\left( \mu_\mathcal{V}\left(\left|f_\varepsilon-\chi \right|^2\right)\right) + \varepsilon\\
    &\ll \varepsilon,
\end{align*}
and \eqref{eq:starrr} then gives
\begin{align}\label{eq:starrrr}
\left\| \mathsf{F}_{N,M}^{(\varepsilon,L)}\right\|_{L^2(\cY)}^2
&= \,
\Theta_\infty(0)+2\sum_{t=1}^{K-1}\Theta_\infty(t) \nonumber\\
& + O \left(N^{-1}(M+K)K \,+ (N^{-1} e^{-\delta M}  e^{\xi K}+e^{-\delta K}) L^{2}\,\varepsilon^{-2l} +KL^{-\frac{\tau-2}{2}} + \varepsilon K \right),
\end{align}
which implies \eqref{eq:final estimate variance F_N}, provided that, besides the conditions on $L,\varepsilon,M,K$ assumed previously, we also have 
\begin{equation}
\label{eq:condition epsilon K}
    \varepsilon K \rightarrow 0, \quad \text{ as }N\rightarrow \infty.
\end{equation}

It remains to show that the series $\sum_{t=1}^{K-1}\Theta_\infty (t)$ converges as $K\rightarrow\infty$.\\ \\
We write
\begin{align}
\label{eq:series of theta_infinity}
    \Theta_\infty (t) = \int_{\cX} \widehat{ \chi_t}\cdot \widehat \chi \;d\mu_\cX - \mu_\cX\left( \widehat{ \chi_t}\right)\mu_\cX\left( \widehat \chi\right).
\end{align}
with $ \chi_t:=\chi\circ a_t $.\\ \\

Using Lemma \ref{siegel transform as incomplete eisenstein series}, we can express the correlations in \eqref{eq:series of theta_infinity} in terms of incomplete  Eisenstein series attached to $\chi_t$ and $\chi$:
\begin{equation}
\label{correlation to eisenstein series}
    \int_{\cX} \widehat \chi_t\cdot \widehat \chi \;d\mu_\cX = \sum_{i,j=1}^{m} \langle E_i(\cdot,\chi_{t,\lambda_i}),E_j(\cdot,\chi_{\lambda_j}) \rangle,
\end{equation}
where the inner product is with respect to $\mu_\cX$.

From \eqref{correlation to eisenstein series}, using Lemma \ref{estimates for M_lambda,f and mu(lambda.f)}, \eqref{def c_Q} and that $s_d$ is the unique exeptional pole in $(\frac{d}{2},d)$, it follows that
\begin{align}
   & \int_{\cX} \widehat \chi_t\cdot \widehat \chi \;d\mu_\cX = \sum_{i,j=1}^{m} \langle E_i(\cdot,\chi_{t,\lambda_i}),E_j(\cdot,\chi_{\lambda_j}) \rangle_{\mu_\cX} \nonumber\\
    &= \sum_{i,j=1}^{m} \Bigg( \omega_\Gamma^2 \mu_\mathcal{V}(\chi_{t,\lambda_i})\mu_\mathcal{V}(\chi_{\lambda_j}) +\omega_\Gamma \langle \delta_{ij}\chi_{t,\lambda_i}+\cT_{ij}(\chi_{t,\lambda_i}),\chi_{\lambda_j}\rangle_{\mu_\mathcal{V}}\nonumber\\
&\qquad\qquad\qquad\qquad\qquad\qquad\qquad\qquad\qquad+\omega_\Gamma \sum_{s_l \in C_\Gamma}M_{\chi_{t,\lambda_i},\chi_{\lambda_j}}(s_l)\text{Res}_{s=s_l}\varphi_{ij}(s)\Bigg)\nonumber\\
&= \omega_\Gamma^2\left(\sum_{i=1}^m\lambda_i^d\right)^2\mu_\mathcal{V}(\chi_t)\mu_\mathcal{V}(\chi) + \omega_\Gamma\sum_{i=1}^m\lambda_i^d\langle \chi_t, \chi\rangle_{\mu_\mathcal{V}} + \omega_\Gamma\sum_{i,j=1}^{m} \langle\cT_{ij}(\chi_{t,\lambda_i}),\chi_{\lambda_j}\rangle_{\mu_\mathcal{V}}     \nonumber \\
&\qquad\qquad\qquad\qquad\qquad\qquad\qquad\qquad\qquad+\omega_\Gamma \sum_{i,j=1}^m \lambda_{i}^{s_d}\lambda_j^{s_d}M_{\chi_{t},\chi}(s_d)\text{Res}_{s=s_d}\varphi_{ij}(s)\nonumber\\
    &= \omega_Q^2 \mu_\mathcal{V}(\chi_t)\mu_\mathcal{V}(\chi) +\omega_Q \langle \chi_t, \chi\rangle_{\mu_\mathcal{V}} +\omega_\Gamma\sum_{i,j=1}^{m} \langle\cT_{ij}(\chi_{t,\lambda_i}),\chi_{\lambda_j}\rangle_{\mu_\mathcal{V}}+c_Q M_{\chi_t,\chi}(s_d). \label{eq:correlation of siegel transforms}
\end{align}
Since  $\chi_t $ and $\chi$ have disjoint supports for all $t\geq 1$, it follows from the mean value identity in Theorem \ref{thm: second moment siegel tranform} and  from \eqref{eq:correlation of siegel transforms} that \eqref{eq:series of theta_infinity} reduces for all $t\geq 1$ to 
\begin{equation}
\label{eq:reduced form variance}
    \Theta_\infty(t) = \omega_\Gamma\sum_{i,j=1}^{m} \langle\cT_{ij}(\chi_{t,\lambda_i}),\chi_{\lambda_j}\rangle_{\mu_\mathcal{V}}+c_Q M_{\chi_t,\chi}(s_d).
\end{equation}

We shall estimate the terms $M_{\chi_t,\chi}(s_d)$ and $\langle\cT_{ij}(\chi_{t,\lambda_i}),\chi_{\lambda_j}\rangle_{\mu_\mathcal{V}}$ next. 

\begin{lem}
 \label{convergence of M_f,f^K} 
 For every $s\in (\frac{d}{2},d)$, there exists $\sigma>0$ (depending on $d$ and $s$) such that for every  $t\geq 1$ we have
 $$\left|M_{\chi_t,\chi}(s)\right| \ll e^{-\sigma t}.
 $$
\end{lem}
\begin{proof} 
By definition of the Mellin transform, we have for any $s\in (\frac{d}{2},d)$ 
\begin{align*}
    &M_{\chi_t,\chi}(s) = \sum_{n,l}P_n(s)\widetilde{{\chi}_{_{n,l}}}(s)\widetilde{{(\chi_t)}_{n,l}}(s)\\
    &= \sum_{n,l}P_n(s) \left(\int_0^{+\infty} {\chi}_{_{n,l}}(y)y^{-(s+1)}\,dy\right)\left(\int_0^{+\infty} {(\chi_t)}_{n,l}(y)y^{-(s+1)}\,dy\right)\\
    &= \sum_{n,l}P_n(s) \left(\int_0^{+\infty} \left(\int_{K/M}\chi(ka_ye_0)\overline{\psi_{n,l}}(k)\, dk\right)y^{-(s+1)}\,dy\right)\cdot\\
    & \qquad\qquad\qquad\qquad\qquad\qquad \cdot\left(\int_0^{+\infty} \left(\int_{K/M}\chi_t(ka_ye_0)\overline{\psi_{n,l}}(k)\, dk\right)y^{-(s+1)}\,dy\right).
\end{align*}
Using the decomposition $\mathcal{V} \cong  K/M \times \mathbb{R}_+ $ given by $x=ka_ye_0$ with the spherical coordinates $y\in \mathbb{R}_+$ and $k\in K/M$, we write 
\begin{equation} \label{decomposition in polar coordinates}
    \chi(ka_ye_0)=\phi_{y}(k)\rho(y)\quad \text{and}\quad \chi_t(ka_ye_0)=\phi_{t,y}(k)\rho_t(y),
\end{equation}
where $\rho$ and $\rho_t$ are the characteristic functions of intervals given by the projections on the last coordinate of the domains $F_{1,c}$ and $a_{e^{-t}}(F_{1,c})$ respectively (i.e. intervals of the form $[\alpha,\beta]$ and $[\alpha e^{-t},\beta e^{-t}]$ for some positive constants $\alpha$ and $\beta$ depending on the Diophantine constant $c$ from \eqref{intrinsic Diophantine}), while $\phi_y$ and $\phi_{y,t}$ are characteristic functions on $K/M$ resulting from the decomposition $\mathcal{V} \cong  K/M \times \mathbb{R}_+$.\\

Using that $(\psi_{n,l})_{l\geq 0} $ is an orthonormal basis of $L^2(K/M,n)$ for every $n\geq 0$, it follows
\begin{align}
    &M_{\chi_t,\chi}(s) = \sum_{n,l\geq 0}P_n(s) \left(\int_0^{+\infty} \left(\int_{K/M}\phi_{y}(k)\rho(y)\overline{\psi_{n,l}}(k)\, dk\right)y^{-(s+1)}\,dy\right)\cdot \nonumber\\
    & \qquad\qquad\qquad\qquad\qquad\qquad \cdot\left(\int_0^{+\infty} \left(\int_{K/M}\phi_{t,y}(k)\rho_t(y)\overline{\psi_{n,l}}(k)\, dk\right)y^{-(s+1)}\,dy\right) \nonumber\\
    &=\int_0^{+\infty}\int_0^{+\infty}\sum_{n\geq 0}P_n(s)\langle {(\phi_{y_1})}_n, {(\phi_{t,y_2})}_n\rangle_{_{K/M}}\; \rho(y_1)\rho_t(y_2)y_1^{-(s+1)}y_2^{-(s+1)}\, dy_1dy_2.\label{computation of M_F,f 1}
\end{align}
Using Lemma \ref{estimate P(s)-weighted correlations in K/M} for $s\in (\frac{d}{2},d)$ then gives 
\begin{align}
    \left|M_{\chi_t,\chi}(s)\right|
    \ll_s \int_0^{+\infty}\int_0^{+\infty}\lVert\phi_{y_1}\rVert_2\lVert\phi_{t,y_2}\rVert_1^{\frac{1}{2}+\frac{2s-d}{d+1}}\; \rho(y_1)\rho_t(y_2)y_1^{-(s+1)}y_2^{-(s+1)}\, dy_1dy_2.\label{computation of M_F,f}
\end{align}
We observe further that, since $e_0a_y=y^{-1}e_0$ and since $F_{1,c} \subseteq \{x \in \mathcal{V}: x_1^2+\dots + x_{d}^2 <c^2 \} $, we have for any $y>0$,
\begin{align}
    \left\| \phi_y \right\|_1 &= \int_{K/M}\chi (ka_ye_0)dk\nonumber\\
    &=\int_{K/M}\chi (y^{-1}ke_0)dk\nonumber\\
    &\leq \int_{K/M}\chi_{\{x_1^2+\dots +x_{d}^2 <y^2c^2\}}(ke_0)dk\nonumber\\
    &\ll_c y^d.\label{estimate phi_y}
\end{align}
We write $\sigma:= \frac{1}{2}-\frac{2s-d}{d+1}-\frac{s}{d}$, such that $\sigma>0$ for any $s\in(\frac{d}{2},d)$ and $d\geq 2$, then combine \eqref{computation of M_F,f} and \eqref{estimate phi_y}, and use that supp$(\rho)$ is bounded away from 0 and that supp$(\rho_t)$ $=[\alpha e^{-t},\beta e^{-t}]$. We obtain 
\begin{align*}
  \left| M_{\chi_t,\chi}(s)\right| &\ll \left(\int_0^{+\infty}y_1^{\frac{d}{2}}y_1^{-(s+1)}\rho(y_1)\;dy_1\right)\left(\int_0^{+\infty}y_2^{s+d\sigma}y_2^{-(s+1)}\rho_t(y_2)\, dy_2\right)\\
  &\ll_{\text{supp}(\rho)} \int_0^{+\infty}y_2^{-1+d\sigma}\; \rho_t(y_2)dy_2 \\
  &\ll \int_{\alpha e^{-t}}^{\beta e^{-t}}y_2^{-1+d\sigma}\;dy_2\\
  &\ll e^{-d\sigma t}.
\end{align*}
\end{proof}

We give next an estimate of the term $\sum_{i,j=1}^m\langle\cT_{ij}(\chi_{\lambda_i}),\chi_{t,\lambda_j}\rangle_{\mu_\mathcal{V}}$. In order to simplify the notations, we will omit without loss of generality the scaling coefficients $\lambda_1,\dots, \lambda_m$ and consider only a single cusp.

\begin{lem}
\label{convergence of <tau(X_t),X>}
    There exists $\gamma >0$ such that for every $t\geq 1$ we have
    $$\left|\langle\cT(\chi),\chi_t\rangle_{\mu_\mathcal{V}}\right| \ll e^{-\gamma t}.
    $$
\end{lem}
\begin{proof}
    
    For a smooth and compactly supported function $f\in C_c^\infty(\mathcal{V})$, the operator $\cT_{}$ can be expressed explicitly in terms  of the Mellin transform of the spherical harmonic coefficients of $f$ as follows (see \cite[Proof of Theorem 2.3]{secondmomentKY23}):
    $$
    \cT_{}(f)(ka_ye_0) = \sum_{n,l}\left(\frac{1}{2\pi i}\int_{\left( \frac{d}{2}\right)} P_n(s)\varphi(s)\widetilde{{f}_{n,l}}(s)y^{d-s}ds\right)\psi_{n,l}(k)
    $$
    where the contour integration is along the line of complex numbers with real part $\frac{d}{2}$.\\
    We shall approximate $ \chi$ by a smooth and compactly supported function $f_{\varepsilon_t}$ in the sense of \eqref{epsilon approximation of chi}, with a parameter $\varepsilon_t>0$ to be fixed later. We note that $\varepsilon_t$ is independent from the parameter $\varepsilon=\varepsilon(N)$ introduced in \eqref{epsilon approximation of chi}. Since $\cT$ is a bounded linear operator on $L^2(\mathcal{V}) $ with operator norm $\lVert\cT \rVert_{\text{op}} \leq 1$, we have:
    \begin{align*}
        \left| \langle \cT_{}(\chi),\chi_t\rangle\right| &= \left| \langle \cT_{}(\chi -f_{\varepsilon_t})+\cT_{}(f_{\varepsilon_t}),\chi_t\rangle\right|\\
        &\leq \left\| \chi -f_{\varepsilon_t}\right\|_2\left\|\chi_t \right\|_2+\left| \langle \cT_{}(f_{\varepsilon_t}),\chi_t\rangle\right| \\
        &\ll  \varepsilon_t^{1/2} + \left| \langle \cT_{}(f_{\varepsilon_t}),\chi_t\rangle\right| \; .
    \end{align*}
    Using the decomposition $d\mu_\mathcal{V}=y^{-(d+1)}dydk$, the spherical expansion $ f(ka_ye_0)=\sum_{n,l}f_{n,l}(y)\psi_{n,l}(k)$ and the decomposition $L^2(K/M)= \bigoplus_{n\geq 0} L^2(K/M,n)$, it follows:
    \begin{align*}
       &\langle \cT_{}(f_{\varepsilon_t}),\chi_t\rangle_{\mu_\mathcal{V}} \\
        &= \int_{K/M}\int_0^{+\infty}\left(\sum_{n,l}\frac{1}{2\pi i} \left( \int_{\left(\frac{d}{2}\right)}P_n(s)\varphi(s)\widetilde{{(f_{\varepsilon_t})}_{n,l}}(s)y^{d-s}ds\right)\psi_{n,l}(k)\right)\cdot\\ & \qquad \qquad \qquad \qquad \qquad \qquad \qquad \qquad \cdot\left(\overline{\sum_{n',l'}{(\chi_t)}_{n',l'}(y)\psi_{n',l'}(k)} \right)y^{-(d+1)}dydk\\
        &= \int_{K/M}\frac{1}{2\pi i} \int_{\left(\frac{d}{2}\right)}\left( \sum_{n,l}P_n(s)\varphi(s)\widetilde{{(f_{\varepsilon_t})}_{n,l}}(s)\psi_{n,l}(k)\right)\cdot\\ & \quad \quad \qquad \qquad \qquad \qquad \qquad \cdot\left( \sum_{n',l'}\left( \int_0^{+\infty}\overline{{(\chi_t)}_{n',l'}(y)}y^{-(s+1)}dy\right)\overline{\psi_{n',l'}(k)}\right)dsdk\\
        &=\sum_{n,l}\frac{1}{2\pi i}\left(\int_{\left(\frac{d}{2}\right)}P_n(s)\varphi(s)\widetilde{{(f_{\varepsilon_t})}_{n,l}}(s)\widetilde{{(\chi_t)}_{n,l}}(\overline{s}) ds \right).
    \end{align*}
    
    We use again the same decomposition as in \eqref{decomposition in polar coordinates} $$\chi_t(ka_ye_0)=\phi_{t,y}(k)\rho_{t}(y)$$  and introduce the function
     $$F_{\varepsilon_t}(y,k):= f_{\varepsilon_t}(ka_ye_0).$$ Moreover, using the fact that there is at most one exceptional pole at $s_d=\lfloor \frac{d+2}{2}\rfloor$ in $(\frac{d}{2},d)$, we can move the contour of integration to the line $\left( r\right)$ for some $\frac{d}{2}<r<s_d$. By expanding the integrand similarly to \eqref{computation of M_F,f 1}, then using integration by parts for the $l$-th partial derivative with respect to $y_1$, with $l\geq 1$ large enough to be fixed later, we have 
\begin{align}
    &\langle \cT_{}(f_{\varepsilon_t}),\chi_t\rangle_{\mu_\mathcal{V}} \nonumber\\
    &=\frac{1}{2\pi i}\int_{\left(r\right)}\varphi(s)\sum_{n\geq0}P_n(s)\Bigg(\int_{\mathbb{R}_+}\int_{\mathbb{R}_+}\langle F_{\varepsilon_t}(y_1,\cdot ), {(\phi_{t,y_2})}_n\rangle_{\mu_{K/M}}\; \cdot\nonumber\\
&\qquad\qquad\qquad\qquad\qquad\qquad\qquad\qquad\qquad
\cdot \rho_t(y_2)y_1^{-(s+1)}y_2^{-(\bar s+1)}\, dy_1dy_2 \Bigg) ds\nonumber\\
    &=\frac{(-1)^{l}}{2\pi i}\int_{\left(r\right)}\varphi(s)\sum_{n\geq0}P_n(s)
    \Bigg(\int_{\mathbb{R}_+}\int_{\mathbb{R}_+}\langle \frac{\partial^l F_{\varepsilon_t}(y_1,\cdot ) }{\partial^l y_1}, {(\phi_{t,y_2})}_n\rangle_{\mu_{K/M}}\; \cdot\nonumber\\
&\qquad\qquad\qquad\qquad\qquad\qquad\qquad\qquad
\cdot\rho_t(y_2)\frac{y_1^{-(s+1)+l}}{\prod_{j=0}^{l-1}(s+j)}y_2^{-(\bar s+1)}\, dy_1dy_2 \Bigg) ds.\label{integral Tf,chi}
\end{align}
Using the orthogonality of the projections pr$_n$ and Lemma \ref{estimate P(s)-weighted correlations in K/M} with \eqref{integral Tf,chi}, it follows, for any $0<\delta<2r-d$ and any $r\in (\frac{d}{2},s_d)$, with $s=r+it\in\mathbb{C}$,
\begin{align*}
    &\left|\langle \cT_{}(f_{\varepsilon_t}),\chi \rangle_{\mu_\mathcal{V}}\right| \\
&\ll_l \int_{\left(r\right)}\left|\varphi(s)s^{-l}\right|
    \Bigg(\int_{\mathbb{R}_+}\int_{\mathbb{R}_+}\left|\sum_{n\geq0}P_n(s)\langle \left(\frac{\partial^l F_{\varepsilon_t}(y_1,\cdot ) }{\partial^l y_1}\right)_n, {(\phi_{t,y_2})}_n\rangle_{\mu_{K/M}}\right|\; \cdot\nonumber\\
&\qquad\qquad\qquad\qquad\qquad\qquad\qquad\qquad
\cdot\rho_t(y_2)y_1^{-(r+1)+l}y_2^{-( r+1)}\, dy_1dy_2 \Bigg) ds\\
     &\ll_\delta\int_{\left(r\right)}\left|\varphi(s)s^{-l}\right||t|^{1/2}ds~\left(\int_{\mathbb{R}_+}\int_{\mathbb{R}_+}\left\lVert \frac{\partial^l F_{\varepsilon_t}(y_1,\cdot ) }{\partial^l y_1}\right\rVert_{L^2_{K/M}}\left\lVert\phi_{t,y_2}\right\rVert_{L^1_{K/M}}^{\frac{1}{2}+\frac{2r-d-\delta}{d+1}}\;\cdot \right. \\
&\qquad\qquad\qquad\qquad\qquad\qquad\qquad\qquad\qquad\cdot~ \rho_t(y_2)y_1^{-(r+1)+l}y_2^{-(r+1)} dy_1dy_2 \Bigg).
\end{align*}
We use further that the projection $I:=\text{pr}_{\mathbb R_+}$(supp$(F_{\varepsilon_t}))$ is uniformly bounded away from $y=0$, and, similarly to \eqref{estimate phi_y}, that $ \lVert \phi_{t,y_2}\rVert_{L^1_{K/M}}\ll y_2^d $. It follows,
\begin{align*}
&\int_{\mathbb{R}_+}\int_{\mathbb{R}_+}\left\lVert \frac{\partial^l F_{\varepsilon_t}(y_1,\cdot ) }{\partial^l y_1}\right\rVert_{L^2_{K/M}}\left\lVert\phi_{t,y_2}\right\rVert_{L^1_{K/M}}^{\frac{1}{2}+\frac{2r-d-\delta}{d+1}}\;\cdot  \rho_t(y_2)y_1^{-(r+1)+l}y_2^{-(r+1)} dy_1dy_2\\
     &\ll \left\lVert f_{\varepsilon_t} \right\rVert_{C^l}\left(\int_I y_1^{-(r+1)+l}dy_1\right)\left(\int_{\mathbb{R}_+}y_2^{d(\frac{1}{2}+\frac{2r-d-\delta}{d+1})}\; \rho_t(y_2)y_2^{-(r+1)} dy_2\right)\\
     &\ll_r\left\lVert f_{\varepsilon_t} \right\rVert_{C^l}\left(\int_{\mathbb{R}_+}\; \rho_t(y_2)y_2^{d(\frac{1}{2}+\frac{2r-d-\delta}{d+1})-r-1} dy_2\right)\\
    &\ll \varepsilon_t^{-l}~ ~\int_{\alpha e^{-t}}^{\beta e^{-t}}\; y_2^{d(\frac{1}{2}+\frac{2r-d-\delta}{d+1})-r-1} dy_2\\
    &\ll\varepsilon_t^{-l}e^{-t\left(d(\frac{1}{2}+\frac{2r-d-\delta}{d+1})-r\right)} \;.
\end{align*}
We write $\sigma:=d(\frac{1}{2}+\frac{2r-d-\delta}{d+1})-r $ and choose $0<\delta<\frac{(2r-d)(d-1)}{2d}$ such that $\sigma >0$. We use further the following estimates for the scattering matrix $ \varphi(s)$ near the critical line $(\frac{d}{2})$ in terms of a function $W(t)\geq 1$, with $s=r+it$, introduced in \cite[Propositions 7.11 and 7.13]{sarnak1980}. We have
\begin{align}
    &|\varphi(s)|^2=1+O\left( (r-\frac{d}{2})W(t)\right)~,~ \text{for} ~ \frac{d}{2}<r<r_0<d,\label{estimate phi(s)}\\
\qand \qquad &\int_0^T W(t) dt\ll T^{d+1}~.\label{estimate W(t)}
\end{align}
For fixed $r\in(\frac{d}{2},s_d)$ and $l\geq 1$ large enough, using Cauchy-Schwarz inequality then the estimate \eqref{estimate phi(s)}, we have
\begin{align*}
&\left(\int_{\left(r\right)}\left|\varphi(s)s^{-l+\frac{1}{2}} \right|ds \right)^2  \leq  \left(\int_{\left(r\right)}\left|\varphi(s)s^{-\frac{l}{2}} \right|^2ds \right)\left(\int_{\left(r\right)}\left|s^{-l+1} \right|ds \right)\\ 
 &= \left(\int_{\mathbb R}\left(1+O\left(rW(t)\right)\right)|t^2+r^2|^{-\frac{l}{2}}dt\right)\left(\int_{\mathbb R}|t^2+r^2|^{\frac{-l+1}{2}} dt \right) \\
 &\ll_{r,l}\int_{\mathbb R}W(t)|t^2+r^2|^{-\frac{l}{2}}dt \\
&\ll \int_{\mathbb R} |t|^{d+2}|t^2+r^2|^{-\frac{l+2}{2}}dt , \quad \text{ by integration by parts, with the estimate\eqref{estimate W(t)}},\\
&= O(1),\quad \text{ for some fixed }l> d+1~.
\end{align*}
Altogether, we obtain
$$\left| \langle \cT_{}(\chi),\chi_t\rangle\right| \ll \varepsilon_t^{1/2} + \varepsilon_t^{-l}e^{-\sigma t},
$$
and choose $\varepsilon_t = e^{-2\gamma t}$ with $\gamma :=\frac{\sigma}{1+2l} $.
\end{proof}
Putting together \eqref{eq:final estimate variance F_N}, \eqref{eq:reduced form variance} and the estimates showed in Lemmas \ref{convergence of M_f,f^K} and \ref{convergence of <tau(X_t),X>}, we obtain

\begin{align*}
    \left\lVert \mathsf{F}_N^{(\varepsilon,L)} \right\rVert_{L^2_\mathcal{Y}}^2 &= \sum_{t=-K+1}^{K-1}\Theta_\infty(t) + o(1)\\
    &= \sum_{t=-K+1}^{K-1} \left(\sum_{i,j=1}^{m} \langle\cT_{ij}(\chi_{t,\lambda_i}),\chi_{\lambda_j}\rangle_{\mu_\mathcal{V}}+c_Q M_{\chi_t,\chi}(s_d)\right) + o(1)
    \end{align*}
and 
\begin{equation}
\label{eq: explicite variance}
    \left\lVert \mathsf{F}_N^{(\varepsilon,L)} \right\rVert_{L^2_\mathcal{Y}}^2\stackrel{N\rightarrow \infty}{\longrightarrow} \; \sigma^2 \;:=\;\sum_{i,j=1}^{m} \langle\cT_{ij}(\chi_{\infty,\lambda_i}),\chi_{\lambda_j}\rangle_{\mu_\mathcal{V}}+c_Q M_{\chi_\infty,\chi}(s_d),
\end{equation}
where we denote by $\chi_\infty$ the characteristic function of the domain
$$
 \{ x \in \mathcal{V} : x_{d+2}^{2}- x_{d+1}^{2} < c^{2} \} \; = \; \bigcup_{t=-\infty}^{\infty}a_{-t}\left(F_{1,c}\right).
$$

\subsection{Proof of the effective estimate for the counting function}	
\label{sec:proof of the second main theorem}
To obtain an effective estimate for the counting function $\mathsf{N}_{T,c}$, the central argument in our approach is to derive an almost-everywhere bound for averages $\sum_{t=0}^{T-1}(\widehat\chi\circ a_t-\mu_\mathcal{X}(\widehat\chi))$ from an $L^2$-bound on these averages. We generalized this argument in \cite{ouaggag2022effective} to $L^p$-bounds, $p>1$, following the approach in \cite{kleinbock2017pointwise} based on an original idea of Schmidt in \cite{Schmidt1960AMT}. We generalize in the following proposition our result from \cite[Proposition 4.2.]{ouaggag2022effective} in order to take into account the approximation of $\chi_{1,c}$ by the sequence $(\chi_{1,c,t})_{t\geq 0}$ introduced by the "sandwiching" \eqref{resandwiching N_T,c}.

\begin{prop}
\label{pointwise equidistribution}
Let $(Y,\nu)$ be a probability space, and let $(f_t)_{t\geq0}$ be a sequence of measurable functions $f_t: Y \rightarrow \mathbb{R}$. Suppose there exist constants $ p> 1$ and $C>0$ such that for any integers $0\leq a<b$,
\begin{equation}\label{Lp bound}
     \int_{Y}\left|  \sum_{t=a}^{b-1}f_t(y) \right|^{p}d\nu(y) \leq C(b-a)~.
 \end{equation}

Then, there exists $C_p>0$, depending only on $p>1$, such that for every $\varepsilon >0$ and $\nu$-almost every $y \in Y$, we have
\begin{equation}\label{a.e. bound}
     \sum_{t=0}^{N-1} f_t(y) \leq C_p
 N^{\frac{1}{p}} \log ^{1+\frac{1}{p}+ \frac{\varepsilon}{p}} N.
 \end{equation}
\end{prop}

\begin{remark}
\label{rem:dyadic subsets}\emph{In the argument as formulated in \cite{kleinbock2017pointwise}, the estimate in \eqref{Lp bound} is satisfied for all pairs $(a,b)$ of the form $(2^i j,2^i(j+1))$ building the (non-disjoint) dyadic decomposition of $N-1$. In fact, it is enough for our argument (see proof of Lemma \ref{lem:sets for Borel Cantelli lemma}) to consider only a smaller selection of such pairs, denoted below by $L(N)$, which builds a partition of $\{1,\dots,N-1\}$ and allows moreover to satisfy the conditions $t\geq \kappa \log L$ and $t\geq -\frac{1}{\theta} \log \varepsilon$ required by Proposition \ref{bounds truncated siegel transform} and Proposition \ref{smooth approximation}, for $t \in \{a,\dots,b-1\}$ and for parameters $L$ and $\varepsilon$ to be defined later as functions of $(a,b)$.\\
 For non-negative integers $a<b$ we write $[a..b)\coloneqq [a,b)\cap \mathbb{N}$. For an integer $s\geq 1$ we consider the following set of dyadic subsets,
\begin{align*}
L_s\coloneqq &\Big\{ \left[2^i..2^{i+1}\right)~:~ 0\leq i\leq s-2 \Big\}\; \cup \\
&\qquad\qquad\Big\{ \left[2^ij..2^{i}(j+1)\right)~:~ 2^ij\geq 2^{s-1}, 2^i(j+1)\leq 2^s \Big\} \;\cup\; \Big\{ 0 \Big\},
\end{align*}
where the sets of the first type $[2^i..2^{i+1})$, $0\leq i\leq s-2$, together with $[0..1)$, are a decomposition of the set $[0..2^{s-1})$.\\
We observe that for any integer $N\geq 2$ with $2^{s-1}< N \leq 2^s$, the set $[0..N)$ is the disjoint union of at most $2s-1$ subsets in $L_s$, namely $[0..1)$, the $s-1$ subsets of the first type and at most $s-1$ sets of the second type which can be constructed from the binary expansion of $N-1$. We denote by $L(N)$ this smaller selection of subsets from $L_s$, i.e. }
$$L(N)\subset L_s\;, \quad \emph{\text{with }} \quad |L(N)|\leq 2s-1 \quad \emph{\text{ and }}\quad [0..N)=\bigsqcup_{I\in L(N)} I.$$
\end{remark}

In the following lemmas, the notations and assumptions are the same as in Proposition \ref{pointwise equidistribution}~.

\begin{lem} For any $s\geq 1$, we have
$$\sum_{I\in L_s}\int_{Y}\left|\sum_{t\in I} f_t(y) \right|^p d\nu(y) \leq Cs2^{s}. 
$$
\label{dyadic estimate}
\end{lem}
\begin{proof} Since $L_s$ is a subset of the set of all dyadic sets  $[2^ij..2^i(j+1))$ where $i,j$ are non-negative integers and $2^i(j+1)\leq 2^s$, we have 
\begin{align*}
\sum_{I\in L_s}\int_{Y}\left|\sum_{t\in I} f_t(y) \right|^p d\nu(y)  &\leq \sum_{i=0}^{s-1}\sum_{j=0}^{2^{s-i}-1} \int_{Y}\left|\sum_{t\in I} f_t(y) \right|^p d\nu(y) \\
&\leq \sum_{i=0}^{s-1}\sum_{j=0}^{2^{s-i}-1} C2^i\\
&\leq Cs2^s.
\end{align*}
\end{proof}

\begin{lem}
For every $\varepsilon >0$, there exists a sequence of measurable subsets $\left( Y_{s}\right)_{s \in \mathbb{N}}$ of $Y$ such that:
\begin{enumerate}
\item $\nu(Y_s) \leq Cs^{-(1+p\varepsilon)}.$
\item For every integer $N\geq 2$ with $2^{s-1}\leq N-1<2^s$ and every $y \notin Y_s$ one has
\begin{equation}
\left| \sum_{t=0}^{N-1}f_t(y) \right| \ll_p 
2^{\frac{s}{p}}s^{1+\frac{1}{p}+\varepsilon}.
\end{equation}
\end{enumerate}
\label{lem:sets for Borel Cantelli lemma}
\end{lem}

\begin{proof}
For every $s\geq1 $, consider 
\begin{equation}\label{eq: definition sets Y_s}
Y_{s}=
\left\lbrace y \in Y~:~ \sum_{I\in L_s} \left| \sum_{t\in I}f_t(y)\right|^p ~>~2^{s} s^{2+ p\varepsilon} \right\rbrace.
\end{equation}

The first assertion follows from Lemma \ref{dyadic estimate} and Markov's Inequality.\\
Further, for any $N\geq 2$ such that $2^{s-1}\leq N-1<2^s$ and any $y\notin Y_s$, using the partition $[0..N)=\bigsqcup_{I\in L(N)}I$ with $L(N)$ of cardinality at most $2s-1$, we have

\begin{align*}
\left|\sum_{t=0}^{N-1} f_t(y)\right|^p &= \left|\sum_{I\in L(N)}\sum_{t\in I} f_t(y)\right|^p\\
&\leq (2s-1)^{p-1}\sum_{I\in L(N)}\left|\sum_{t\in I} f_t(y)\right|^p &\text{(by Hölder's Inequality)}\\
&\leq (2s-1)^{p-1}\sum_{I\in L_s}\left|\sum_{t\in I} f_t(y)\right|^p \\
&\ll_p s^{1+p+p\varepsilon}2^{s} &\text{(since }y\notin Y_{s})
\end{align*}
which yields the claim by raising to the power $\frac{1}{p}$.
\end{proof}
\begin{proof}[Proof of Proposition \ref{pointwise equidistribution}]
Let $\varepsilon>0$ and choose a sequence of measurable subsets $(Y_s)_{s\in \mathbb{N}}$ as in \eqref{eq: definition sets Y_s}. Observe that
$$\sum_{s=1}^\infty \nu(Y_s) \leq \sum_{s=1}^\infty C s^{-(1+p\varepsilon)} < \infty.
$$
The Borel-Cantelli lemma implies that there exists a full-measure subset $Y(\varepsilon)\subset Y$ such that for every $y \in Y(\varepsilon) $ there exists $s_y \in \mathbb{N}$ such that for all $s > s_y$ we have $y \notin Y_s$.\\
Let $N\geq 2$ and $s= 1 + \left\lfloor \log (N-1) \right\rfloor$, so that $2^{s-1}\leq N-1 < 2^s $. Then, for $N-1 \geq 2^{s_y}$ we have $s> s_y$ and $y \notin Y_s$, thus
\begin{align*}
\left|\sum_{t=0}^{N-1}f_t(y) \right| &\ll_p 2^{\frac{s}{p}}s^{1+\frac{1}{p} +\varepsilon}
\\
&\leq
(2N)^{\frac{1}{p}}\log^{1+\frac{1}{p}+ \varepsilon}(2N).
\end{align*}
This implies the claim for $y \in \cap_{m\in \mathbb{N}}Y(1/m)$.
\end{proof}

We now apply Proposition \ref{pointwise equidistribution} to the sequence of functions defined for $t\geq 0$ by $$f_{t}=\widehat{\chi_t}\circ a_t-\int_{\mathcal{X}}\widehat{\chi_t},$$ where $\chi_t=\chi_{_{F_{1,c,t}}}$ is the characteristic function of the set $F_{1,c,t}$ defined in \eqref{eq:defintion F_1,c_l}. We also recall that, by Proposition \ref{siegel mean value thrm}, we have
$$\text{vol}(F_{1,c,t}) \coloneqq \int_{\mathcal{V}} \chi_t(z)dz = \int_{\mathcal{X}} \widehat{\chi_t}(\Lambda)d\mu_\mathcal{X}(\Lambda).$$

\begin{prop}
\label{pointwise equidistribution counting function}
Let $d\geq 3$. For all $\varepsilon>0$, for almost every $k \in K$ we have
$$
 \sum_{t=0}^{N-1} \widehat{\chi_t}(a_t k \Lambda_0) = \sum_{t=0}^{N-1}\emph{vol}(F_{1,c,t}) + O_{k,\varepsilon} (N^{\frac{1}{2}+\varepsilon}).
$$
\end{prop}
\begin{proof}
Using Proposition \ref{pointwise equidistribution}, it is enough to show that for every $1<p<2$, for all pairs of integers $(a,b)$ from the dyadic decomposition of $N$ specified in Remark \ref{rem:dyadic subsets}, we have
\begin{equation}
\label{L2-delta estimate}
\left|\left| \sum_{t=a}^{b-1} \left( \widehat{\chi_t}-\mu_{\mathcal{X}}(\widehat{\chi_t})\right)\circ a_t \right|\right|_{L^{p}(\mathcal{Y})}^{p} \ll (b-a)~.
\end{equation}
We first observe that in all the estimates used in this proof, the dependence on the test functions $\chi_t$ involves only the support of the functions, and the estimates are therefore uniform in $t>0$.\\
Let $1<p<2$. Using the estimates for the truncated Siegel transform from Proposition \ref{bounds truncated siegel transform}, we have, for $\frac{2p}{3p-2}\leq \tau < d$ and $t\geq \kappa \log L$,
\begin{align}
&\left|\left|  \left( \widehat{\chi_t}\circ a_t-\mu_{\mathcal{X}}(\widehat{\chi_t})\right)- \left( \widehat{\chi_t}^{(L)}\circ a_t-\mu_{\mathcal{X}}(\widehat{\chi_t}^{(L)})\right)\right|\right|_{L^{p}_\mathcal{Y}} \\
&\qquad\qquad\qquad\qquad\qquad\leq \left|\left|   \widehat{\chi_t}\circ a_t- \widehat{\chi_t}^{(L)}\circ a_t\right|\right|_{L^{p}_\mathcal{Y}}+ \int_{\mathcal{X}}\left|\widehat{\chi_t}-\widehat{\chi_t}^{(L)} \right| \nonumber\\
& \qquad\qquad\qquad\qquad\qquad\ll L^{-\frac{\tau(2-p)}{2p}}~+~L^{-(\tau-1)} \nonumber\\
& \qquad\qquad\qquad\qquad\qquad\ll L^{-\frac{\tau(2-p)}{2p}}.
\label{estimate truncation}
\end{align}
Further, using Proposition \ref{smooth approximation} and the estimates from Proposition \ref{bounds truncated siegel transform} and (\ref{epsilon approximation of chi}), we have for $t\geq -\frac{1}{\theta}\log \varepsilon$,
\begin{align}
&\left\lVert  \left( \widehat{\chi_t}^{(L)}\circ a_t-\mu_{\mathcal{X}}(\widehat{\chi_t}^{(L)})\right) - \left( \widehat{f}_{t,\varepsilon}^{(L)}\circ a_t-\mu_{\mathcal{X}}(\widehat{f}_{t,\varepsilon}^{(L)})\right)\right\rVert_{L^{p}_\mathcal{Y}}  \\
&\qquad \qquad \qquad\leq \left|\left|   \widehat{\chi_t}^{(L)}\circ a_t- \widehat{f}_{t,\varepsilon}^{(L)}\circ a_t\right|\right|_{L^{p}_\mathcal{Y}}+ \int_{\mathcal{X}}\left|\widehat{\chi_t}^{(L)}-\widehat{f}_{t,\varepsilon}^{(L)} \right| \nonumber\\
&\qquad \qquad \qquad\leq \left|\left|(\widehat{\chi_t -f_{t,\varepsilon}})^{(L)}\circ a_t \right|\right|_{\infty}^{\frac{p-1}{p}}\cdot \left|\left|(\widehat{\chi_t -f_{t,\varepsilon}})^{(L)}\circ a_t \right|\right|_{L^1_\mathcal{Y}}^{\frac{1}{p}}+ \int_{\mathcal{V}}\left|\chi_t-f_{t,\varepsilon} \right| \nonumber \\
&\qquad \qquad \qquad \ll L^{\frac{p-1}{p}}~\varepsilon^{\frac{1}{p}}~ +~ \varepsilon \nonumber \\
& \qquad \qquad \qquad\ll L^{\frac{p-1}{p}}~\varepsilon^{\frac{1}{p}}.
\label{estimate smooth}
\end{align}
Further, using effective equidistribution for smooth and compactly supported functions (Theorem \ref{double equidistribution K orbits} for $r=1$), we have 
\begin{align}
&\left\lVert  \sum_{t=a}^{b-1}\left( \widehat{f}_{t,\varepsilon}^{(L)}-\mu_{\mathcal{X}}(\widehat{f}_{t,\varepsilon}^{(L)})\right)\circ a_t \right\rVert_{L^{p}_\mathcal{Y}} \leq \left\lVert  \sum_{t=a}^{b-1}\left( \widehat{f}_{t,\varepsilon}^{(L)}-\mu_{\mathcal{X}}(\widehat{f}_{t,\varepsilon}^{(L)})\right)\circ a_t \right\rVert_{L^{2}_\mathcal{Y}} \nonumber \\
&\qquad\leq \left\lVert  \sum_{t=a}^{b-1}\left( \widehat{f}_{t,\varepsilon}^{(L)}\circ a_t-\mu_{\mathcal{Y}}(\widehat{f}_{t,\varepsilon}^{(L)}\circ a_t)\right) \right\rVert_{L^{2}_\mathcal{Y}} + \left\lvert  \sum_{t=a}^{b-1}\left( \mu_{\mathcal{Y}}(\widehat{f}_{t,\varepsilon}^{(L)}\circ a_t)-\mu_{\mathcal{X}}(\widehat{f}_{t,\varepsilon}^{(L)})\right) \right\rvert \nonumber \\
&\qquad\ll (b-a)^{\frac{1}{2}}\left\lVert \tilde{\mathsf{F}}_{b,a}^{(\varepsilon,L)} \right\rVert_{L^2_\mathcal{Y}} + \sum_{t=a}^{b-1} e^{-\delta t} \left\lVert \widehat{f}_{t,\varepsilon}^{(L)}\right\rVert_l,
\label{estimate equidistribution}
\end{align}
where we define similarly to \eqref{eq:def F_N,M} the averages
$$\tilde{\mathsf{F}}_{b,a}^{(\varepsilon, L)} = \frac{1}{\sqrt{b-a}} \sum_{t=a}^{b-1}\left( \widehat{f}_{t,\varepsilon}^{(L)} \circ a_t -\mu_\mathcal{Y}(\widehat{f}_{t,\varepsilon}^{(L)}\circ a_t) \right).$$
We stress again, that the dependence on the test functions $f_{\varepsilon}$ in the estimate \eqref{eq:starrrr} only involves the support of $f_\varepsilon$, and therefore, replacing $f_\varepsilon$ by $f_{t,\varepsilon}$ yields a similar estimate. We have 
\begin{align}\label{estimate F_ba}
    \left\lVert  \tilde{\mathsf{F}}_{b,a}^{(\varepsilon, L)}\right\rVert_{L^2_\mathcal{Y}} &= \sum_{-K}^{K}\Theta_\infty(t)\nonumber\\
    &+ O\left( (b-a)^{-1}K^2+\left((b-a)^{-1}e^{-\delta a}e^{\xi K}+e^{-\delta K}\right)L^2\varepsilon^{-2l}+KL^{-\frac{\tau-2}{2}}+K\varepsilon\right).
\end{align}
Putting \eqref{estimate equidistribution}, \eqref{estimate truncation}, \eqref{estimate smooth} and \eqref{estimate F_ba} together, and considering moreover that $\sum_{-K}^{K}\Theta_\infty(t)$ is bounded uniformly in $K$ (by the convergence shown in Section \ref{subsec:estimating variance}), we obtain, for $a\geq 1$ large enough such that the conditions $t\geq \kappa \log L$ and $t\geq -\frac{1}{\theta}\log \varepsilon$ are satisfied (for parameters $L$ and $\varepsilon$ to be specified later),
\begin{align}
    &\left|\left| \sum_{t=a}^{b-1} \left( \widehat{\chi_t}-\int_{\mathcal{X}}\widehat{\chi_t}\right)\circ a_t \right|\right|_{L^{p}(\mathcal{Y})}\nonumber\\
    &\ll (b-a)\left(L^{-\frac{\tau(2-p)}{2p}}+L^{\frac{p-1}{p}}\varepsilon^{\frac{1}{p}}\right)\;+\;e^{-\delta a}L\varepsilon^{-l} \label{first term in all estimates L_p together}\\
    &+(b-a)^{\frac{1}{2}}\left( 1+(b-a)^{-1}K^2+\left((b-a)^{-1}e^{-\delta a}e^{\xi K}+e^{-\delta K}\right)L^2\varepsilon^{-2l}+KL^{-\frac{\tau-2}{2}}+K\varepsilon \right)\label{all estimates L_p together}.
\end{align}
In order to bound the first term in \eqref{first term in all estimates L_p together} by $ (b-a)^{\frac{1}{p}}$, we choose 
\begin{align*}
    &L=(b-a)^{q_2} ,\;\text{ with }\;q_2= \frac{2(p-1)}{\tau (2-p)},\\
  \text{and}\qquad &\varepsilon= L^{-\frac{\tau(2-p)+2(p-1)}{2}}= (b-a)^{-q_1},\\
  \text{ with }\quad &q_1=q_2\frac{\tau(2-p)+2(p-1)}{2}= (p-1)(1+q_2).
  \end{align*}
  We choose further
  \begin{align*}
         \quad &K=c_p \log (b-a),\;\text{ with some }\;  c_p>0,
\end{align*}
and, in order to satisfy the condition $K\leq a$, we verify that for all but finitely many pairs $(a,b)$ from the dyadic decomposition of $N$, i.e. for pairs of the forms $(2^i,2^{i+1})$ and $(2^ij,2^i(j+1)) $ with $i\geq i_0(c_p)\geq 1$, we have
\begin{equation*}
    K=c_p \log(2^i) \leq a.
\end{equation*}
With this choice of $L$, $\varepsilon$ and $K$, we verify next that the terms in \eqref{all estimates L_p together} are bounded by $(b-a)^{\frac{1}{p}}$, for all but finitely many pairs $(a,b)$ from the dyadic decomposition. We have indeed
\begin{equation*}
    (b-a)^{-\frac{1}{2}}e^{-\delta a}e^{\xi K}L^2\varepsilon^{-2l} \leq e^{-\delta 2^{i_0}}(b-a)^{-\frac{1}{2}+\xi c_p+2q_2+2lq_1}\leq  (b-a)^{\frac{1}{p}},
\end{equation*}
and for some $c_p>0$ large enough, we also have
\begin{equation*}
    (b-a)^{\frac{1}{2}}e^{-\delta K}L^2\varepsilon^{-2l} \leq (b-a)^{\frac{1}{2}-\delta c_p+2q_2+2lq_1}\leq  (b-a)^{\frac{1}{p}}.
\end{equation*}
One verifies easily that the other terms in \eqref{all estimates L_p together} are also bounded by $(b-a)^\frac{1}{p}$.\\
Finally, we verify that the conditions $t\geq \kappa \log  L$ and $t\geq -\frac{1}{\theta}\log \varepsilon$ are also satisfied, since we have, for all but finitely many pairs $(a,b)$,
\begin{align*}
    &\kappa \log L = \kappa q_2 \log (b-a) = \kappa q_2 \log (2^i) \leq a \leq t\\
    \text{and} \quad &-\frac{1}{\theta}\log \varepsilon = \frac{q_1}{\theta} \log (b-a) = \frac{q_1}{\theta} \log (2^i) \leq a \leq t.
\end{align*}
We obtain, for $d>\tau>2$,
$$
\left|\left| \sum_{t=a}^{b-1} \left( \widehat{\chi_t}-\int_{\mathcal{X}}\widehat{\chi_t}\right)\circ a_t \right|\right|_{L^{p}_\mathcal{Y}} \ll_p (b-a)^{\frac{1}{p}},
$$
which ends the proof.
\end{proof}
 Using the same argument as in the proof of Proposition \ref{pointwise equidistribution counting function}, with $\chi$ instead of $\chi_t$ and with the estimate from Proposition \ref{smooth approximation for c} instead of Proposition \ref{smooth approximation}, yields the same asymptotic for the ergodic averages $ \sum_t \widehat\chi(a_tk\Lambda_0)$. 
\begin{prop}
\label{pointwise equidistribution counting function l}
Let $d\geq 3$. For almost every $k \in K$ and all $\varepsilon>0$, we have
\begin{equation*}
  \sum_{t=0}^{N-1} \widehat{\chi}(a_t k \Lambda_0) = N\,\emph{vol}(F_{1,c}) + O_{k,\varepsilon} (N^{\frac{1}{2}+\varepsilon}).  
\end{equation*}
\end{prop}

\paragraph{\textbf{Proof of the first main theorem:}}
\begin{proof}[Proof of Theorem \ref{thm: first result}]
\label{proof of main theorem}
We recall the estimate from (\ref{resandwiching N_T,c}). We have
\begin{align*}
\sum_{t=0}^{\left\lfloor T\right\rfloor-r_0 -1} \widehat{\chi}_{1,c,t}(a_t k\Lambda_0) +O\left( 1\right) ~\leq ~\mathsf{N}_{T,c}(\alpha_k) +O(1) ~\leq~ \sum_{t=0}^{\left\lceil T\right\rceil+r_0 -1} \widehat{\chi}_{1,c}(a_t k\Lambda_0),
\end{align*}
which gives, together with Propositions \ref{pointwise equidistribution counting function} and \ref{pointwise equidistribution counting function l}, the following estimate, for almost every $k\in K$, for all $T> T_k$ for some $T_k\geq 1$ large enough, and for all $\varepsilon>0$. We have
\begin{equation*}
\sum_{t=0}^{\left\lfloor T \right\rfloor -r_0-1} \text{vol}(F_{1,c,t})+O\left(T ^{\frac{1}{2}+\varepsilon} \right)\leq  \mathsf{N}_{T,c}(\alpha_k)  \leq T\text{vol}(F_{1,c})+O\left(T ^{\frac{1}{2}+\varepsilon}\right).
\end{equation*}
From (\ref{volum F_1,c approximation}), we have
\begin{align*}
\sum_{t=0}^{\left\lfloor T \right\rfloor-r_0-1} \text{vol}(F_{1,c,t})&= \sum_{t=0}^{\left\lfloor T \right\rfloor-r_0-1} \left(\text{vol}(F_{1,c})+O(t^{-1})\right)\\
&\geq \; T \text{vol}(F_{1,c}) +O\left(\log T\right).
\end{align*}
Altogether, we obtain
\begin{equation*}
\mathsf{N}_{T,c}(\alpha_k)
=T\text{vol}(F_{1,c})+O\left(T ^{\frac{1}{2}+\varepsilon}\right).
\end{equation*}
Since full-measure sets in $K$ correspond to full-measure sets in $S^d$, we conclude that this last estimate holds for almost every $\alpha \in$ $S^d$.
\end{proof}

\newpage

\section{CLT for counting Diophantine approximations on spheres}\label{sec:clt counting}

We carry on in this section our analysis of the counting function of rational approximations on the sphere, defined, for $T>0$, $c>0$, $\alpha \in S^{d}$ and $d\geq 3$, by 
\begin{equation*}
    \mathsf{N}_{T,c}(\alpha) \coloneqq  \lvert \lbrace (p,q) \in \mathbb{Z}^{d+1}\times \mathbb{N} : \frac{p}{q} \in S^d,~1 \leq q < \cosh T \text{ and } (\ref{intrinsic Diophantine 2}) \text{ holds } \rbrace  \rvert~,
\end{equation*}
for the intrinsic Diophantine approximation problem (with the critical Dirichlet exponent for intrinsic Diophantine approximation on $S^{d}$), given by
\begin{equation}
\label{intrinsic Diophantine 2}
 \left\lVert \alpha - \frac{p}{q}\right\rVert~ < ~  \frac{c}{q} \;.
\end{equation}


We showed in the previous section, using an analog of the classical Dani correspondence for the space of orthogonal lattices $\mathcal{X}=\SO(n+1,1)/\SO_{\mathbb{Z}}(n+1,n)$, that the counting function $\mathsf{N}_{T,c}$ can be related to the light-cone Siegel transform of the characteristic function $\chi$ of an elementary domain $F_{1,c}\subset \mathcal{V}$, using averages of the form
\begin{equation*}
    \mathsf{N}_{T,c}(\alpha) \approx \sum_{t=0}^{T-1} \widehat\chi \circ a_t (k_\alpha \Lambda_0),
\end{equation*}
i.e. ergodic averages of the light-cone Siegel transform $\widehat{\chi}$ along compact orbits $k_\alpha \Lambda_0$ pushed by the hyperbolic subgroup $\{a_t\}$. We carried out the analysis of these averages and, using equidistribution and non-divergence results on the space $\mathcal{X}$, showed that we can approximate the test function $\chi_t$ by a family of smooth and compactly supported approximations $f_{t,\varepsilon}^{(L)}$ and still control the approximation of $\mathsf{N}_{T,c}$, i.e. 
\begin{equation}\label{smooth truncated approx}
    \mathsf{N}_{T,c}(\alpha) \approx \sum_{t=0}^{T-1} \widehat f_{t,\varepsilon}^{(L)} \circ a_t (k_\alpha \Lambda_0).
\end{equation}
We investigate the limit distribution of $N_{T.c}$ with respect to the probability measure on $S^d$, using again the approximation \eqref{smooth truncated approx} to reduce the problem to the analysis of the limit distribution of ergodic averages for truncated and smooth approximations of the Siegel transform of $\widehat\chi$. We prove our second main result concerning the counting function $N_{T.c}$.

\begin{thm}[] 
\label{thm: second result}
Let $d\geq 3$. Then, for every $\xi \in \mathbb{R}$, 
\begin{equation}
\mu_{_{S^d}}\left( \left\lbrace \alpha \in S^d:\frac{\mathsf{N}_{T,c}(\alpha)-C_{c,d}\cdot T}{T^{1/2}}<\xi \right\rbrace\right) \rightarrow \emph{Norm}_\sigma(\xi)~,\quad \text{as }T\rightarrow \infty,
\end{equation}
where $\emph{Norm}_\sigma$ denotes the normal distribution with variance $\sigma^2\geq 0$.
\end{thm}	

Before developing this analysis, we verify in what follows that this approximation still gives the same limit distribution as $\mathsf{N}_{T,c}$.

\subsection{Approximation of the averaging function}\label{approximations of F_N}
In order to study the distribution of the counting function $\mathsf{N}_{T,c}$, we introduce the following averaging function.
\begin{equation}
\label{averages counting function}
\mathsf{F}_N \coloneqq \frac{1}{\sqrt{N}} \sum_{t=0}^{N-1}\left( \widehat{\chi_t} \circ a_t -\mu_\mathcal{Y}(\widehat{\chi_t}\circ a_t) \right)~, \text{ for } N\geq 1.
\end{equation}

We shall use in the following arguments the basic observation that if we approximate $\mathsf{F}_N$ by a sequence $\widetilde{\mathsf{F}}_N$ in such a way that $||\mathsf{F}_N-\widetilde{\mathsf{F}}_N ||_{L^1_\mathcal{Y}} \xrightarrow{N\rightarrow \infty} 0$ and the limit distribution of $\widetilde{\mathsf{F}}_N$ is continuous, then $\mathsf{F}_N$ and $\widetilde{\mathsf{F}}_N$ have the same convergence in distribution.

\subsubsection{Truncated averages}
We first observe that $\mathsf{F}_N$ has the same convergence in distribution as the truncated averages
\begin{equation}
\mathsf{F}_{N,M} \coloneqq \frac{1}{\sqrt{N-M}} \sum_{t=M}^{N-1}\left( \widehat{\chi_t} \circ a_t -\mu_\mathcal{Y}(\widehat{\chi_t}\circ a_t) \right),
\end{equation}
for some $M=M(N)\rightarrow \infty$ to be specified later.\\
Indeed, we have
\begin{align*}
    \left\lVert \mathsf{F}_N -\mathsf{F}_{N,M}  \right\rVert_{L^1_\mathcal{Y}} &\leq \frac{1}{\sqrt{N}}\sum_{t=0}^{M-1}\left\lVert \widehat{\chi_t} \circ a_t -\mu_\mathcal{Y}(\widehat{\chi_t}\circ a_t) \right\rVert_{L^1_\mathcal{Y}} \\
    &\qquad\qquad+ \left( \frac{1}{\sqrt{N-M}}-\frac{1}{\sqrt{N}}\right) \sum_{t=M}^{N-1}\left\lVert \widehat{\chi_t} \circ a_t -\mu_\mathcal{Y}(\widehat{\chi_t}\circ a_t) \right\rVert_{L^1_\mathcal{Y}}\\
    &\ll \frac{M}{\sqrt{N}}\sup_{t\geq 0}\left\lVert \widehat{\chi_t} \circ a_t  \right\rVert_{L^1_\mathcal{Y}},
\end{align*}
hence, by (\ref{third estimate}) and provided that 
\begin{equation}
\label{eq:cc1}
M=o(N^{1/2})
\end{equation}
we have 
$$
\left\lVert \mathsf{F}_N -\mathsf{F}_{N,M}  \right\rVert_{L^1_\mathcal{Y}} \rightarrow 0,\quad \text{as }N\rightarrow \infty.
$$

\subsubsection{Averages for the Siegel transform of a smooth approximation}
Further, we observe that the averages $\mathsf{F}_{N,M}$ have the same convergence in distribution if the characteristic function $\chi$ is replaced by the smooth approximation $f_\varepsilon$ introduced earlier. Indeed, if we consider the averages
\begin{equation}
\mathsf{F}_{N,M}^{(\varepsilon)} \coloneqq \frac{1}{\sqrt{N-M}} \sum_{t=M}^{N-1}\left( \widehat{f}_{\varepsilon} \circ a_t -\mu_\mathcal{Y}(\widehat{f}_{\varepsilon}\circ a_t) \right),
\end{equation}
with the parameter $\varepsilon = \varepsilon(N)$, $\varepsilon(N) \xrightarrow{N \to \infty} 0 $ to be specified later, then, Proposition \ref{smooth approximation} implies
\begin{equation}
 \left\lVert  \mathsf{F}_{N,M}- \mathsf{F}_{N,M}^{(\varepsilon)} \right\rVert_{L^1_\mathcal{Y}} \leq \frac{2}{\sqrt{N-M}}\sum_{t=M}^{N-1}\left\lVert \widehat{f}_{\varepsilon} \circ a_t -\widehat{\chi_t}\circ a_t \right\rVert_{L^1_\mathcal{Y}} \ll (N-M)^{1/2}(\varepsilon+e^{-\theta M})~.
\end{equation}
We will choose $\varepsilon$ and $M$ such that
\begin{equation}
\label{eq:cc2}
(N-M)^{1/2}\varepsilon \rightarrow 0 \quad \text{and} \quad (N-M)^{1/2}e^{-\theta M}\rightarrow 0~,
\end{equation}
which yields
$$
\left\lVert  \mathsf{F}_{N,M}- \mathsf{F}_{N,M}^{(\varepsilon)} \right\rVert_{L^1_\mathcal{Y}} \rightarrow 0 \quad \text{as}\quad N\rightarrow \infty~.
$$

\subsubsection{Averages for the truncated Siegel transform}
Finally, we also have the same convergence in distribution for the averages of the truncated Siegel transform
\begin{equation}\label{eq:def F_N,M}
\mathsf{F}_{N,M}^{(\varepsilon, L)} \coloneqq \frac{1}{\sqrt{N-M}} \sum_{t=M}^{N-1}\left( \widehat{f}_{\varepsilon}^{(L)} \circ a_t -\mu_\mathcal{Y}(\widehat{f}_{\varepsilon}^{(L)}\circ a_t) \right),
\end{equation}
defined for parameters $\varepsilon(N) \xrightarrow{N \to \infty} 0$ and  $L(N) \xrightarrow{N \to \infty} \infty $ to be specified later. \\
We assume that 
\begin{equation}
\label{eq:cc3}
M\gg \log L
\end{equation}
such that Proposition \ref{non-espace of mass} applies when $t\geq M$. Since the family of functions $f_\varepsilon$ is uniformly bounded by a compactly supported function, the estimate (\ref{L1 bound}) gives
\begin{align*}
\left\lVert \mathsf{F}_{N,M}^{(\varepsilon)} -\mathsf{F}_{N,M}^{(\varepsilon, L)}\right\rVert_{L^1_\mathcal{Y}} &\leq \frac{1}{\sqrt{N-M}}\sum_{t=M}^{N-1} \left\lVert \left(\widehat{f}_{\varepsilon}\circ a_t -\widehat{f}_{\varepsilon}^{(L)}\circ a_t\right)-\mu_\mathcal{Y}\left(\widehat{f}_{\varepsilon}\circ a_t -\widehat{f}_{\varepsilon}^{(L)}\circ a_t\right) \right\rVert_{L^1_\mathcal{Y}}\\
&\leq \frac{2}{\sqrt{N-M}}\sum_{t=M}^{N-1} \left\lVert \widehat{f}_{\varepsilon}\circ a_t -\widehat{f}_{\varepsilon}^{(L)}\circ a_t \right\rVert_{L^1_\mathcal{Y}}\\
&\ll_\tau (N-M)^{1/2}L^{-\tau/2},\qquad \text{for all } ~\tau<d.
\end{align*}
We will choose $L(N) \xrightarrow{N \to \infty} \infty $ such that 
\begin{align}
\label{eq:cc4}
N-M&=o(L^p) \quad \text{for some }p<d,
\end{align}
to obtain 
$$
\left\lVert \mathsf{F}_{N,M}^{(\varepsilon)} -\mathsf{F}_{N,M}^{(\varepsilon, L)}\right\rVert_{L^1_\mathcal{Y}} \rightarrow 0,\quad \text{as }N\rightarrow \infty.
$$

Hence, if we prove the CLT for the sequence $(\mathsf{F}_{N,M}^{(\varepsilon, L)})$, then the CLT for $(\mathsf{F}_N)$ would follow.

\subsection{Cumulants of the counting function}\label{sec:cumulants}

\subsubsection{The method of cumulants}
\label{subsec:method of cumulants}
We recall in this section the general approach of the method of cumulants (presented in \cite{bjoerklund2020} and \cite{bjoerklund2018central}) to establish the convergence to a normal distribution using a characterization by the cumulants.\\
Given a probability space $(X,\mu)$ and bounded measurable functions $\varphi_1,\cdots,\varphi_r$ on $X$, we define their \emph{joint cumulant} as
$$
\text{Cum}_\mu^{(r)}(\varphi_1,\cdots,\varphi_r)= \sum_{\mathcal{P}}(-1)^{|\mathcal{P}|-1}(|\mathcal{P}|-1)!\displaystyle\prod_{I\in \mathcal{P}}\int_X\left(\displaystyle\prod_{i\in I}\varphi_i \right)d\mu ~,
$$
where the sum is over all partitions $\mathcal{P}$ of the set $\{1,\cdots,r\}$. For a bounded measurable function $\varphi$ on $X$ we write 
$$
\text{Cum}_\mu^{(r)}(\varphi)=\text{Cum}_\mu^{(r)}(\varphi,\cdots,\varphi)~.
$$
We will use the following classical CLT-criterion (see \cite{frechet1931cumulants}).

\begin{prop}
\label{prop:method cumulants}
Let $(f_N)_{N\geq 1}$ be a sequence of real-valued bounded measurable functions such that
\begin{equation}
\label{characterization 1 law Norm}
\int_X f_N\,d\mu=0 ~ , \qquad \sigma^2\coloneqq \lim_{N\rightarrow \infty}\int_X f_N^2\,d\mu < \infty
\end{equation}
and
\begin{equation}
\label{characterization 2 law Norm}
\lim_{N\rightarrow \infty}\emph{Cum}_\mu^{(r)}(f_N)=0~, \quad \text{for all } r\geq 3~.
\end{equation}
Then for every $\xi \in \mathbb{R}$,
$$
\mu\left( \left\lbrace f_N<\xi\right\rbrace\right) \rightarrow \emph{Norm}_{\sigma}(\xi) \quad \text{as } N\rightarrow \infty~.
$$
\end{prop} 
The method of cumulants is equivalent to the more widely known ``method of moments", but the cumulants offer the following convenient cancelation property.\\

For a partition $\mathcal{Q}$ of $\{1,\cdots,r \}$, we define the conditional \emph{joint cumulant} with respect to $\mathcal{Q}$ by
$$
\text{Cum}_\mu^{(r)}(\varphi_1,\cdots,\varphi_r|\mathcal{Q})= \sum_{\mathcal{P}}(-1)^{|\mathcal{P}|-1}(|\mathcal{P}|-1)!\displaystyle\prod_{I\in \mathcal{P}}\displaystyle\prod_{J\in \mathcal{Q}}\int_X\left(\displaystyle\prod_{i\in I\cap
J}\varphi_i \right)d\mu ~.
$$

\begin{prop} \cite{bjoerklund2020}
\label{prop:cancelation cumulants} For any partition $\mathcal{Q}$ with $|\mathcal{Q}|\geq 2$, 
\begin{equation}
\label{cancelation cumulants}
\emph{Cum}_\mu^{(r)}(\varphi_1,\cdots,\varphi_r|\mathcal{Q})=0~,
\end{equation}
for all $\varphi_1,\cdots,\varphi_r \in L^\infty(X,\mu)$.
\end{prop}

\subsubsection{Estimating the cumulants}
\label{subsec:estimating cumulants}

It will be convenient to write
$$
\psi^{(\varepsilon,L)}_t(y):=\widehat{f}_\varepsilon^{(L)}(a_ty)-\mu_\mathcal{Y}(\widehat{f}_\varepsilon^{(L)}\circ a_t),
$$
so that the averaging function is 
$$
\mathsf{F}_{N,M}^{(\varepsilon, L)}=\frac{1}{\sqrt{N-M}}\sum_{t=M}^{N-1}\psi^{(\varepsilon,L)}_t \quad \text{with}\quad \int_\mathcal{Y}\mathsf{F}_{N,M}^{(\varepsilon, L)}~d\mu_\mathcal{Y}=0.$$
Our aim in this section is to estimate the following joint cumulants for $r\geq 3$,
\begin{equation}
\label{equ:cumulants F_N}
\text{Cum}_{\mu_\mathcal{Y}}^{(r)}\left(\mathsf{F}_{N,M}^{(\varepsilon, L)}\right)=\frac{1}{(N-M)^{r/2}}\sum_{t_1,\dots,t_r=M}^{N-1}\text{Cum}_{\mu_\mathcal{Y}}^{(r)}\left(\psi^{(\varepsilon,L)}_{t_1},\dots, \psi^{(\varepsilon,L)}_{t_r}\right).
\end{equation}
We reproduce below the argument as developed in  \cite{bjoerklund2020} and \cite{bjoerklund2018central}, taking into account the dependence on the parameters $L$ and $\varepsilon$ coming from the truncated Siegel transform and the smooth approximation respectively. The main idea in estimating these joint cumulants is to decompose (\ref{equ:cumulants F_N}) into sub-sums corresponding to ``separated" or ``clustered" tuples $t_1,\dots,t_r$ and to control their sizes.\\

\paragraph{\textbf{Separated and clustered times $t_1,\dots,t_r$}}
It will be convenient to consider $\{0,\dots,d-1\}^r$ as a subset of $\mathbb{R}_+^{r+1}$ with the embedding $(t_1,\dots,t_r)\rightarrow (0,t_1,\dots,t_r)$.

Following the approach developed in \cite{bjoerklund2020}, we define for 
non-empty subsets $I$ and $J$ of $\{0,\ldots, r\}$ and $\overline{t} = (t_0,\ldots,t_r) \in \bR_+^{r+1}$, 
$$
\rho^{I}(\overline{t}) := \max\big\{ |t_i-t_j| \, : \, i,j \in I \big\}\quad\hbox{and}
\quad
\rho_{I,J}(\overline{t}):= \min\big\{ |t_i-t_j| \, : \, i \in I, \enskip j \in J \big\},
$$
and if $\cQ$ is a partition of $\{0,\ldots,r\}$, we set
$$
\rho^{\cQ}(\overline{t}) := \max\big\{ \rho^{I}(\overline{t}) \, : \, I \in \cQ \big\}
\quad\hbox{and}\quad
\rho_{\cQ}(\overline{t}) := \min\big\{ \rho_{I,J}(\overline{t}) \, : \, I \neq J, \enskip I, J \in \cQ \big\}.
$$
For $0 \leq \alpha < \beta$, we define
$$
\Delta_{\cQ}(\alpha,\beta) 
:= 
\big\{ 
\overline{t} \in \bR_+^{r+1} \, : \, 
\rho^{\cQ}(\overline{t}) \leq \alpha, 
\qen
\rho_{\cQ}(\overline{t}) > \beta \big\}
$$
and
$$
\Delta(\alpha):= 
\big\{ 
\overline{t} \in \bR_+^{r+1} \, : \, 
\rho(t_i,t_j) \leq \alpha \hbox{ for all $i,j$}\big\}.
$$
The following decomposition of $\bR_+^{r+1}$ was established in \cite[Prop. 6.2]{bjoerklund2020}: given 
\begin{equation}\label{eq:inequality}
0=\alpha_0<\beta_1<\alpha_1=(3+r)\beta_1<\beta_2<\cdots <\beta_r<\alpha_r=(3+r)\beta_r<\beta_{r+1},
\end{equation}
we have
\begin{equation}
\label{eq:decomp}
\bR_+^{r+1} = \Delta(\beta_{r+1}) \cup \Big( \bigcup_{j=0}^{r} \bigcup_{|\cQ| \geq 2} \Delta_{\cQ}(\alpha_j,\beta_{j+1}) \Big),
\end{equation}
where the union is taken over the partitions $\cQ$ of $\{0,\ldots,r\}$ with $|\cQ|\ge 2$. Upon taking restrictions, we also have 
\begin{equation}\label{eq:decomp_000}
\{M,\ldots,N-1\}^r =\Omega(\beta_{r+1};M,N) \cup \Big( \bigcup_{j=0}^{r} \bigcup_{|\cQ| \geq 2} \Omega_{\cQ}(\alpha_j,\beta_{j+1};M,N) \Big),
\end{equation}
for all $N> M \geq 0$, where 
\begin{align*}
\Omega(\beta_{r+1};M,N) &:=\{M,\ldots,N-1\}^r \cap \Delta(\beta_{r+1}),\\
\Omega_{\cQ}(\alpha_j,\beta_{j+1};M,N) &:=\{M,\ldots,N-1\}^r \cap \Delta_{\cQ}(\alpha_j,\beta_{j+1}).
\end{align*}
In order to estimate the cumulant \eqref{equ:cumulants F_N}, we shall separately estimate the sums over $\Omega(\beta_{r+1};M,N)$ and $\Omega_{\cQ}(\alpha_j,\beta_{j+1};M,N)$; the exact choices of the sequences $(\alpha_j)$ and $(\beta_j)$ will be fixed later.\\
We may first choose 
\begin{equation}
M>\beta_{r+1}
\end{equation}
so that $\Omega(\beta_{r+1};M,N)=\emptyset$ and does not contribute to the sum.\\

\paragraph{\textbf{\underline{Case 1:} Summing over $(t_1,\ldots,t_r)\in \Omega_{\cQ}(\alpha_j,\beta_{j+1};M,N)$ with $\cQ=\{\{0\},\{1,\ldots,r\}\}$}}.
We shall first show that, in this case, we have
\begin{equation}
\label{eq:cc2_0}
\hbox{Cum}_{\mu_\cY}^{(r)}\left(\psi^{(\varepsilon,L)}_{t_1},\ldots,\psi^{(\varepsilon,L)}_{t_r}\right)\approx
\hbox{Cum}_{\mu_\cX}^{(r)}\left(\phi^{(\varepsilon,L)}\circ a_{t_1},\ldots,\phi^{(\varepsilon,L)}\circ a_{t_r}\right)
\end{equation}
where $\phi^{(\varepsilon,L)}:= \widehat f_\varepsilon^{(L)}-\mu_\cX(\widehat f_\varepsilon^{(L)})$.
This reduces to estimating the integrals
\begin{align}
&\int_{\cY} \left(\prod_{i\in I} \psi^{(\varepsilon,L)}_{t_i}\right)\,d\mu_\cY \label{eq:psi}\\
=&
\sum_{J\subset I} (-1)^{|I\backslash J|} 
\left(\int_{\cY} \left(\prod_{i\in J} \widehat f_\varepsilon^{(L)}\circ a_{t_i}\right)\,d\mu_\cY\right)
\prod_{i\in I\backslash J} \left(\int_{\cY} (\widehat f_\varepsilon^{(L)}\circ a_{t_i})\,d\mu_\cY\right).\nonumber
\end{align}
If $(t_1,\ldots,t_r)\in \Omega_\cQ(\alpha_j,\beta_{j+1};M,N)$, and thus
$$
|t_{i_1}-t_{i_2}|\le \alpha_j \quad\hbox{and}\quad t_{i_1}\ge \beta_{j+1}\quad\quad\hbox{for all $1\le i_1,i_2\le r,$}
$$
it follows from Corollary \ref{cor:multiple equidistribution K orbits} with $r=1$ that there exists $\delta>0$ such that
\begin{equation}\label{eq:psi1}
\int_{\cY} (\widehat f_\varepsilon^{(L)}\circ a_{t_i})\,d\mu_\cY
=\mu_{\cX}\left(\widehat f_\varepsilon^{(L)}\right)+O\left(e^{-\delta \beta_{j+1} }\,\left\|\widehat f_\varepsilon^{(L)}\right\|_{C^l}\right).
\end{equation}
For a fixed $J \subset I$, we define
$$
\Phi^{(\varepsilon,L)}:=\prod_{i\in J} \widehat f_\varepsilon^{(L)}\circ a_{t_i-t_1},
$$
and note that for some $\xi=\xi(d,l)>0$, we have
$$
\left\|\Phi^{(\varepsilon,L)}\right\|_{C^l}\ll \prod_{i\in J} \left\|\widehat f_\varepsilon^{(L)}\circ a_{t_i-t_1}\right\|_{C^l} \ll e^{|J|\xi\, \alpha_j}\,\left\|\widehat f_\varepsilon^{(L)}\right\|_{C^l}^{|J|}.
$$
If we again apply Corollary \ref{cor:multiple equidistribution K orbits} to the function $\Phi^{(\varepsilon,L)}$, we obtain 
\begin{align}
\int_{\cY} \left(\prod_{i\in J} \widehat f_\varepsilon^{(L)}\circ a_{t_i}\right)\,d\mu_\cY
&=\int_{\cY} (\Phi^{(\varepsilon,L)}\circ a_{t_1})\,d\mu_\cY \label{eq:PPsi}\\
&=\int_{\cX}\Phi^{(\varepsilon,L)}\, d\mu_\cX+O\left(e^{-\delta \beta_{j+1} }\,\left\|\Phi^{(\varepsilon,L)}\right\|_{C^l}\right) \nonumber\\
&=\int_{\cX}\left(\prod_{i\in J} \widehat f_\varepsilon^{(L)}\circ a_{t_i}\right)\, d\mu_\cX+O\left(e^{-\delta \beta_{j+1} } e^{r\xi\, \alpha_j}\,\left\|\widehat f_\varepsilon^{(L)}\right\|_{C^l}^{|J|}\right), \nonumber
\end{align}
where we used that $\mu_\cX$ is invariant under the transformation $a_t$.
Let us now choose the exponents $\alpha_j$ and $\beta_{j+1}$ so that
$\delta \beta_{j+1}-r\xi\alpha_j>0$. 
Combining \eqref{eq:psi}, \eqref{eq:psi1} and \eqref{eq:PPsi}, we deduce that 
\begin{align}
\int_{\cY} \left(\prod_{i\in I} \psi^{(\varepsilon,L)}_{t_i}\right)\,d\mu_\cY \label{eq:prod1}
=&
\sum_{J\subset I} (-1)^{|I\backslash J|} 
\left(\int_{\cX} \left(\prod_{i\in J} \widehat f_\varepsilon^{(L)}\circ a_{t_i}\right)\,d\mu_\cX\right)
\mu_{\cX}\left(\widehat f_\varepsilon^{(L)}\right)^{|I\backslash J|}\\
&+O\left(e^{-\delta \beta_{j+1} } e^{r\xi\, \alpha_j}\,\left\|\widehat f_\varepsilon^{(L)}\right\|_{C^l}^{|I|}\right) \nonumber\\
=& \int_{\cX} \prod_{i\in I} \left(\widehat f_\varepsilon^{(L)}\circ a_{t_i}-\mu_\cX(\widehat f_\varepsilon^{(L)})\right)\,d\mu_\cX
+O\left(e^{-(\delta \beta_{j+1}-r\xi\alpha_j )}\,\left\|\widehat f_\varepsilon^{(L)}\right\|_{C^l}^{|I|}\right), \nonumber
\end{align}
and thus, for any partition $\cP$,
$$
\prod_{I\in \cP}\int_{\cY} \left(\prod_{i\in I} \psi^{(\varepsilon,L)}_{t_i}\right)\,d\mu_\cY
=
\prod_{I\in\cP}\int_{\cX} \left(\prod_{i\in I} \phi^{(\varepsilon,L)}\circ a_{t_i}\right)\,d\mu_\cX
+O\left(e^{-(\delta \beta_{j+1}-r\xi\alpha_j )}\,\left\|\widehat f_\varepsilon^{(L)}\right\|_{C^l}^{r}\right),
$$
and consequently,
\begin{align}\label{eq:cumm0}
\hbox{Cum}_{\mu_\cY}^{(r)}\left(\psi^{(\varepsilon,L)}_{t_1},\ldots,\psi^{(\varepsilon,L)}_{t_r}\right)=&
\;\hbox{Cum}_{\mu_\cX}^{(r)}\left(\phi^{(\varepsilon,L)}\circ a_{t_1},\ldots,\phi^{(\varepsilon,L)}\circ a_{t_r}\right)\\
&+O\left(e^{-(\delta \beta_{j+1}-r\xi\, \alpha_j)}\,\left\|\widehat f_\varepsilon^{(L)}\right\|_{C^l}^{r}\right)\nonumber
\end{align}
whenever $(t_1,\ldots,t_r)\in \Omega_\cQ(\alpha_j,\beta_{j+1};M,N)$ with $\cQ=\{\{0\},\{1,\ldots,r\}\}$,
from which \eqref{eq:cc2_0} follows. \\

We now claim that 
\begin{equation}
\label{eq:cum0_bound}
\left|\hbox{Cum}_{\mu_\cX}^{(r)}\left(\phi^{(\varepsilon,L)}\circ a_{t_1},\ldots,\phi^{(\varepsilon,L)}\circ a_{t_r}\right)\right| \ll  \left\|\widehat f_\varepsilon^{(L)}\right\|_{C^0}^{(r-(d-1))^+} \left\|\widehat f_\varepsilon^{(L)} \right\|_{L^{d-1}(\cX)}^{\min(r,d-1)},
\end{equation}
where we use the notation $x^+=\max(x,0)$. The implied constants in \eqref{eq:cum0_bound} and below depend only on $\hbox{supp}(f_\varepsilon)$, so that it is uniform in $\varepsilon$. By the definition of the cumulant, to prove  \eqref{eq:cum0_bound}, it suffices to show 
that for every $z\ge 1$ and indices $i_1,\ldots,i_z$,
\begin{equation}
\label{eq:cum0_bound_0}
\int_{\cX} \left|\left(\phi^{(\varepsilon, L)}\circ a_{t_{i_1}}\right)\cdots\left(\phi^{(\varepsilon,L)}\circ a_{t_{i_z}}\right)\right|\, d\mu_\cX
\ll \left\|\widehat f_\varepsilon^{(L)}\right\|_{C^0}^{(z-(d-1))^+}
\left\|\widehat f_\varepsilon^{(L)} \right\|_{L^{d-1}(\cX)}^{\min(z,d-1)}.
\end{equation}
Using 
the generalized H\"older inequality, we deduce that when $z\le d-1$,
\begin{align*}
\int_{\cX} \left|\left(\phi^{(\varepsilon,L)}\circ a_{t_{i_1}}\right)\cdots\left(\phi^{(\varepsilon,L)}\circ a_{t_{i_z}}\right)\right|\, d\mu_\cX&\le 
\left\|\phi^{(\varepsilon,L)}\circ a_{t_{i_1}}\right\|_{L^{d-1}(\cX)}\cdots\left\|\phi^{(\varepsilon,L)}\circ a_{t_{i_z}}\right\|_{L^{d-1}(\cX)}\\
&\ll \left\|\widehat f_\varepsilon^{(L)}\right\|_{L^{d-1}(\cX)}^z.
\end{align*}
Also, when $z>d-1$, 
\begin{align*}
&\int_{\cX} \left|\left(\phi^{(\varepsilon,L)}\circ a_{t_{i_1}}\right)\cdots\left(\phi^{(\varepsilon,L)}\circ a_{t_{i_z}}\right)\right|\, d\mu_\cX\\
\le& \left\|\phi^{(\varepsilon,L)}\right\|_{C^0}^{z-(d-1)}
\int_{\cX} \left|\left(\phi^{(\varepsilon,L)}\circ a_{t_{i_1}}\right)\cdots\left(\phi^{(\varepsilon,L)}\circ a_{t_{i_{d-1}}}\right)\right|\, d\mu_\cX\\
\ll & \left\|\widehat f_\varepsilon^{(L)}\right\|_{C^0}^{z-(d-1)} \left\|\widehat f_\varepsilon^{(L)}\right\|_{L^{d-1}(\cX)}^{d-1}.
\end{align*}
This implies \eqref{eq:cum0_bound_0} and \eqref{eq:cum0_bound}. \\

Finally, we recall that if $(t_1,\ldots,t_r)\in \Omega_\cQ(\alpha_j,\beta_{j+1};M,N)$ with $\cQ=\{\{0\},\{1,\ldots,r\}\}$, then we have $|t_{i_1}-t_{i_2}|\le \alpha_j$
for all $i_1\ne i_2$, and thus
\begin{equation}
|\Omega_{\cQ}(\alpha_j,\beta_{j+1};M,N)|\ll (N-M)\alpha_j^{r-1}\label{eq:card}.
\end{equation}

Combining  \eqref{eq:cumm0}, \eqref{eq:cum0_bound} and \eqref{eq:card} in \eqref{equ:cumulants F_N}, and using Proposition \ref{bounds truncated siegel transform} with (\ref{epsilon approximation of chi}), we conclude that
\begin{align*}
 &\frac{1}{(N-M)^{r/2}}  \sum_{\overline{t}\in \Omega_{\cQ}(\alpha_j,\beta_{j+1};M,N)}\text{Cum}_{\mu_\mathcal{Y}}^{(r)}\left(\psi^{(\varepsilon,L)}_{t_1},\dots, \psi^{(\varepsilon,L)}_{t_r}\right)   \\
&\ll \; (N-M)^{r/2}\, e^{-(\delta\beta_{j+1}- r\alpha_j\xi)}\, \left\|\widehat f_\varepsilon^{(L)}\right\|_{C^l}^r +(N-M)^{1-r/2} \alpha_j^{r-1} 
\left\|\widehat f_\varepsilon^{(L)}\right\|_{C^0}^{(r-(d-1))^+} \left\|\widehat f_\varepsilon^{(L)} \right\|_{L^{d-1}(\cX)}^{\min(r,d-1)} \\
&\ll \; (N-M)^{r/2}\, e^{-(\delta\beta_{j+1}- r\alpha_j\xi)}\, L^r \varepsilon^{-rl} 
+~(N-M)^{1-r/2} \alpha_j^{r-1} 
L^{(r-(d-1))^+}.
\end{align*}

\paragraph{\textbf{\underline{Case 2:} Summing over $(t_1,\ldots,t_r)\in \Omega_{\cQ}(\alpha_j,\beta_{j+1};M,N)$ with $|\cQ|\ge 2$ and 
	$\cQ\ne \{\{0\},\{1,\ldots,r\}\}$}.}
	
In this case, the partition $\cQ$ defines a non-trivial partition $\cQ'=\{I_0,\ldots,I_\ell\}$ of $\{1,\ldots,r\}$ such that 
for all $(t_1,\ldots, t_r)\in \Omega_Q(\alpha_j,\beta_{j+1};M,N)$, we have 
\begin{equation}
\label{eq:ineq0}
|t_{i_1}-t_{i_2}|\le \alpha_j \;\; \hbox{if $i_1\sim_{\cQ'} i_2$} \qand
|t_{i_1}-t_{i_2}|> \beta_{j+1} \;\; \hbox{if $i_1\not\sim_{\cQ'} i_2$},
\end{equation}
and
\[
t_i\le \alpha_j\;\; \hbox{for all $i\in I_0$},\nonumber \qand t_i>\beta_{j+1} \;\; \hbox{for all $i\notin I_0$}. 
\]
In particular,
\begin{equation}
\label{eq:Dcase2}
D(t_{i_1},\ldots, t_{i_\ell})\ge \beta_{j+1},
\end{equation}
Let $I$ be an arbitrary subset of $\{1,\ldots,r\}$; we shall show that 
\begin{equation}
\label{eq:approx}
\int_{\cY} \left( \prod_{i\in I} \psi_{t_i}^{(\varepsilon,L)}\right)\, d\mu_\cY
\approx \prod_{h=0}^\ell \left(\int_{\cY} \left(\prod_{i\in I\cap I_h} \psi_{t_i}^{(\varepsilon,L)} \right)\,d\mu_\cY\right),
\end{equation}
where we henceforth shall use the convention that the product is equal to one when $I\cap I_h=\emptyset$. \\

Let us estimate  the right hand side of \eqref{eq:approx}. We begin by setting
$$
\Phi_0^{(\varepsilon,L)}:=\prod_{i\in I\cap I_0} \psi_{t_i}^{(\varepsilon,L)}.
$$
It is easy to see that there exists $\xi=\xi(n,l)>0$ such that
\begin{equation}
\label{eq:Phi_0}
\|\Phi_0^{(\varepsilon,L)}\|_{C^l}\ll \prod_{i\in I\cap I_0} \|\widehat f_\varepsilon^{(L)}\circ a_{t_i}-\mu_\cY(\widehat f_\varepsilon^{(L)}\circ a_{t_i})\|_{C^l} \ll e^{|I\cap I_0|\xi\, \alpha_j}\,\|\widehat f_\varepsilon^{(L)}\|_{C^l}^{|I\cap I_0|}.
\end{equation}
To prove \eqref{eq:approx}, 
we expand $\psi_{t_i}^{(\varepsilon,L)}=\widehat f_\varepsilon^{(L)}\circ a_{t_i}-\mu_\cY(\widehat f_\varepsilon^{(L)}\circ a_{t_i})$
for $i\in I\backslash I_0$ and get
\begin{eqnarray}
\label{eq:psi00}
\int_{\cY} \left(\prod_{i\in I} \psi^{(\varepsilon,L)}_{t_i}\right)\,d\mu_\cY 
&=& 
\sum_{J\subset I\backslash I_0} (-1)^{|I\backslash (J\cup I_0)|} \cdot \\
&\cdot&
\left(\int_{\cY} \Phi_0^{(\varepsilon,L)} \left(\prod_{i\in J} \widehat f_\varepsilon^{(L)}\circ a_{t_i}\right)\,d\mu_\cY\right) 
\prod_{i\in I\backslash (J\cup I_0)} \left(\int_{\cY} (\widehat f_\varepsilon^{(L)}\circ a_{t_i})\,d\mu_\cY\right).\nonumber
\end{eqnarray}
We recall that when $i\notin I_0$, we have $t_i\ge \beta_{j+1}$, and thus it follows from Corollary \ref{cor:multiple equidistribution K orbits} with $r=1$ that
\begin{equation}\label{eq:psi_0}
\int_{\cY} (\widehat f_\varepsilon^{(L)}\circ a_{t_i})\,d\mu_\cY
=\mu_{\cX}(\widehat f_\varepsilon^{(L)})+O\left(e^{-\delta \beta_{j+1} }\,\|\widehat f_\varepsilon^{(L)}\|_{C^l}\right), \quad \hbox{ with $i\notin I_0$}.
\end{equation}
To estimate the other integrals in \eqref{eq:psi00}, we also apply Corollary \ref{cor:multiple equidistribution K orbits}.
Let us  first fix a subset $J \subset I \setminus I_0$ and for each $1 \leq h \leq l$, we pick $i_h\in I_h$, 
and set
$$
\Phi_h^{(\varepsilon,L)}:=\prod_{i\in J\cap I_h} \widehat f_\varepsilon^{(L)}\circ a_{t_i-t_{i_h}}.
$$
Then 
$$
\int_{\cY} \Phi_0^{(\varepsilon,L)} \left(\prod_{i\in J} \widehat f_\varepsilon^{(L)}\circ a_{t_i}\right)\,d\mu_\cY=
\int_{\cY} \Phi_0^{(\varepsilon,L)} \left( \prod_{h=1}^\ell \Phi_h^{(\varepsilon,L)}\circ a_{t_{i_h}}
\right)\,d\mu_\cY.
$$
We note that for $i\in I_h$, we have $|t_i-t_{i_h}|\le \alpha_j$, and 
thus there exists $\xi=\xi(n,l)>0$ such that
\begin{equation}
\label{eq:sss}
\|\Phi_h^{(\varepsilon,L)}\|_{C^l}\ll \prod_{i\in J\cap I_h} \|\widehat f_\varepsilon^{(L)}\circ a_{t_i-t_{i_h}}\|_{C^l} \ll e^{|J\cap I_h|\xi\, \alpha_j}\,\|\widehat f_\varepsilon^{(L)}\|_{C^l}^{|J\cap I_h|}.
\end{equation}
Using \eqref{eq:Dcase2}, Corollary \ref{cor:multiple equidistribution K orbits} implies that
\begin{align*}
\int_{\cY} \Phi_0^{(\varepsilon,L)} \left( \prod_{h=1}^\ell \Phi_h^{(\varepsilon,L)}\circ a_{t_{i_h}}\right)\,d\mu_\cY
= &\left(\int_{\cY}\Phi_0^{(\varepsilon,L)} \, d\mu_\cY \right)\prod_{h=1}^\ell \left(\int_{\cX}\Phi_h^{(\varepsilon,L)} \, d\mu_\cX\right)\\
 &+O\left( e^{-\delta \beta_{j+1}}\, \prod_{h=0}^\ell \|\Phi_h^{(\varepsilon,L)}\|_{C^l}\right).
\end{align*}
Using \eqref{eq:Phi_0} and \eqref{eq:sss} and invariance of  the measure $\mu_\cX$, we deduce that
\begin{align*}
\int_{\cY} \Phi_0^{(\varepsilon,L)} \left( \prod_{h=1}^\ell \Phi_h^{(\varepsilon,L)}\circ a_{t_{i_h}}\right)\,d\mu_\cY
= &\left(\int_{\cY}\Phi_0^{(\varepsilon,L)} \, d\mu_\cY \right)\prod_{h=1}^\ell \left(\int_{\mathcal{X}}\left(\prod_{i\in J\cap I_h} \widehat f_\varepsilon^{(L)}\circ a_{t_i} \right)\, d\mu_\cX\right)\\
&+O\left( e^{-(\delta \beta_{j+1}-r\xi \alpha_j)}\, \|\widehat f_\varepsilon^{(L)}\|_{C^l}^{|(I\cap I_0)\cup J|}\right).
\end{align*}
Hence, we conclude that
\begin{align}
\int_{\cY} \Phi_0^{(\varepsilon,L)} \left(\prod_{i\in J} \widehat f_\varepsilon^{(L)}\circ a_{t_i}\right)\,d\mu_\cY
= &\left(\int_{\mathcal{Y}}\Phi_0^{(\varepsilon,L)} \, d\mu_\cY \right)\prod_{h=1}^\ell \left(\int_{\mathcal{X}}\left(\prod_{i\in J\cap I_h} \widehat f_\varepsilon^{(L)}\circ a_{t_i} \right)\, d\mu_\cX\right) \label{eq:p000}\\
&+O\left( e^{-(\delta \beta_{j+1}-r\xi \alpha_j)}\,  \|\widehat f_\varepsilon^{(L)}\|_{C^l}^{|(I\cap I_0)\cup J|}\right).\nonumber
\end{align}
We shall choose the parameters $\alpha_j$ and $\beta_{j+1}$ so that 
\begin{equation}\label{eq:i1}
\delta\beta_{j+1}-r\xi\alpha_j>0.
\end{equation}
Substituting \eqref{eq:psi_0} and \eqref{eq:p000} in \eqref{eq:psi00}, we deduce that 
\begin{align} \label{eq:ppsi}
&\int_{\cY} \left(\prod_{i\in I} \psi_{t_i}^{(\varepsilon,L)}\right)\,d\mu_\cY  \\ 
=&
\sum_{J\subset I\backslash I_0} (-1)^{|I\backslash (J\cup I_0)|} 
\left(\int_{\mathcal{Y}}\Phi_0^{(\varepsilon,L)} \, d\mu_\cY \right)\prod_{h=1}^\ell \left(\int_{\mathcal{X}}\left(\prod_{i\in J\cap I_h} \widehat f_\varepsilon^{(L)}\circ a_{t_i} \right)\, d\mu_\cX\right) \mu_\cX(\widehat f_\varepsilon^{(L)})^{|I\backslash (J\cup I_0)|} \nonumber \\
& +O \left( e^{-(\delta \beta_{j+1}-r\xi \alpha_j)}\, \|\widehat f_\varepsilon^{(L)}\|_{C^l}^{|I|}\right).\nonumber
\end{align}

Next, we estimate the right hand side of \eqref{eq:approx}.
Let us fix $1 \leq h \leq l$ and for a subset $J \subset I \cap I_h$, we define 
$$
\Phi_J^{(\varepsilon,L)}:=\prod_{i\in J} \widehat f_\varepsilon^{(L)}\circ a_{t_i-t_{i_h}}.
$$
As in \eqref{eq:sss}, for some $\xi>0$,
$$
\|\Phi_J^{(\varepsilon,L)}\|_{C^k}\ll \prod_{i\in J} \|\widehat f_\varepsilon^{(L)}\circ a_{t_i-t_{i_h}}\|_{C^l} \ll e^{|J|\xi\, \alpha_j}\,\|\widehat f_\varepsilon^{(L)}\|_{C^l}^{|J|}.
$$
Applying Corollary \ref{cor:multiple equidistribution K orbits} to the function $\Phi_J^{(\varepsilon,L)}$ and using that $t_{i_h}>\beta_{j+1}$, we get 
\begin{align}
\int_{\cY} \left(\prod_{i\in J} \widehat f_\varepsilon^{(L)}\circ a_{t_i}\right)\,d\mu_\cY
&=\int_{\cY} (\Phi_J^{(\varepsilon,L)}\circ a_{t_{i_h}})\,d\mu_\cY \label{eq:PPsi_0}\\
&=\int_{\cX}\Phi_J^{(\varepsilon,L)}\, d\mu_\cX+O\left(e^{-\delta \beta_{j+1} }\,\|\Phi_J^{(\varepsilon,L)}\|_{C^l}\right) \nonumber\\
&=\int_{\cX}\left(\prod_{i\in J} \widehat f_\varepsilon^{(L)}\circ a_{t_i}\right)\, d\mu_\cX+O\left(e^{-\delta \beta_{j+1} } e^{r\xi\, \alpha_j}\,\|\widehat f_\varepsilon^{(L)}\|_{C^l}^{|J|}\right), \nonumber
\end{align}
where we have used $a$-invariance of $\mu_\cX$.
Combining \eqref{eq:psi_0} and \eqref{eq:PPsi_0}, we deduce that 
\begin{align}
&\int_{\cY} \left(\prod_{i\in I\cap I_h} \psi_{t_i}^{(\varepsilon,L)}\right)\,d\mu_\cY \nonumber\\
&=
\sum_{J\subset I\cap I_h} (-1)^{|(I\cap I_h)\backslash J|} 
\left(\int_{\cX} \left(\prod_{i\in J} \widehat f_\varepsilon^{(L)}\circ a_{t_i}\right)\,d\mu_\cX\right)
\mu_{\cX}(\widehat f_\varepsilon^{(L)})^{|(I\cap I_h)\backslash J|}\nonumber \\
&\qquad\qquad+O\left(e^{-\delta \beta_{j+1} } e^{r\xi\, \alpha_j}\,\|\widehat f_\varepsilon^{(L)}\|_{C^l}^{|I\cap I_h|}\right) \label{eq:prod2}\\
=& \int_{\cX} \prod_{i\in I\cap I_h} \left(\widehat f_\varepsilon^{(L)}\circ a_{t_i}-\mu_\cX(\widehat f_\varepsilon^{(L)})\right)\,d\mu_\cX
+O\left(e^{-(\delta \beta_{j+1}-r\xi\alpha_j )}\,\|\widehat f_\varepsilon^{(L)}\|_{C^l}^{|I\cap I_h|}\right),\nonumber
\end{align}
which implies
\begin{align*}
&\prod_{h=0}^\ell \left(\int_{\cY} \left(\prod_{i\in I\cap I_h} \psi_{t_i}^{(\varepsilon,L)} \right)\,d\mu_\cY\right)\\
=& \left(\int_{\mathcal{Y}}\Phi_0^{(\varepsilon,L)} \, d\mu_\cY \right) \prod_{h=1}^\ell \left(\int_{\cX} \prod_{i\in I\cap I_h} \left(\widehat f_\varepsilon^{(L)}\circ a_{t_i}-\mu_\cX(\widehat f_\varepsilon^{(L)})\right)\,d\mu_\cX\right)\\
&+O\left(e^{-(\delta \beta_{j+1} - r\xi\alpha_j)} \,\|\widehat f_\varepsilon^{(L)}\|_{C^l}^{r}\right).
\end{align*}
Furthermore, multiplying out the products over $I\cap I_h$, we get
\begin{align}\label{eq:ppp}
&\prod_{h=0}^\ell \left(\int_{\cY} \left(\prod_{i\in I\cap I_h} \psi_{t_i}^{(\varepsilon,L)} \right)\,d\mu_\cY\right)\\
=& \left(\int_{\mathcal{Y}}\Phi_0^{(\varepsilon,L)} \, d\mu_\cY \right) 
\sum_{J\subset I\backslash I_0} (-1)^{|I\backslash (I_0\cup J)|}\prod_{h=1}^\ell \left(\int_{\cX} \prod_{i\in I_h\cap J} \widehat f_\varepsilon^{(L)}\circ a_{t_i}\,d\mu_\cX\right)
\mu_\cX(\widehat f_\varepsilon^{(L)})^{|I\backslash (I_0\cup J)|} \nonumber \\
&+O\left(e^{-(\delta \beta_{j+1} - r\xi\alpha_j)} \,\|\widehat f_\varepsilon^{(L)}\|_{C^l}^{|I|}\right). \nonumber
\end{align}
Comparing \eqref{eq:ppsi} and \eqref{eq:ppp}, we finally conclude that
\begin{align*}
\int_{\cY} \left(\prod_{i\in I} \psi_{t_i}^{(\varepsilon,L)}\right)\,d\mu_\cY =& \prod_{h=0}^\ell \left(\int_{\cY} \left(\prod_{i\in I\cap I_h} \psi_{t_i}^{(\varepsilon,L)} \right)\,d\mu_\cY\right)\\
&+O\left(e^{-(\delta \beta_{j+1} - r\xi\alpha_j)} \,\|\widehat f_\varepsilon^{(L)}\|_{C^l}^{|I|}\right)
\end{align*}
when $(t_1,\ldots,t_r)\in\Omega_{\cQ}(\alpha_j,\beta_{j+1};M,N)$.
This establishes \eqref{eq:approx} with an explicit error term.
This estimate implies that for the partition $\cQ'=\{I_0,\ldots, I_\ell\}$,
$$
\hbox{Cum}_{\mu_\cY}^{(r)}(\psi_{t_1}^{(\varepsilon,L)},\ldots,\psi_{t_r}^{(\varepsilon,L)})=
\hbox{Cum}_{\mu_\cY}^{(r)}(\psi^{(\varepsilon,L)}_{t_1},\ldots,\psi^{(\varepsilon,L)}_{t_r}|\cQ')+
O\left(e^{-(\delta \beta_{j+1} - r\xi\alpha_j)} \,\|\widehat f_\varepsilon^{(L)}\|_{C^l}^{r}\right)
$$
By Proposition \ref{prop:cancelation cumulants},
$$
\hbox{Cum}_{\mu_\cY}^{(r)}(\psi^{(\varepsilon,L)}_{t_1},\ldots,\psi^{(\varepsilon,L)}_{t_r}|\cQ')=0,
$$
so it follows that for all $(t_1,\ldots, t_r)\in \Omega_Q(\alpha_j,\beta_{j+1};M,N)$,
\begin{equation}
\label{eq:cum3}
\left|\hbox{Cum}_{\mu_\cY}^{(r)}(\psi^{(\varepsilon,L)}_{t_1},\ldots,\psi_{t_r})\right|
\ll e^{-(\delta \beta_{j+1} - r\xi\alpha_j)} \,\|\widehat f_\varepsilon^{(L)}\|_{C^l}^{r}.
\end{equation}
It follows
\begin{align*}
\frac{1}{(N-M)^{r/2}}\sum_{(t_1,\ldots, t_r)\in \Omega_{\cQ}(\alpha_j,\beta_{j+1};M,N)} &\left|\hbox{Cum}_{\mu_\cY}^{(r)}(\psi^{(\varepsilon, L)}_{t_1},\ldots,\psi^{(\varepsilon, L)}_{t_r})\right|\\
&\ll (N-M)^{r/2}e^{-(\delta \beta_{j+1} - r\xi\alpha_j)} \,\left\|\widehat f_\varepsilon^{(L)}\right\|_{C^l}^{r}\\
&\ll (N-M)^{r/2}e^{-(\delta \beta_{j+1} - r\xi\alpha_j)} L^r\,\varepsilon^{-rl},
\end{align*}
where we used Proposition \ref{bounds truncated siegel transform} and \eqref{epsilon approximation of chi}.\\

\subsubsection{Final estimates on the cumulants}
Finally, we combine the established bounds to get the following estimate
\begin{align}
\left|\hbox{Cum}_{\mu_\cY}^{(r)}(\mathsf{F}_{N,M}^{(\varepsilon, L)})\right|\ll 
&~~ (N-M)^{1-r/2}\left({\max}_j\, \alpha_j^{r-1}\right) L^{(r-(d-1))^+} \label{eq:cumm_last} \\
&\;\; +
  (N-M)^{r/2} \left({\max}_{j}\,\, e^{-(\delta \beta_{j+1} - r\xi\alpha_j)}\right) L^r\,\varepsilon^{-rl}.\nonumber
\end{align}
This estimate holds provided that \eqref{eq:inequality} and \eqref{eq:i1} hold, namely when
\begin{equation}
\label{eq:alpha}
\alpha_j=(3+r)\beta_{j}<\beta_{j+1}\quad\hbox{and}\quad \delta\beta_{j+1}-r\xi\alpha_j>0\quad\quad\hbox{ for $j=1,\ldots,r$.}
\end{equation}
Given any $\gamma>0$, we define the parameters $\beta_j$ inductively by the formula
\begin{equation}
\label{eq:beta}
\beta_1=\gamma\quad\hbox{and}\quad \beta_{j+1}=\max\left(\gamma+(3+r)\beta_{j}, \gamma+\delta^{-1}r(3+r)\xi\beta_{j}\right).
\end{equation}
It easily follows by induction that $\beta_{r+1}\ll_r \gamma$, and choosing 
\begin{equation}
\label{eq:cc5}
M\gg_r \gamma
\end{equation}
we deduce from \eqref{eq:cumm_last} that
\begin{align}
\label{eq:cumm_smooth}
\left|\hbox{Cum}_{\mu_\cY}^{(r)}(\mathsf{F}_{N,M}^{(\varepsilon, L)})\right|
\ll (N-M)^{r/2}  e^{-\delta \gamma} L^{r} \varepsilon^{-rl}
+(N-M)^{1-r/2}   \gamma^{r-1} L^{(r-(d-1))^+}.
\end{align}
We observe that since $d\ge 2$, 
$$
\frac{(r-(d-1))^+}{d}<r/2-1
\quad\hbox{for all $r\ge 3$,}
$$
Hence, we can choose $q>1/n$ such that
$$
q(r-(d-1))^+<r/2-1
\quad\hbox{for all $r\ge 3$.}
$$
Then we select
$$
L=(N-M)^q,
$$
so that, in particular, the condition (\ref{eq:cc4}) is satisfied.\\ \\
We recall that $\delta=\delta(r)$ and $l=l(r)$ and write  \eqref{eq:cumm_smooth} as 
\begin{equation}
\label{eq:cumm_smooth_0}
\left|\hbox{Cum}_{\mu_\cY}^{(r)}(\mathsf{F}_{N,M}^{(\varepsilon, L)})\right|
\ll (N-M)^{r/2+rq} e^{-\delta\gamma}\, \varepsilon^{-rl}
+(N-M)^{q(r-(d-1))^+-(r/2-1)} \gamma^{r-1}~.
\end{equation}

Choosing $\gamma$ of the form
$$
\gamma=c_r\cdot \log (N-M)
$$
with sufficiently large $c_r>0$, and assuming
\begin{equation}
\label{eq:cc6}
(N-M)^{r/2}L^r\varepsilon^{-rl}=o(e^{\delta \gamma})
\end{equation}
we conclude that
$$
\hbox{Cum}_{\mu_\cY}^{(r)}(\mathsf{F}_{N,M}^{(\varepsilon, L)})\to 0\quad\hbox{ as $N\to\infty$}
$$
for all $r\ge 3$. \\

The choice of the parameters $L$, $\varepsilon$, $M$, $K$ and $\gamma$ satisfying all the conditions mentioned earlier is discussed in section \ref{sec:proof of the first main theorem}.

\subsection{Proof of the CLT for the counting function}	
\label{sec:proof of the first main theorem}
Using the characterization by the cumulants (Proposition \ref{prop:method cumulants}), we first show that the sequences $(\mathsf{F}_{N,M}^{(\varepsilon,L)})_{N\geq 1}$, and hence also the sequence $(\mathsf{F}_N)_{N\geq 1}$, converge in distribution to the normal law Norm$_\sigma$.

\begin{thm}
\label{thm:normal distribution of F_N}
Let $m\geq 2$. For every $\xi \in \mathbb{R}$,
$$
\mu_\mathcal{Y}\left(\{y \in \mathcal{Y}:\mathsf{F}_N(y)<\xi \} \right)~\rightarrow \text{Norm}_\sigma(\xi)~
$$
as $N\rightarrow \infty$, with a variance $\sigma^2<\infty$ given in \eqref{eq: explicite variance}.
\end{thm}

\begin{proof}

By Proposition \ref{prop:method cumulants} and considering that $\mathsf{F}_N$ and $\mathsf{F}_{N,M}^{(\varepsilon,L)}$ have the same limit distribution, it is enough to show that
$$
\text{Cum}_{\mu_{\mathcal{Y}}}^{(r)}\left(\mathsf{F}_{N,M}^{(\varepsilon,L)}\right)\rightarrow 0~ \quad \text{as }~ N\rightarrow \infty
$$
when $r\geq 3$, and
$$
\left\lVert \mathsf{F}_{N,M}^{(\varepsilon,L)} \right\rVert_{L^2_\mathcal{Y}}^2\rightarrow \sigma^2~ \quad \text{as }~N\rightarrow\infty~.
$$
We showed in sections \ref{subsec:final estimate variance} and \ref{subsec:estimating cumulants} that these two conditions hold, provided that the parameters 
$$
\varepsilon=\varepsilon(N),~L=L(N),~M=M(N),~\gamma=\gamma(N),~K=K(N)\leq M(N)
$$
satisfy the conditions we recall here,
\begin{align}
M&=o(N^{1/2})~,\tag{\ref{eq:cc1}}\\
(N-M)^{1/2}\varepsilon &\rightarrow 0~,\tag{\ref{eq:cc2}}\\
(N-M)^{1/2}e^{-\theta M}&\rightarrow 0~,\tag{\ref{eq:cc2}}\\
M&\gg \log L~,\tag{\ref{eq:cc3}}\\
(N-M)&=o\left(L^p\right) ~,\text{ for some }p<d~,\tag{\ref{eq:cc4}}\\
M&\gg_r \gamma~,\tag{\ref{eq:cc5}}\\
(N-M)^{r/2}L^r\varepsilon^{-rl}&=o(e^{\delta \gamma})~,\tag{\ref{eq:cc6}}\\
e^{-\delta K} L^{2}\varepsilon^{-2l} &\to 0~, \tag{\ref{eq:cc7}}\\
(N-M)^{-1} e^{-\delta M}  e^{\xi K} L^{2}\varepsilon^{-2l} &\to 0~,\tag{\ref{eq:cc8}}\\
(N-M)^{-1}K^2 &\to 0\tag{\ref{eq:cc9}},\\
KL^{-\frac{\tau-2}{2}}&\rightarrow 0\quad \hbox{, for some}\quad \tau<d \tag{\ref{eq:condition KL}}\\
\varepsilon \, K &\rightarrow 0\tag{\ref{eq:condition epsilon K}}.
\end{align}
One can easily verify that the following choice of parameters, with $ d\geq 3$,
\begin{align}
M&=(\log N) (\log \log N), \label{eq:choice M}\\
\varepsilon&=(N-M)^{-q_1}, ~\text{ for some }q_1>\frac{1}{2},\label{eq:choice epsilon}\\
L&=(N-M)^{q_2} ~ \text{ for some }q_2>0\text{ large enough to satisfy }\eqref{eq:cc4}\label{eq:choice L},\\
K&=c_1\log (N-M)~\text{ for some }~ c_1>0 \text{ large enough to satisfy }\eqref{eq:cc7},\label{eq:choice K}\\
\gamma&=c_r\log (N-M)~\text{ for some }~ c_r>0 \text{ large enough to satisfy }\eqref{eq:cc6}\label{eq:choice gamma}
\end{align}
satisfy the required conditions.\\ 
Hence, Theorem \ref{thm:normal distribution of F_N} follows from Proposition \ref{prop:method cumulants}.
\end{proof}

Next we relate the function $\mathsf{F}_N $ to the counting function $\mathsf{N}_{T,c}$ and show that $(\mathsf{F}_N)_{N\geq 1} $ has the same limit distribution as $(\mathsf{D}_T)_{T>0}$ defined in \eqref{def D_T} below.\\
For $k\in K$ and $\alpha_k \in $ $S^d$ defined by $k(\alpha_k,1)=(0,\dots,0,1,1)\in S^d$, we consider
\begin{equation}\label{def D_T}
    \mathsf{D}_T(k)\: :=\frac{\mathsf{N}_{T,c}(\alpha_k)-C_{c,d}\cdot T}{T^{1/2}},
\end{equation}
where $C_{c,d}:=\text{vol}(F_{1,c})=\mu_\mathcal{V}(\chi_{1,c})$.\\

We recall the equivalent notations introduced in section \ref{subsec:The counting function as ergodic averages on the light-cone},
\begin{align*}
    &\chi=\chi_{1,c}=\chi_{F_{1,c}}\\
    \qand &\chi_t=\chi_{1,c,t}=\chi_{F_{1,c,t}}, \quad\text{for all } t\geq 0.
\end{align*}
The following estimate will be useful for the analysis of the function $\mathsf{D}_T$.
\begin{lem}
\label{lem:asymptote elementary counting function}
We have, for all $N\geq 1$,
\begin{align}
    &\sum_{t=0}^{N-1}\int_{\mathcal{Y}} \widehat{\chi_t}\circ a_t d\mu_{\mathcal{Y}} =  \sum_{t=0}^{N-1}\emph{vol}(F_{1,c,t})+O(1)~,\label{asymptote Y-average of chi}
    \end{align}
    and, in particular,
    \begin{align}\sum_{t=0}^{N-1}\int_{\mathcal{Y}} \widehat{\chi}\circ a_t d\mu_{\mathcal{Y}} = N\cdot C_{c,d}+O(1)~.\label{asymptote Y-average of chi_t}
\end{align}

\end{lem}

\begin{proof}
    By the mean value identity in \eqref{mean value identity}, we have, for all $t\geq 0$,
    $$\text{vol}(F_{1,c,t})=\int_\mathcal{V} \chi_t d\mu_\mathcal{V}= \int_\mathcal{X} \widehat \chi_t d \mu_\mathcal{X}=\mu_{\mathcal{X}}(\widehat{\chi_t})$$
    It follows
\begin{align}
    \left| \sum_{t=0}^{N-1}\int_{\mathcal{Y}} \widehat{\chi_t}\circ a_t d\mu_{\mathcal{Y}} - \sum_{t=0}^{N-1}\text{vol}(F_{1,c,t})\right| 
    &\leq \sum_{t=0}^{N-1}\int_{\mathcal{Y}}\left|  \widehat{\chi_t}\circ a_t- \mu_\mathcal{X}(\widehat \chi_t) \right|d\mu_{\mathcal{Y}}.\label{sum of averages}
\end{align}
We introduce a parameter $L_t>0$ such that $L_t \xrightarrow[ ]{t \rightarrow \infty } \infty$ and use the estimates for the truncated Siegel transform from Proposition \ref{bounds truncated siegel transform}. We note that the implicit constants involving the test function $chi_t$ only depend on its support, and are therefore uniform in $t\geq0$. We have for any $2\leq \tau <d$ and $t\geq \kappa \log L_t$,
\begin{align}
    &\left\lVert \left(\widehat{\chi_t}\circ a_t - \mu_\mathcal{X}(\widehat \chi_t )\right) -\left(\widehat{\chi_t}^{(L_t)}\circ a_t - \mu_\mathcal{X}(\widehat \chi_t^{(L_t)} )\right) \right\rVert_{L^1(\mathcal{Y})}\nonumber\\
    &\qquad\qquad\leq \left\lVert \widehat{\chi_t}\circ a_t  -\widehat{\chi_t}^{(L_t)}\circ a_t\right\rVert_{L^1(\mathcal{Y})} + \mu_{\mathcal{X}}\left( \left| \widehat \chi_t - \widehat \chi_t ^{(L_t)}\right|\right)\nonumber\\
    &\qquad\qquad\ll L_t^{-\frac{\tau}{2}}+ L_t^{-(\tau-1)}\nonumber\\
    &\qquad\qquad\ll L_t^{-\frac{\tau}{2}}.\label{estimate truncation Y}
\end{align}
Introducing further a parameter $\varepsilon_t>0$ such that $\varepsilon_t \xrightarrow[ ]{t \rightarrow \infty } 0$ and using the estimates for the smooth approximation of $\chi_t$ from Proposition \ref{smooth approximation} and from \eqref{epsilon approximation of chi_t} we have
\begin{align}
    &\left\lVert \left(\widehat{\chi_t}^{(L_t)}\circ a_t - \mu_\mathcal{X}(\widehat \chi_t ^{(L_t)})\right) -\left(\widehat f_{t,\varepsilon_t}^{(L_t)}\circ a_t - \mu_\mathcal{X}(\widehat f_{t,\varepsilon_t}^{(L_t)} )\right) \right\rVert_{L^1(\mathcal{Y})} \nonumber\\
    & \qquad \qquad \qquad\leq \left\lVert \widehat{\chi_t}^{(L_t)}\circ a_t - \widehat f_{t,\varepsilon_t}^{(L_t)}\circ a_t\right\rVert_{L^1(\mathcal{Y})} + \mu_{\mathcal{X}}\left( \left| \widehat \chi_t ^{(L_t)}- \widehat f_{t,\varepsilon_t}^{(L_t)}\right|\right)\nonumber\\
    &\qquad \qquad \qquad\ll \varepsilon_t +e^{-\theta t}.\label{estimate smoothening}
\end{align}
Using further the effective equidistribution estimate from Proposition \ref{cor:multiple equidistribution K orbits}, and noting that the implicit constants are uniform in $t>0$, we have
\begin{align}
    \left\lVert \widehat{f_{t,\varepsilon_t}^{(L_t)}}\circ a_t - \mu_\mathcal{X}(\widehat{f_{t,\varepsilon_t}^{(L_t)}} ) \right\rVert_{L^1(\mathcal{Y})} &\ll e^{-\delta t} \left\lVert  \widehat{f_{t,\varepsilon_t}^{(L_t)}}\right\rVert_{C^l} \nonumber\\
    &\ll e^{-\delta t} \varepsilon_t^{-l}L_t.\label{estimate equidistribution 2}
\end{align}
We choose $L_t=t^a$ and $\varepsilon_t=t^{-b}$ for some $a>\frac{2}{\tau}$ and $b>1$, then fix an integer $N_0=N_0(\kappa,a)\geq 1$ such that $t\geq \kappa \log L_t$ for all $t\geq N_0$. \\
Combining \eqref{estimate truncation Y}, \eqref{estimate smoothening}, \eqref{estimate equidistribution 2} with \eqref{sum of averages}, we obtain
\begin{align*}
    \left| \sum_{t=0}^{N-1}\int_{\mathcal{Y}} \widehat{\chi_t}\circ a_t d\mu_{\mathcal{Y}} - N\text{vol}(F_{1,c,t})\right| &\ll \left| \sum_{t=0}^{N_0-1}\left(\int_{\mathcal{Y}} \widehat{\chi_t}\circ a_t d\mu_{\mathcal{Y}} - \text{vol}(F_{1,c,t})\right)\right|\\
    &\qquad + \sum_{t=N_0}^{N-1} \left( L_t^{-\frac{\tau}{2}}+\varepsilon_t+ e^{-\theta t} + e^{-\delta t} \varepsilon_t^{-l}L_t\right)\\
    &= O(1)~.
\end{align*}
\end{proof}

We will also need the following estimate related to the approximation in (\ref{resandwiching N_T,c}).

\begin{lem}
\label{lem:asymptote sandwishing}
We have, for all positive integers $N\geq N_0$ large enough,
$$\int_K \left| \mathsf{N}_{N,c}(\alpha_k)-\sum_{t=0}^{N-1} \widehat{\chi}_{1,c}\circ a_t(k\Lambda_0)  \right| d\mu_K(k) = O(\log N)~.
$$
\end{lem}
\begin{proof}
We have, by (\ref{resandwiching N_T,c}), for all positive integers $N\geq N_0$ large enough, $k\in K$ and $\alpha_k \in S^d$ as in \eqref{eq:definition alpha_k},
\begin{equation}
    \sum_{t=0}^{N-r_0 -1} \widehat{\chi}_{1,c,t}(a_t k\Lambda_0) +O\left( 1\right) ~\leq ~\mathsf{N}_{N,c}(\alpha_k) +O(1) ~\leq~ \sum_{t=0}^{N+r_0 -1} \widehat{\chi}_{1,c}(a_t k\Lambda_0).\label{inequalities a}
\end{equation}
We have moreover, 
\begin{equation}\label{inequalities b}
    \sum_{t=0}^{ N-r_0 -1} \widehat{\chi}_{1,c,t}\circ a_t(k\Lambda_0)\leq \sum_{t=0}^{N-1} \widehat{\chi}_{1,c}\circ a_t(k\Lambda_0) \leq \sum_{t=0}^{ N+r_0-1 } \widehat{\chi}_{1,c}\circ a_t(k\Lambda_0).
\end{equation}
Subtracting \eqref{inequalities a} and \eqref{inequalities b} and averaging over $K$ gives
\begin{align}
&\int_K \left| \mathsf{N}_{N,c}(\alpha_k)-\sum_{t=0}^{N-1} \widehat{\chi}_{1,c}\circ a_t(k\Lambda_0)  \right| d\mu_K(k)\nonumber\\ 
&\leq \int_K \left| \sum_{t=0}^{N+r_0 -1} \widehat{\chi}_{1,c}\circ a_t(k\Lambda_0)-\sum_{t=0}^{ N-r_0 -1} \widehat{\chi}_{1,c,t}\circ a_t(k\Lambda_0)  \right|d\mu_K(k)+O(1)\nonumber\\
&\leq \sum_{t=0}^{ N-r_0 -1} \int_\mathcal{Y}\widehat{\chi}_{_{F_{1,c}\setminus F_{1,c,t}}}\circ a_td\mu_\mathcal{Y}+ \sum_{t= N-r_0 }^{ N+r_0 -1} \int_\mathcal{Y}\widehat{\chi}_{1,c}\circ a_td\mu_\mathcal{Y} +O(1).\label{average discrepancy}
\end{align}
By Lemma \ref{lem:asymptote elementary counting function} we have
\begin{equation}
    \sum_{t= N-r_0}^{N+r_0-1} \int_\mathcal{Y}\widehat{\chi}_{1,c}\circ a_td\mu_\mathcal{Y}=O(1).\label{second sum}
\end{equation}
We have further, by linearity of the Siegel transform, then applying again lemma \ref{lem:asymptote elementary counting function} and the estimate \eqref{volum F_1,c approximation},
\begin{align}
    \sum_{t=0}^{ N-r_0 -1} \int_\mathcal{Y}\widehat{\chi}_{_{F_{1,c}\setminus F_{1,c,t}}}\circ a_td\mu_\mathcal{Y} &= \sum_{t=0}^{ N-r_0 -1} \int_\mathcal{Y}\left(\widehat{\chi}_{_{F_{1,c}}}\circ a_t - \widehat{\chi}_{_{F_{1,c,t}}}\circ a_t\right) d\mu_\mathcal{Y}\nonumber\\
    &= \sum_{t=0}^{ N-r_0 -1} \Big(\text{vol}(F_{1,c}) - \text{vol}(F_{1,c,t})\Big) + O(1)\nonumber\\
    &= O\left(\sum_{t=0}^{ N-r_0 -1} t^{-1}\right) = O\left(\log N\right).\label{average of difference}
\end{align}
Combining \eqref{second sum} and \eqref{average of difference} with \eqref{average discrepancy} yields the claim.
\end{proof}

It follows from Lemmas \ref{lem:asymptote elementary counting function} and \ref{lem:asymptote sandwishing} that
\begin{align*}
&\int_K\left| \mathsf{D}_{N}(k) -\mathsf{F}_N(k\Lambda_0) \right|d\mu_K(k) \\&= \frac{1}{N^{1/2}} \int_K\left| \mathsf{N}_{N,c}(\alpha_k)-\sum_{t=0}^{N-1} \widehat{\chi}\circ a_t(k\Lambda_0)  +\sum_{t=0}^{N-1} \int_\mathcal{Y} \widehat{\chi}\circ a_t -C_{c,d}N \right|d\mu_K(k) \\
&=o(1)~,
\end{align*}
hence, the sequences $(\mathsf{D}_{N})$ and $(\mathsf{F}_N)$ have the same limit distribution, i.e. for all $\xi \in \mathbb{R}$, we have
$$
| \{k \in K:\mathsf{D}_{N}(k)<\xi \}|\rightarrow \text{Norm}_\sigma (\xi)~, \quad\text{ as } N\rightarrow \infty. 
$$ 
If we take $N_T=\left\lfloor T \right\rfloor$, then $N_T \leq T< N_T+1$, hence
\begin{align*}
    &\mathsf{D}_T(k)= \frac{\mathsf{N}_{T,c}(\alpha_k)-C_{c,d}T}{T^{1/2}}\leq \frac{\mathsf{N}_{N_T+1,c}(\alpha_k)-C_{c,d}N_T}{T^{1/2}}=a_T\mathsf{D}_{N_T+1}+b_T~, \\
    \quad \text{with } \qquad &a_T\rightarrow 1 \text{ and } b_T\rightarrow 0.
    \end{align*}
It follows
$$
| \{k \in K:\mathsf{D}_{T}(k)<\xi \}| \geq | \{k \in K:\mathsf{D}_{N_T+1}(k)<(\xi-b_T)/a_T \}|.
$$
Therefore, for any $\varepsilon>0$ and sufficiently large $T$,
$$
| \{k \in K:\mathsf{D}_{T}(k)<\xi \}| \geq | \{k \in K:\mathsf{D}_{N_T+1}(k)<\xi-\varepsilon \}|,
$$
thus
$$
\liminf_{T\rightarrow\infty}| \{k \in K:\mathsf{D}_{T}(k)<\xi \}| \geq \text{Norm}_\sigma(\xi-\varepsilon)~, 
$$
for all $\varepsilon>0$, which implies
$$
\liminf_{T\rightarrow\infty}| \{k \in K:\mathsf{D}_{T}(k)<\xi \}| \geq \text{Norm}_\sigma(\xi)~.
$$
One shows similarly 
$$
\limsup_{T\rightarrow\infty}| \{k \in K:\mathsf{D}_{T}(k)<\xi \}| \leq \text{Norm}_\sigma(\xi)~,
$$
which finishes the proof of Theorem \ref{thm: second result}.

\newpage

\printbibliography

\end{document}